# Tamed to compatible: Symplectic forms via moduli space integration

Clifford Henry Taubes[†]

Department of Mathematics
Harvard University
Cambridge, MA 02138

chtaubes@math.harvard.edu

ABSTRACT: Fix a compact 4-dimensional manifold with self-dual 2$^{nd}$ Betti number one and with a given symplectic form. This article proves the following: The Frêchet space of tamed almost complex structures as defined by the given symplectic form has an open and dense subset whose complex structures are compatible with respect to a symplectic form that is cohomologous to the given one. The theorem is proved by constructing the new symplectic form by integrating over a space of currents that are defined by pseudo-holomorphic curves.

[†]Supported in part by the National Science Foundation and the Clay Math Institute

# 1. Introduction

Suppose that X is a compact, oriented manifold and suppose that ω is a symplectic form on X that is compatible with the orientation. An endomorphism, J, of TX is said to be an almost complex structure when $J^2 = -1$. Such an almost complex structure is said to be *tamed* by ω when the bilinear form ω(·, J(·)) is positive definite. The almost complex structure J is said to be *compatible* with ω when this same bilinear form is also symmetric. Gromov [Gr] observed that tamed almost complex structures and also compatible almost complex structures always exist. Simon Donaldson [D1] posed the following question: If an almost complex structure is tamed by ω, must it be compatible with a symplectic form?

Of particular interest is the case when X is 4-dimensional. The theorem below says something about this question in the case when X is 4-dimensional and its self-dual, second Betti number is 1. By way of a reminder, this Betti number, $b^{2+}$, is the number of positive eigenvalues of the intersection pairing on $H_2(X; \mathbb{R})$. The theorem speaks of a *generic* compatible, almost complex structure. This condition is met by all compatible complex structures from a certain dense and open set. (These terms are defined using the $C^\infty$-Frechêt space topology.) The requirements for membership in this set are given in the body of this article.

**Theorem 1**: *Suppose that X is a compact, oriented 4-dimensional manifold with $b^{2+} = 1$ and with a given symplectic form, ω. A generic, ω-tamed almost complex structures on X is compatible with a symplectic form on X. Moreover, the class in $H^2(X; \mathbb{R})$ of this form can be taken to be that of ω if the latter's class comes from $H^2(X; \mathbb{Q})$.*

Note that in the case when $X = \mathbb{CP}^2$, the conclusions of Theorem 1 hold for every tamed almost complex structure. This follows from what Gromov says in [Gr] about pseudoholomorphic subvarieties and what the author says in [T1] about symplectic forms on $\mathbb{CP}^2$.

What follows is meant to indicate how Theorem 1 is proved. To start, remark that a 2-dimensional current is, by definition, a bounded linear functional on the space of smooth 2-forms on X. All currents here are understood to be 2-dimensional. A current is said to be closed when it is zero on the exterior derivative of any 1-form.

Let J denote an almost complex structure on TX. The complexification of the cotangent bundle of X decomposes as $TX_\mathbb{C}{}^* = T^{1,0}X \oplus T^{0,1}X$ where $T^{1,0}X$ annihilates the subspace in $TX_\mathbb{C}$ where J acts as -i. The bundle of $\mathbb{C}$-valued 2-forms on X decomposes analogously as $T^{2,0}X \oplus T^{1,1}X \oplus T^{0,2}X$ where $T^{2,0}X = \wedge^2 T^{1,0}$, and $T^{0,2}X = \wedge^2 T^{0,1}$. Meanwhile, $T^{1,1}X = T^{1,0}X \otimes T^{0,1}X$. A current is said to be of *type 1-1* when it annihilates $T^{2,0}X$ and $T^{0,2}X$. A current of type 1-1 is said to be *non-negative* when its value is non-negative on any section of $\wedge^2 T^*X$ that can be written as $if\, \sigma \wedge \bar\sigma$ with σ a section of $T^{1,0}X$ and $f \geq 0$ a non-negative function. A current is *positive* when its value on such a form is positive when $f$ is not identically zero.

Suppose that $\Omega$ is a symplectic form on X such that $\Omega \wedge \Omega$ defines the given orientation. If J is $\Omega$-compatible, then the linear functional

$$\upsilon \to \int_X \upsilon \wedge \Omega$$

(1.1)

is a closed, positive current of type 1-1. Note however that a closed, positive current of type 1-1 need not be given by integration against a symplectic form. An example of a non-negative current of this sort is given next.

A closed set $C \subset X$ with finite, non-zero 2-dimensional Hausdorff measure is said to be a J-holomorphic subvariety if it has no isolated points, and if the complement of a finite set of points in C is a smooth submanifold with J-invariant tangent space. Such a subvariety has a canonical orientation given by J. Integration on the smooth part of C is defined using this orientation. This understood, the linear functional

$$\upsilon \to \int_C \upsilon$$

(1.2)

defines a closed, non-negative current of type 1-1.

The set of closed, non-negative, type 1-1 currents is convex. As a consequence, a non-negatively weighted average of currents that are defined as in (1.2) by J-holomorphic subvarieties also defines a closed, non-negative, type 1-1 current. When J is suitably generic, and tamed by a given symplectic form, such an average is described momentarily whose current is very nearly given by (1.1) with $\Omega$ a compatible symplectic form. A final step in the proof makes the needed modifications to obtain the desired symplectic form.

Section 1a to follow describes in some detail the current that is used to prove Theorem 1. Section 1b proves Theorem 1 given what is said in Section 1a. The subsequent sections of this article are organized as follows: Section 4 has the proofs of the three key propositions in Section 1a. Section 2 supplies what is needed in Section 4 with regards to the properties of pseudoholomorphic curves, and Section 3 supplies what is needed about their moduli spaces. A lengthy appendix to this article supplies the proofs of the various genericity assertions that are made in Section 3.

The motivation for the use of J-holomorphic subvarieties to say something about compatible symplectic forms comes from a construction in Donaldson's seminal paper [D2] on the existence of symplectic subvarieties in a compact manifold with a given symplectic form. In particular, Donaldson used such subvarieties to construct a sequence of currents that converges to the current defined by integration against the given symplectic form. The construction introduced by Donaldson is not used here to find the needed subvarieties. The desired J-holomorphic subvarieties are found by using versions of the Seiberg-Witten equations by using the main theorem in [T2] and the wall crossing formulae in [LL1]. The application of these equations requires the $b^{2+} = 1$ condition for X. What is done here is reminiscent of what is done by Gromov for the case when X = $\mathbb{CP}^2$ [G].

Donaldson suggests in [D1] a very different approach to his question, one along the lines used by Yau in [Ya] to prove the Calabi conjecture. This alternate approach is considered by Weinkove in [W], with Tosatti and Yau in [TWY] and with Tosatti in [TW]. Tosatti's Harvard Ph. D. thesis says more about this approach as well. Meanwhile, Li and Zhang in[LV] discuss cohomological issues that concern this question of tamed versus compatible almost complex structures.

Before reading the detailed arguments, the reader should take note of the following conventions: First, no generality is lost by assuming that ω defines a class in $H^2(X; \mathbb{R})$ that comes from $H^2(X; \mathbb{Q})$. This is because the taming condition is an open condition on the set of almost complex structures. This rationality condition is imposed below. The second convention concerns notation: In all cases, $c_0$ denotes a real number that is greater than 1, and also independent of any parameters that are germaine to the discussion at hand. Its value can be assumed to increase between subsequent appearances.

The author owes debts of gratitude to T-J Li for pointing out some key features of the $b^{2+} = 1$ Seiberg-Witten invariants, and to Rosa Sena-Dias for asking probing questions. The author also thanks the Clay Math Institute for their generosity, and the Mathematical Sciences Research Institute in Berkeley for their hospitality.

**a) The current**

The subsection describes in more detail a non-negative 1-1 current that is used to obtain a symplectic form which is compatible with an almost complex structure chosen from a certain generic subset of those that are ω-tamed. The description is given in five parts.

*Part 1*: This part sets the stage. To begin introduce [ω] to denote the cohomology class of the symplectic form. Fix a large integer N, and in particular one such that e = N[ω] is a class from $H^2(X; \mathbb{Z})$. Use c ∈ $H^2(X; \mathbb{Z})$ to denote the first Chern class of $T^{1,0}X$. When ê and ô are classes in $H^2(X; \mathbb{Z})$, use ê·ô to denote their cup product pairing. Granted this notation, introduce

$$\iota_e = e \cdot e - c \cdot e \quad \textit{and} \quad g = \tfrac{1}{2}(e \cdot e + c \cdot e) + 1.$$

(1.3)

It is assumed in what follows that N is chosen so that $\iota_e$ is very large. Keep in mind for what follows that $\iota_e$ is in all cases an even number.

View the set of ω-tamed almost complex structures as a topological space whose topology comes from the Frêchet topology on $C^\infty(X; \text{End}(TX))$. The almost complex structure is taken from a certain dense, open subset of ω-tamed almost complex structures. This subset is denoted here by $\mathcal{J}_{e*}$; it is characterized in the subsequent sections. Let J denote the chosen almost complex structure. The polarization of the quadratic form given by $\omega(\cdot, J(\cdot))$ defines a Riemannian metric on X, and this metric is used implicitly in all that follows to define, for example, norms, distances on X, integration over open sets in X and self-dual and anti-self dual 2-forms.

*Part 2*: Use J to define the set $\mathcal{M}_{e,g}$ whose elements are irreducible, J-holomorphic subvarieties with the following two properties: First, each has fundamental class Poincaré dual to e. Second, each is the image via an almost everywhere 1-1, pseudoholomorphic map to X of a connected, complex curve with genus g. The definition of $\mathcal{J}_{e*}$ guarantees that $\mathcal{M}_{e,g}$ is a smooth manifold whose dimension is $\iota_e$.

Set $d = \tfrac{1}{2}\iota_e$ and use $\mathcal{M}_{e,g,d}$ to denote the subset in $\mathcal{M}_{e,g} \times (\times_d X)$ that consists of elements of the form $(C, x_1, \ldots, x_d)$ with each $x_k \in C$. The definition of $\mathcal{J}_{e*}$ guarantees that $\mathcal{M}_{e,g,d}$ is smooth submanifold in $\mathcal{M}_{e,g} \times (\times_d X)$ whose dimension is $4d = 2\iota_e$.

*Part 3*: The $b^{2+} = 1$ condition on X enters in this part. Introduce $\pi_d: \mathcal{M}_{e,g,d} \to \times_d X$ to denote the restriction of the projection map onto the $\times_d X$ factor of $\mathcal{M}_{e,g} \times (\times_d X)$. The $b^{2+}$ condition guarantees the following:

**Proposition 1.1**: *Take N very large so as to define* $e = N[\omega]$, *choose J from the set* $\mathcal{J}_{e*}$, *and use J to define* $\mathcal{M}_{e,g,d}$. *Then the map* $\pi_d: \mathcal{M}_{e,g,d} \to \times_d X$ *is onto the complement of a compact, measure zero subset. Moreover,* $\pi_d^{-1}(\eta)$ *is non-empty and finite if* $\eta \in \times_d X$ *is a regular value of* $\pi_d$ *in the complement of this subset.*

This proposition is proved in Section 4a.

A set $\mathcal{K} \subset \mathcal{M}_{e,g,d}$ is said in what follows to be *fiberwise finite* when the following condition holds: There exists $K \geq 1$ such that $\pi_d^{-1}(\eta) \cap \mathcal{K}$ has at most K elements when $\eta$ comes from the complement of the set of critical values of $\pi_d$.

*Part 4*: Fix a fiberwise finite subset $\mathcal{K} \subset \mathcal{M}_{e,g,d}$. Fix $K \geq 1$ and a point $\eta \in \pi_d(\mathcal{K})$ such that $\pi_d^{-1}(\eta) \cap \mathcal{K}$ has at most K elements. Let $\upsilon$ denote a 2-form on X and define

$$\phi_\eta(\upsilon) = \sum_{C \in \pi_d^{-1}(\eta) \cap \mathcal{K}} \int_C \upsilon \ .$$

(1.4)

With $\eta$ so chosen, the assignment $\upsilon \to \phi_\eta(\upsilon)$ defines a current, which is to say a bounded, linear functional on $C^\infty(X; \wedge^2 T^*M)$. The bound on the norm comes courtesy of the taming assumption because the latter leads to a bound on the area of any given subvariety in $\mathcal{M}_{e,g}$; and such a uniform area bound leads directly to a bound $|\phi_\eta(\upsilon)| \leq c_0 K \sup_X |\upsilon|$. .

With $\upsilon$ fixed, the assignment $\eta \to \phi_\eta(\upsilon)$ defines a function on $\pi_d(\mathcal{K})$. The following proposition says something about this function for certain choices of $\mathcal{K}$.

**Proposition 1.2**: *Fix $N \gg 1$ to define $e = N[\omega]$, and take the almost complex structure from $\mathcal{J}_{e*}$. Suppose that $\mathcal{K} \subset \mathcal{M}_{e,g,d}$ is a non-empty, open, fiberwise finite set. Then*
- $\pi_d(\mathcal{K}) \subset \times_d X$ *is a measurable subset.*
- *Let $\upsilon$ denote a given smooth 2-form on X. Then the function $\eta \to \phi_\eta(\upsilon)$ on $\pi_d(\mathcal{K})$ is a bounded, measurable function with bound a $\mathcal{K}$-dependent multiple of $\sup_X |\upsilon|$.*
- *The assignment*

$$\upsilon \to \Phi_\mathcal{K}(\upsilon) = \int_{\eta \in \pi_d(\mathcal{K})} \phi_\eta(\upsilon)$$

*defines a closed and non-negative type 1-1 current on X which is non-trivial if $\pi_d(\mathcal{K})$ has positive measure.*

This proposition is proved in Section 4b.

*Part 5*: The next proposition says more about certain versions of $\Phi_\mathcal{K}$. In particular, it asserts that $\mathcal{K}$ can be chosen so as to make $\Phi_\mathcal{K}$ a positive current which is very nearly given by integration as in (1.1) with a form $\Omega$ that is non-degenerate and bounded. This proposition is at the heart of the proof of Theorem 1.

**Proposition 1.3**: *Fix $N \gg 1$ so as to define $e = N[\omega]$. The set $\mathcal{J}_{e*}$ is such that if $J \in \mathcal{J}_{e*}$, then there exists a fiberwise finite set $\mathcal{K} \subset \mathcal{M}_{e,g,d}$ and a constant $\kappa > 1$ with the following two properties:*
- *The conclusions of Proposition 1.2 hold and so the current $\Phi_\mathcal{K}$ is well defined.*
- *Fix $t > 0$ but small, and let $B \subset X$ denote a ball of radius $t$. Let $\sigma$ denote a unit length section of $T^{1,0}X|_B$ and let $f_B$ denote the characteristic function of B. Then*

$$\kappa^{-1} t^4 < \Phi_\mathcal{K}(i f_B \sigma \wedge \bar\sigma) < \kappa t^4.$$

The proof of this proposition occupies Sections 4d-f of this article.

What follows directly says something about Proposition 1.3's residual set constraint on the almost complex structure. The proof of the proposition uses this constraint to insure that certain pseudoholomorphic curve moduli spaces have 'generic' properties. To give an idea of what can go wrong, suppose that all subvarieties that contribute to $\Phi_\mathcal{K}$ contain a particular point in X. If such is the case, then the upper bound in the proposition must be replaced by $\kappa t^2$. The resulting current will not define a smooth form on X, although it may well define one on the blow up of X at the given point. This example is the extreme case of what might be called the key concern: The bounds claimed by the proposition require that the relevant subvarieties have uniform 'density' across X: Regions where the density is low determine the lower bound; whereas high density regions, those where the subvarieties concentrate, determine the upper bound. The needed genericity properties of the relevant moduli spaces are described by the various propositions in the upcoming Section 3 and their proofs are given in the appendix to this article. The proofs are lengthy, but at the heart of each is the Sard-Smale theorem [Sm]; this used to deduce the existence of a residual set of regular values for a relevant map.

**b) The proof of Theorem 1**

The proof of this theorem has four parts.

*Part 1*: Fix a version of the current $\Phi_\mathcal{K}$ as given by Proposition 1.3. This first part of the proof defines a smoothing of $\Phi_\mathcal{K}$ so as to give a current that is defined as in (1.1) by integration against a smooth 2-form. The construction that follows is along the lines of one that appears in Section 2 of [Su]. To start, fix an exponential map exp: $TX \to X$ whose differential along the zero section is the identity. There exists $\delta > 0$ such that exp embeds the radius $\delta$ ball about the origin in any given fiber of TX. Use $BX \subset TX$ in what follows to denote the space of vectors with norm less than $\delta$. Let $\pi$: $TX \to X$ denote the bundle projection. If $x \in X$ has been specified, then $TX|_x$ and $BX|_x$ are used to denote the respective fibers at x of TX and BX.

Use $\exp_!$ to denote the push-forward homomorphism given by the exponential map; this a homomorphism from the space of smooth, compactly supported sections of $\wedge^p T^*(BX)$ to $C^\infty(X; \wedge^{p-4} T^*X)$.

As $T(TX|_x)$ is canonically isomorphic to $TX|_x$, the fiber tangent space $T(TX|_x)$ inherits a translationally invariant metric, this the metric on TX at x. Let $w_x$ denote the associated translationally invariant volume 4-form on the vector space $TX|_x$. Fix a smooth, non-increasing function, $\chi$: $[0, \infty) \to [0, 1]$ that is equal to 1 on $[0, \frac{5}{16}]$ and is equal to 0 on $[\frac{7}{16}, \infty)$. Use A to denote

$$A = 2\pi^2 \int_0^\infty \chi(t) t^3 dt .$$

(1.5)

Given $\varepsilon > 0$ but much smaller than both $\delta$ and the constant $\rho$ used to define $\Phi$, define the function $\chi_\varepsilon$ on TX by setting $\chi_\varepsilon(v) = \chi(\varepsilon^{-1}|v|)$. A standard construction (see, e.g. Chapter 1.6 of [BT]) produces a compactly supported, closed 4-form on BX whose pull-back to any given fiber $BX|_x$ is the 4-form $A^{-1}\varepsilon^{-4}\chi_\varepsilon w_x$. This form represents the Thom class in the compactly supported cohomology of BX. Denote this 4-form by $w_\varepsilon$. Note for future reference that $w_\varepsilon$ can be constructed so that its derivatives obey $|\nabla w_\varepsilon| \leq c_0 \varepsilon^{-5}$.

The form $w_\varepsilon$ is used in what follows to define the infinitely smoothing operator from $C^\infty(X; \wedge^*T^*X)$ to itself that acts as

$$\upsilon \to \exp_!(w_\varepsilon \wedge \pi^*\upsilon) .$$

(1.6)

This smoothing map is denoted by $\mathcal{T}_\varepsilon$. The map $\mathcal{T}_\varepsilon$ approximates the identity operator in the sense that

$$|\upsilon - \mathcal{T}_\varepsilon(\upsilon)| \leq c_0 \varepsilon |\nabla \upsilon|.$$

(1.7)

Moreover, if $r > 0$ and $\upsilon$ has support in a ball of radius r, then $\mathcal{T}_\varepsilon(\upsilon)$ has support in the concentric, radius $r+\varepsilon$ ball. It is also the case that $\mathcal{T}_\varepsilon d = d\mathcal{T}_\varepsilon$; this because push-forward of compactly supported forms and pull-back both commute with exterior differentiation.

This smoothing operator extends to define a bounded, linear map from distribution valued forms to $C^\infty(X; \wedge^*T^*X)$. For example, fix $x \in X$ and an element, $\beta_x$, in $\wedge^2T^*X$. These define the distribution that acts to send a given 2-form $\upsilon$ to $*(\upsilon|_x \wedge \beta)$ where $*$ here denotes the metric's Hodge star. Denote this distribution valued 2-form by $\beta\delta_x$. Then $\mathcal{T}_\varepsilon(\beta\delta_x)$ is a smooth 2-form with support in the radius $\varepsilon$ ball centered at x. The norm of this form obeys

$$|\mathcal{T}_\varepsilon(\beta\delta_x)| \leq c_0 \varepsilon^{-4} |\beta| .$$

(1.8)

This form approximates the distribution $\upsilon \to *(\upsilon|_x \wedge \beta)$ in the sense that the integral of over X of $\upsilon \wedge \mathcal{T}_\varepsilon(\beta\delta_x)$ differs from $*(\upsilon|_x \wedge \beta)$ by no more that $c_0 \varepsilon |\nabla\upsilon||\beta|$.

A 2-form on X is defined by specifying the Hodge dual of its wedge with any given 2-form. This understood, define the 2-form $\Omega^\varepsilon$ by demanding that

$$*(\beta \wedge \Omega^\varepsilon|_x) = \Phi_\mathcal{K}(\mathcal{T}_\varepsilon(\beta\delta_x))$$

(1.9)

for any given $\beta \in \wedge^2 T^*X|_x$. Note that this implies the formula

$$\int_X \upsilon \wedge \Omega^\varepsilon = \Phi_\mathcal{K}(\exp_!(w_\varepsilon \wedge \pi^*\upsilon))$$

(1.10)

for any given 2-form $\upsilon$.

What is written in (1.9) defines a smooth 2-form, this by virtue of the fact that $w_\varepsilon$ is a smooth, compactly supported form on BX. The form $\Omega_\varepsilon$ is also a closed form. To see why, suppose that $\upsilon = d\alpha$. It follows from (1.10) that

$$\int_X d\alpha \wedge \Omega^\varepsilon = \Phi_\mathcal{K}(\exp_!(w_\varepsilon \wedge \pi^*d\alpha))$$

(1.11)

As $dw_\varepsilon = 0$, the right hand side of (1.11) is $\Phi_\mathcal{K}(\exp_!(d(w_\varepsilon \wedge \alpha)))$. Given that push-forward of compactly supported forms commutes with d, this is the same as $\Phi_\mathcal{K}(d(\exp_!(w_\varepsilon \wedge \alpha)))$. The latter expression is zero by virtue of the fact that $\Phi_\mathcal{K}$ is a closed current.

The small $\varepsilon$ versions of $\Omega^\varepsilon$ supply the desired smoothings of $\Phi_\mathcal{K}$. Indeed, the fact that integration against $\Omega^\varepsilon$ approximates $\Phi_\mathcal{K}$ follows from (1.7) and (1.10).

*Part 2*: This part of the proof supplies an upper bound for $|\Omega^\varepsilon|$, an upper bound for the (2, 0) component of $\Omega^\varepsilon$; and a lower bound for the $\Omega^\varepsilon \wedge \Omega^\varepsilon$. As it turns out, the latter is positive when $\varepsilon$ is small; and so the small $\varepsilon$ versions of $\Omega^\varepsilon$ are symplectic.

A bound for $|\Omega^\varepsilon|$ of the form

$$|\Omega^\varepsilon| \leq c_J$$

(1.12)

follows directly from (1.8) and the upper bound given by Proposition 1.4. Here, and in what follows, $c_J \geq 1$ is an $\varepsilon$-independent constant. Note however that $c_J$ depends on J and the other data used to define $\Phi_\mathcal{K}$. As with $c_0$, the value of $c_J$ can be assumed to increase between subsequent appearances.

Consider next the norm of the (2, 0) part of $\Omega^\varepsilon$. To start, fix $x \in X$ and a coordinate chart centered at x with complex, Gaussian coordinates, (z, w), that are defined in a ball about the origin in $\mathbb{C}^2$ and are such that dz and dw are orthonormal at x and span $T^{1,0}X|_x$. Use Taylor's theorem with remainder to see that the norms of both dz and dw on the -i subspace of J at a point x´ in this coordinate chart are bounded by $c_0 \text{dist}(x, x´)$. Granted the latter bound, and granted (1.8), it then follows from the upper bound in Proposition 1.3 that the (2, 0) part of $\Omega^\varepsilon$ obeys

$$|(\Omega^\varepsilon)_{2,0}| \leq c_J \varepsilon$$

(1.13)

when ε is very small. Here, $(\Omega^\varepsilon)_{2,0}$ denotes the $(2, 0)$ part of $\Omega^\varepsilon$. Given the bound on $|\nabla w_\varepsilon|$ by $c_0 \varepsilon^{-5}$, a similar line of reasoning finds

$$|\nabla(\Omega^\varepsilon)_{2,0}| \leq c_J .$$
(1.14)

To consider $\Omega^\varepsilon \wedge \Omega^\varepsilon$, return to the coordinate patch just described. Let $(\alpha, \beta) \in \mathbb{C}^2$ denote a given unit vector, and let $u = \alpha z + \beta w$. What was just said about the norms of dz and dw on -i eigenvectors of J implies that the $(1, 1)$ portion of $du \wedge d\bar{u}$ can be written as

$$(du \wedge d\bar{u})_{1,1} = f\, \sigma \wedge \bar{\sigma}$$
(1.15)

where $\sigma$ is a unit length section of $T^{1,0}X$ and where $|f - 1| \leq c_0\, \mathrm{dist}(x, \cdot)^2$. This last fact with (1.13) and the lower bound from Proposition 1.3 imply that

$$*(i(du \wedge d\bar{u}) \wedge \Omega^\varepsilon) \geq c_J^{-1}$$
(1.16)

when ε is small. The latter bound implies directly that

$$*(\Omega^\varepsilon \wedge \Omega^\varepsilon) \geq c_J^{-1}$$
(1.17)

when ε is small.

*Part 3*: This part sets the stage for a modification of $\Omega^\varepsilon$ that results in a symplectic $(1, 1)$ form. To start the stage setting, let $\wedge^{2+}T^*X \subset \wedge^2 T^*X$ denote subbundle of self-dual 2-forms. Use $d^+: C^\infty(X; T^*X) \to C^\infty(X; \wedge^{2+}T^*X)$ to denote the composition of first d and then orthogonal projection. This map is surjective onto the complement of the 1-dimensional subspace spanned by the closed self-dual forms. Hodge theory gives an inverse to $d^+$ on the $L^2$-orthogonal complement of this line of self-dual closed forms. The following lemma says something about this inverse.

**Lemma 1.4**: *There is an inverse of $d^+$ and a constant $\kappa \geq 1$ with the following properties: Suppose that $\upsilon \in C^\infty(X; \wedge^{2+}T^*X)$ is $L^2$-orthogonal to the kernel of d. Then $|(d^+)^{-1}\upsilon| \leq \kappa|\upsilon|$. Moreover, if $t > 0$, then $|\nabla((d^+)^{-1}\upsilon)| \leq \kappa(|\ln(t)|\,|\upsilon| + t|\nabla \upsilon|)$.*

*Proof of Lemma 1.4*: A standard parametric construction writes $(d^+)^{-1}$ as integration against an integral kernel. Viewed as integration against a form on $X \times X$, the latter is singular along the diagonal. This singular form, $\mathfrak{a}$, obeys $|(\mathfrak{a}|_{(x,y)})| \leq c_0\, \mathrm{dist}(x, y)^{-3}$ and $|(\nabla \mathfrak{a})|_{(x,y)}| \leq c_0\, \mathrm{dist}(x, y)^{-4}$. These bounds lead directly to the lemma's assertions.

One other fact is needed for the promised modification of $\Omega^\varepsilon$: The form $\omega$ is self dual with respect to the metric that is obtained by polarizing the quadratic form $\omega(\cdot, J(\cdot))$. Since $b^{2+} = 1$, and since $\omega$ is closed, it generates the line of closed, self-dual 2-forms.

*Part 4*: This last part of the proof supplies the promised modification of $\Omega^\varepsilon$. To start, let $K_\mathbb{R} \subset \wedge^{2+}T^*X$ denote the underlying real bundle for $T^{2,0}X$. Use $\Omega^\varepsilon_K$ to denote the orthogonal projection of $\Omega^\varepsilon$ onto the subbundle $K_\mathbb{R}$ and use $x_\varepsilon$ to denote the integral over X of $\Omega^\varepsilon_K \wedge \omega$. It follows from (1.13) that this number obeys $|x_\varepsilon| \leq c_J \varepsilon$. Meanwhile, use $\omega_K$ to denote the projection of $\omega$ into $K_\mathbb{R}$. If this is zero, then J is already $\omega$ compatible, so assume that $\omega_K$ is not zero. Use $z_J$ to denote the integral over X of $\omega_K \wedge \omega$, this the $L^2$ norm of $\omega_K$.

The form $\Omega^\varepsilon_K - z_J^{-1} x_\varepsilon \omega_K$ is $L^2$ orthogonal to $\omega$, and thus to all closed self-dual 2-forms. It is thus in the image of $d^+$ and so can be written as $d^+ a_\varepsilon$ with $a_\varepsilon \in C^\infty(X; T^*X)$. Use Lemma 1.4 with (1.13) and (1.14) to see that $a_\varepsilon$ can be chosen so that

$$|\nabla a_\varepsilon| \leq c_J(t + \varepsilon |\ln(t)|) .$$

(1.18)

In particular, take $t = \varepsilon$ to see that the latter choice for $a_\varepsilon$ obeys

$$|\nabla a_\varepsilon| \leq c_J \varepsilon |\ln \varepsilon| .$$

(1.19)

Set $\Omega = \Omega^\varepsilon - z_J^{-1} x_\varepsilon \omega - da_\varepsilon$ with $a_\varepsilon$ as just described. The (2, 0) part of this 2-form is zero; and it follows from (1.17) and (1.19) that the small $\varepsilon$ versions are symplectic. Any such small $\varepsilon$ version of this 2-form is a J compatible symplectic form.

It remains now only to explain why the class of $\Omega$ in $H^2(X; \mathbb{R})$ is a multiple of the symplectic class $[\omega]$. Of course, this is the case if and only if such is the case for the class defined by $\Omega^\varepsilon$. To see about the latter, note that it is sufficient to verify that the linear functional $\upsilon \to \int_X \upsilon \wedge \Omega^\varepsilon$ on $C^\infty(X; \wedge^2 T^*X)$ descends to $H^2(X; \mathbb{R})$ as a multiple of the functional given by cup product with $[\omega]$. This follows from (1.6) and (1.10) from the fact that $w_\varepsilon$ represents the Thom class in the compactly supported cohomology of BXs. See, e.g Chapter 1.6 in [BT].

## 2. Pseudoholomorphic subvarieties

The purpose of this section is to summarize various known results about pseudoholomorphic subvarieties for use in proving the various claims in Sections 1a,b.

### a) The existence of J-holomorphic subvarieties

The cup product pairing between classes e and e´ in $H^2(X; \mathbb{R})$ is denoted subsequently by e·e´. Assume in what follows that $b^{2+} = 1$ and that ω is a given symplectic form on X such that $\omega \wedge \omega$ gives the orientation. A class in $H^2(X; \mathbb{R})$ is said to lie in the *positive cone* when it has positive cup product pairing with itself and with the class that is defined by ω. Introduce $c \in H^2(X; \mathbb{Z})$ to denote the first Chern class of the complex line bundle $T^{2,0}X$. The class c depends only on the component of ω in the subspace of 2-forms with nowhere zero square. This class is used to associate the integer $\iota_e = e \cdot e - c \cdot e$ to a given class $e \in H^2(X; \mathbb{Z})$. Note that this integer $\iota_e$ is even.

Let J denote a given ω-tamed almost complex structure. If $C \subset X$ is a J-holomorphic subvariety, then it defines via (1.2) a linear form on $H^2(X; \mathbb{R})$. This is an integer valued linear functional. Its Poincaré dual in $H^2(X; \mathbb{Z})$ is denote by $e_C$. Since ω is positive on the smooth part of C, the form $e_C$ is a non-zero class in the positive cone.

The proposition below refers to an *irreducible* J-holomorphic subvariety. A subvariety is irreducible if its smooth locus is connected. Any given J-holomorphic subvariety is a union of a finite set of irreducible J-holomorphic subvarieties.

**Proposition 2.1**: *Fix a class $\hat{e} \in H^2(X; \mathbb{Z})$ in the positive cone with $\hat{e} \cdot \hat{e} > 0$. There exists $\kappa > 1 + |c \cdot \hat{e}|/\hat{e} \cdot \hat{e}$ such that if n is an integer greater than κ, then the following is true: Suppose that J is a given, ω-tamed almost complex structure. Fix $\frac{1}{2}\iota_{n\hat{e}}$ points in X and there exists a finite set, Θ, of pairs of the form (C, m) with $C \subset X$ an irreducible, J-holomorphic subvariety and m a positive integer. Moreover,*
- $\sum_{(C,m) \in \Theta} m\, e_C = n\hat{e}$.
- $\cup_{(C,m) \in \Theta} C$ *contains the chosen set of points.*

*Proof of Proposition 2.1*: Suppose first that the chosen points are distinct, and that J is taken from a certain dense and open set of compatible almost complex structures. The assertion in this case follows directly from Lemma 3.3 in [LL2] using the main theorem [LL3]. See also [T2]. The theorems in [LL3] and [T2] equate certain Seiberg-Witten invariants of X with invariants of the symplectic form that are obtained by making a suitably weighted count of pseudoholomorphic curves that go through the chosen points. Lemma 3.3 in [LL2] is used to prove that the desired Seiberg-Witten invariants are non-zero. The statement of the theorem follows when J is not from this set, or when the points are not distinct follows using the standard convergence theorems for sequences of pseudoholomorphic curves. A proof is given in [Wo]. See also [Ye], or the more recent book [H]. The assertion can also be proved by arguments that differ only in notation from those in Section 6a of [T3]; a proof along the latter lines can be found in Section 4 of [T4].

Given what was just said about ω-compatible almost complex structures, the case when J is ω-tamed follows from three observations. Here are the first two observations: The space of ω-tamed almost complex structures is contractible and the space of ω compatible almost complex structures is non-empty. Here is the third observation: The aforementioned weighted curve counts are invariant with respect to variations of the almost complex structure along a path of tamed almost complex structures. Note in this regard that the aforementioned convergence theorems hold for tamed as well as compatible almost complex structures. The point being that these convergence theorems require only an area bound for the sequence in question; and, as noted in the upcoming Lemma 2.2, such a bound is implied by the taming assumption.

**b) Properties of J-holomorphic subvarieties**

Suppose that J is an ω-tamed almost complex structure. The four parts of this subsection summarize some basic facts about irreducible, J-holomorphic subvarieties. In all that follows, C denotes an irreducible, J-holomorphic subvariety.

The assertions made in this subsection are well known to the experts.

*Part 1*: There exists a compact, connected complex curve $C_0$ with an almost everywhere 1-1 map $\varphi: C_0 \to C$ whose differential intertwines the complex structure on $C_0$ with J's restriction to TC. The curve $C_0$ is said to be C's *model curve*. The genus of the model curve of a given irreducible, J-holomorphic subvariety is apriori bounded. In particular, the adjunction formula implies that

$$\mathrm{genus}(C_0) \leq \tfrac{1}{2}(e_C \cdot e_C + c \cdot e_C) + 1$$

(2.1)

with equality if and only if C is embedded. The genus of $C_0$ is one less than the right hand side of (2.1) if and only if $\varphi$ is an immersion with precisely one transversal double point. Note that in general, if C is such that all singularities are transversal double points, then the genus of $C_0$ is given by subtracting the number of such double points from the right hand side of (2.1).

*Part 2*: View $\varphi$ as a map from $C_0$ into X. As such, it has but a finite set of critical points; and these are points where its differential is identically zero. Each point in the complement of the critical locus has a neighborhood that is embedded by $\varphi$. This being the case, there is a rank 2-real vector bundle that is defined on the complement of the critical locus, and whose fiber at a given point p is the normal bundle at $\varphi(p)$ of the image via $\varphi$ of a neighborhood of p that is embedded by $\varphi$. This bundle is called the *normal* bundle. The almost complex structure J gives this bundle the structure of a complex line bundle. The latter is denoted by N. As explained momentarily, the bundle N has a canonical extension over the critical locus.

*Part 3*: Let C denote an irreducible, pseudoholomorphic subvariety and let $C_0$ denote its model curve. Let $p \in C_0$. There is a complex coordinate, u, for a disk in $C_0$ centered on p, and a complex coordinate centered on $\varphi(p)$ in X such that the map $\varphi$ appears as a map from a neighborhood of the origin in $\mathbb{C}$ to a neighborhood of the origin in $\mathbb{C}^2$ having the form

$$u \to (u^{n+1} + \mathfrak{r}_z, c u^{n+k+1} + \mathfrak{r}_w)$$

(2.2)

where $n \geq 0$, $k \geq 1$. Here, $c \in \mathbb{C}$ is zero when $n = 0$; but is non-zero when $n \geq 1$. Meanwhile, $\mathfrak{r}_z$ and $\mathfrak{r}_w$ have the following properties: First, $|\mathfrak{r}_z| \leq c_0 |u|^{n+2}$ and $|d\mathfrak{r}_z| \leq c_0 |u|^{n+1}$. Second, $|\mathfrak{r}_w| \leq c_0 |u|^{n+k+2}$ and $|d\mathfrak{r}_w| \leq c_0 |u|^{n+k+2}$ if $n \neq 0$. If $n = 0$, then $|\mathfrak{r}_w| \leq c_0 |u|^2$ and $|\mathfrak{r}_w| \leq c_0 |u|$. See Proposition 2.6 in [McD]. In any event, (2.2) is proved as a part of the proof of Lemma A.7 in the Appendix.

The point p is a critical point of the map $\varphi$ if and only if $n \geq 1$. It is a consequence of (2.2) that the bundle N extends over p as the pull-back from $\mathbb{C}^2$ via (2.2) of the span of the vector field $\frac{\partial}{\partial w}$.

With regards to complex coordinates, the term *adapted* coordinates is used below to refer to a special sort of coordinate chart. To describe such a chart, fix a point $x \in X$. An adapted a coordinate chart centered at x denotes complex coordinates, (z, w) defined on a radius $c_0^{-1}$ ball centered at x with both vanishing at x, with dz and dw orthonormal at x, with $\{dz, dw\}$ spanning $T^{1,0}X$ at x, and with the norms of $|\nabla dz|$ and $|\nabla dw|$ bounded on the coordinate domain by $c_0$.

*Part 4*: This part of the subsection gives upper and lower bounds for the area of C's intersection with a ball about any given point in C. These bounds are summarized by the next lemma. This lemma uses $[\omega]$ to denote the class in $H^2(X; \mathbb{R})$ defined by $\omega$.

**Lemma 2.2**: *Let J denote an $\omega$-tamed almost complex structure. Fix a metric on X that makes J orthogonal. There exists $\kappa \geq 1$ with the following significance: Let $C \subset X$ denote an irreducible, J-holomorphic subvariety. Then*
- *The area of C is greater than $\kappa^{-1} e_C \cdot [\omega]$ and less than $\kappa e_C \cdot [\omega]$*
- *Fix $r > 0$ and a point $x \in C$. Let $\mathfrak{a}_x(r)$ denote the area of C's intersection with the ball of radius r in X centered at x. Then $\kappa^{-1} r^2 < \mathfrak{a}_x(r) \leq (e_C \cdot [\omega]) \kappa r^2$.*

***Proof of Lemma 2.2***: Use the induced metric from X to define a metric on the tangent bundle of C over the complement of its singular points. Let x denote a point in this smooth locus, and let $v_1$ denote a unit length vector in $TC|_x$. The pair $(v_1, Jv_1)$ is thus an orthonormal basis for $TC|_x$. Then the taming condition implies that $c_0^{-1} \leq \omega(v_1, Jv_1)$, so integration over C of $\omega$ gives an upper bound to the area of C. Thus, area(C) $\leq c_0\, e_C \cdot [\omega]$. The lower bound on the area of C follows from the fact that $|\omega(v_1, Jv_1)| \leq c_0$. The local area bounds follow from a *monotonicity* inequality that asserts the following: If $x \in C$, and if $r_1 \geq r_0 > 0$, then

$$\mathfrak{a}_x(r_1) \geq c_0^{-1}\, \mathfrak{a}_x(r_0)\, r_1^2/r_0^2 \;.$$

(2.3)

The latter with (2.2) imply the lower bound for $\mathfrak{a}_x(r)$ stated in the lemma. The upper bound on $\mathfrak{a}_x(r)$ follows from (2.14) given the aforementioned bound on the area of C. The monotonicity follows, for example, from Theorem 2.1 in [Ye]. However, a much simpler proof is had via integration by parts. To make the latter argument, take local coorinates centered at x where the symplectic form $\omega$ appears as the standard form in $\mathbb{R}^4$, thus $dx_1 \wedge dx_2 + dx_3 \wedge dx_4$. Write this as $\frac{1}{2} d(x_1 dx_2 - x_2 dx_1 + x_3 dx_4 - x_4 dx_3)$. Use integration by parts to express the integral of $\omega$ over the subvariety's intersection with the ball of radius r as a boundary integral. This can be parlayed using the tamed condition into a differential inequality for the function $r \to \mathfrak{a}_x(r)$ that integrates to give (2.3).

### c) The deformation operator

This subsection describes an essentially canonical Fredholm operator that is associated to an irreducible J-holomorphic subvariety and which is used to parametrize its J-holomorphic deformations. The discussion here has two parts. In both, C denotes an irreducible, J-holomorphic subvariety.

*Part 1*: Introduce again C's model curve $C_0$ with its tautological map, $\varphi$, to X. This map is *pseudoholomorphic*, which is to say that its differential intertwines the complex structure on $TC_0$ with the almost complex structure J on TX. To say this differently, introduce $T_{1,0}X \subset TX_\mathbb{C}$ to denote the subspace on which J acts as i, and introduce $T_{1,0}C$ to denote the subspace of $TC_\mathbb{C}$ on which its complex structure also acts by i. The differential of $\varphi$ extends as a $\mathbb{C}$-linear map from $TC_\mathbb{C}$ to $\varphi^*TX_\mathbb{C}$; and $\varphi$ is pseudoholomorphic if and only if this $\mathbb{C}$-linear map sends $T_{1,0}C$ to $T_{1,0}X$. Of interest first are deformations of $\varphi$ that are pseudoholomorphic as maps from the fixed complex curve $C_0$ to X.

Introduce $\exp_X$ to denote the metric's exponential map. Let $\eta$ denote a given section of $C^\infty(C_0, \varphi^*T^{1,0}X)$, and let $\eta_\mathbb{R}$ denote the corresponding section of $\varphi^*TX$. Then the map

$$\exp_X|_{\varphi(\cdot)}(\eta_{\mathbb{R}}(\cdot))$$

(2.4)

defines a smooth map from $C_0$ to $X$.

With $\eta$ fixed, a family of deformations of $\varphi$ that is parametrized by $[0, 1]$ has $s \in [0, 1]$ member given by replacing $\eta_{\mathbb{R}}$ in (2.3) by $s\eta_{\mathbb{R}}$. This family is pseudoholomorphic to first order in s if and only if $\eta$ is annihilated by a certain $\mathbb{R}$-linear, differential operator

$$\mathfrak{D}: C^\infty(C_0; \varphi^*T_{1,0}X) \to C^\infty(C_0; \varphi^*T_{1,0}X \otimes T^{0,1}C_0) .$$

(2.5)

To elaborate, note that there exists a $r_0 > 0$ such that if $|\eta|$ is everywhere bounded by $r_0$, then the map given by (2.4) is pseudoholomorphic if and only if $\eta$ obeys an equation that has the schematic form

$$\mathfrak{D}\eta + \mathfrak{R}_1(\eta)\cdot\nabla\eta + \mathfrak{R}_0(\eta) = 0 .$$

(2.6)

Here, $\mathfrak{R}_1$ and $\mathfrak{R}_2$ obey $|\mathfrak{R}_1(b)| \leq c_0|b|$ and $|\mathfrak{R}_0(b)| \leq c_0|b|^2$.

A deformation of $\varphi$ can have a J-holomorphic image but not define a J-holomorphic map from $C_0$ to $X$. It will, however be pseudoholomorphic with respect to a deformation of the original complex structure on $C_0$. More is said about this next. Assume until told otherwise that $C_0$ has genus at least 2.

To start, tensor the complex vector bundle homomorphism $\varphi_*: T_{1,0}C_0 \to \varphi^*T_{1,0}X$ with the identity homomorphism on $T^{0,1}C_0$ to extend $\varphi_*$ so as to give a bundle homomorphism from $T_{1,0}C_0 \otimes T^{0,1}C_0$ to $\varphi^*T_{1,0}X \otimes T^{0,1}C_0$. This extension is such that

$$\mathfrak{D}\varphi_* = \varphi_* \bar{\partial}$$

(2.7)

where $\bar{\partial}$ here maps $C^\infty(C_0; T_{1,0}C_0)$ to $C^\infty(C_0; T_{1,0}C_0 \otimes T^{1,0}C_0)$. Fix a complex, 3g-3 dimensional space of sections of $T_{1,0}C_0 \otimes T^{0,1}C_0$ that projects isomorphically to the cokernel of $\bar{\partial}$. Use V to denote this space of sections, and let

$$\Pi_V: C^\infty(C_0; \varphi^*T_{1,0}X) \to C^\infty(C_0; \varphi^*T_{1,0}X \otimes T^{0,1}C_0)$$

(2.8)

denote the $L^2$ orthogonal projection to the complement of $\varphi_*V$. Use the latter to define the operator

$$\mathfrak{D}_C = \Pi_V\mathfrak{D}: C^\infty(C_0; \varphi^*T_{1,0}X) \to \Pi_V C^\infty(C_0; \varphi^*T_{1,0}X \otimes T^{0,1}C_0) .$$

(2.9)

This operator appears in the following context: Let $\eta$ denote a given section of $\varphi^*T_{1,0}X$ with $|\eta| \leq r_0$. Then the image of the map in (2.4) is a J-holomorphic subvariety in X if and only if

$$\mathfrak{D}_C\eta + \Pi_V(\mathfrak{R}_1(\eta)\cdot\nabla\eta + \mathfrak{R}_0(\eta)) = 0 .$$
(2.10)

Put a Riemannian metric on $C_0$ which is compatible with its almost complex structure. Use this metric with the metric on $\varphi^*TX$ and the Levi-Civita connection on X to define Sobolev spaces of sections of $\varphi^*T_{1,0}X$ and its tensor product with $T^{0,1}C$. This done, then the operator $\mathfrak{D}_C$ extends to give a bounded, Fredholm operator from the Hilbert space $L^2_1(C_0, \varphi^*T_{1,0}X)$ to the orthogonal complement of V in the Hilbert space $L^2(C_0; \varphi^*T_{1,0}X \otimes T^{0,1}C_0)$.

A slight modification is needed when $C_0$ has genus 0 or 1 for in this case, there are non-trivial, holomorphic sections of $T_{1,0}C_0$. The latter define a complex vector subspace in $C^\infty(C_0; T_{1,0}C_0)$ with complex dimension equal to 3 when $C_0$ has genus 0 and equal to 1 when $C_0$ has genus 1. Use $\ker(\bar{\partial})$ to denote this subspace. In the genus 0 and genus 1 cases, the domain for the operator $\mathfrak{D}_C$ depicted in (2.10) is the quotient vector space $C^\infty(C_0; \varphi^*T_{1,0}X)/\ker(\bar{\partial})$. This last point will be implicit in what follows. In particular, the notation used below will indicate the actual domain only in the cases when the genus of $C_0$ is at least two.

The index of the Fredholm version of $\mathfrak{D}_C$ is denoted by $d_C$. Note that $d_C$ is in all cases even; and in all cases $d_C \leq \iota_{e_C} = e_C\cdot e_C - c\cdot e_C$. This is an equality if and only if $C_0$ has genus equal to $\frac{1}{2}(e_C\cdot e_C + c\cdot e_C) + 1$ and C is embedded. As noted in (2.1), the genus of $C_0$ is no greater than this. In general, the integer $d_C$ can be written as $d_C = \iota_{e_C} - 2\Delta$ when $C_0$ has genus $\frac{1}{2}(e_C\cdot e_C + c\cdot e_C) + 1 - \Delta$. Note that $d_C = \iota_{e_C} - 2$ if and only if C has a single immersion singularity. The case $d_C = \iota_{e_C} - 4$ occurs if and only if C has either two immersion singularities, or a single singularity that is either described by the n =1 version of (2.2), or it has a neighborhood whose intersection with C has three irreducible components, each a smoothly embedded disk.

*Part 2*: The Riemannian metric on $TX|_C$ can be used to realize the bundle N as a subbundle of $\varphi^*T_{1,0}X$. Let $N^\perp$ denote the orthogonal complement. The homomorphism $\varphi_*$ maps $T_{1,0}C$ into $N^\perp$ and can be used to identify the latter with $T_{1,0}C \otimes O(q)$ where q here denotes a sum that is indexed by the critical points of $\varphi$ and whose contribution from any given critical point is the integer n that appears in the relevant version of (2.2).

The operator $\mathfrak{D}$ appears in block diagonal form with respect to the splitting $T_{1,0}X = N \oplus N^\perp$ as

$$\begin{pmatrix} D_C & 0 \\ \mathfrak{R} & \bar{\partial} \end{pmatrix}.$$

(2.11)

where the notation is as follows: First, $D_C \colon C^\infty(C_0; N) \to C^\infty(C_0; N \otimes T^{0,1}C_0)$ is an $\mathbb{R}$-linear operator that acts on a given section $\varsigma$ as

$$D_C \varsigma = \bar{\partial}\varsigma + \nu\varsigma + \mu\bar{\varsigma}.$$

(2.12)

Here, $\bar{\partial}$ in (2.11) is the d-bar operator that is defined using the metric induced on N as a subbundle of $T_{1,0}X$, and $\nu$ is a section of N. Meanwhile, $\mu$ is a section of $N^2 \otimes T^{0,1}C_0$. What is denoted by $\mathfrak{R}$ in (2.11) is an $\mathbb{R}$-linear bundle homomorphism from N to $N^\perp$.

It follows from (2.11) that the respective kernel and cokernel of $\mathfrak{D}_C$ are canonically isomorphic to the kernel and cokernel of $D_C$ when $\varphi$ is an immersion. In fact, more is true if $\varphi$ is an immersion. If such is the case, then a deformation of the subvariety C is given by composing a section of N with a suitable exponential map from N to X. In particular, there exists such a map, $\exp_C$, that is defined on a small radius disk bundle $N_1 \subset N$ and has the following properties:

- $\exp_C$ *maps the zero section to* C; *and its differential along the zero section is an isomorphism from* $TN|_0$ *to* $\varphi^* T_{1,0}X$.
- $\exp_C$ *embeds each fiber of* $N_1$ *as a J-holomorphic disk in* X.

(2.13)

A construction of such a map is described in Section 5d of [T3].

Let $\eta$ denote a section of $N_1$. Then the image in X of the map $\exp_C(\eta(\cdot))$ is a J-holomorphic subvariety if and only if $\eta$ obeys an equation of the form

$$D_C \eta + \mathfrak{r}_1 \cdot \partial\eta + \mathfrak{r}_0 = 0,$$

(2.14)

where $\mathfrak{r}_1$ and $\mathfrak{r}_0$ are smooth, fiber preserving maps from $N_1$ to $\operatorname{Hom}(N \otimes T^{1,0}C; N \otimes T^{0,1}C)$ and to $N \otimes T^{0,1}C$ that obey $|\mathfrak{r}_1(b)| \leq c_0 |b|$ and $|\mathfrak{r}_0(b)| \leq c_0 |b|^2$.

## 3. Moduli spaces

This section summarizes certain facts about moduli spaces of J-holomorphic subvarieties. Of particular concern are the properties of these spaces and their constituent subvarieties when J is suitably generic. Most of what is said will not surprise experts.

### a) Moduli spaces of irreducible, J-holomorphic subvarieties

Fix a tamed almost complex structure J. Associate to a given class $e \in H^2(X; \mathbb{Z})$ a set, $\mathcal{M}_e$, that is defined as follows: Any given element $\Theta \in \mathcal{M}_e$ is a finite set of pairs, where each pair has the form $(C, m)$ with $C \subset X$ denoting an irreducible, J-holomorphic subvariety and m a postive integer. The set of pairs in $\Theta$ is further constrained so that no two of its pairs have the same subvariety component, and so that

$$\sum_{(C,m) \in \Theta} m e_C = e.$$

(3.1)

A topology on $\mathcal{M}_e$ is defined as follows: A sequence $\{\Theta_k\}_{k=1,2,\ldots}$ in $\mathcal{M}_e$ converges to a given element $\Theta$ if the following two conditions are met:

- $\lim_{k \to \infty} \left( \sup_{z \in (\cup_{(C,m) \in \Theta} C)} \text{dist}(z, \cup_{(C',m') \in \Theta_k} C') + \sup_{z' \in (\cup_{(C',m') \in \Theta_k} C')} \text{dist}(\cup_{(C,m) \in \Theta} C, z') \right) = 0.$
- $\lim_{k \to \infty} \sum_{(C',m') \in \Theta_k} m' \int_{C'} \upsilon = \sum_{(C,m) \in \Theta} \int_C \upsilon$ *for any given smooth 2-form* $\upsilon$.

(3.2)

This topology is used implicitly in what follows.

**Proposition 3.1**: *The space $\mathcal{M}_e$ is compact for any given $e \in H^2(X; \mathbb{Z})$. In particular, only finitely many classes in $H^2(X; \mathbb{Z})$ are of the form $e_C$ with $(C, m) \in \Theta$ and $\Theta \in \mathcal{M}_e$.*

*Proof of Proposition 3.1*: Given the area bound from Lemma 2.2, this is a standard application of the compactness theorems for J-holomorphic subvarieties as can be found, for example in [Wo], [Ye], or from what is said in Section 4 of [T4].

With $e \in H^2(X; \mathbb{Z})$ and an integer $k \geq 0$ fixed, denote by $\mathcal{M}_{e,k}$ the subspace of $\mathcal{M}_e$ that consists of sets of the form $\Theta = (C, 1)$ with C an irreducible, J-holomorphic subvariety whose fundamental class is Poincaré dual to e and whose model curve has genus k. The proposition that follows a structure theorem for $\mathcal{M}_{e,k}$ The assertion is only needed in what follows for the cases when $\mathfrak{D}_C$ has trivial cokernel; even so, a proof is given in Subsection Ac) of the appendix.

**Proposition 3.2**: *Given $C \in \mathcal{M}_{e,k}$, there exists a smooth map, $f$, from a neighborhood of $0$ in the kernel of $\mathfrak{D}_C$ to the cokernel of $\mathfrak{D}_C$; and there exists a homeomorphism from $f^{-1}(0)$ to a neighborhood of C in $\mathcal{M}_{e,k}$ sending 0 to C. Moreover, the subset of $\mathcal{M}_{e,k}$ where the cokernel of $\mathfrak{D}_{(\cdot)}$ is trivial has the structure of a smooth manifold; and the smooth structure is such that at any point in this set, the aforementioned homeomorphism from a neighborhood of 0 in the kernel of $\mathfrak{D}_{(\cdot)}$ is a smooth embedding onto an open set.*

**b) Moduli spaces when J is generic**

More can be said in the case when the almost complex structure is generic. To elaborate, introduce the notion of a *residual set* in the space of almost complex structures on TX. Such a set is a countable intersection of open dense sets. As such, it is also dense.

The next proposition is a standard result (see, e.g [McS]). Even so, a proof is given in Subsection Ac) of the appendix.

The proposition that follows introduces g to denote the integer $\frac{1}{2}$ (e·e + c·e) + 1.

**Proposition 3.3**: *Fix* $e \in H^2(X; \mathbb{Z})$ *and an integer* $k \in \{0, \ldots, g\}$. *There is a residual set in the space of almost complex structures such that if* J *comes from this set, then the cokernel of* $\mathfrak{D}_{(\cdot)}$ *is trivial at each point in* $\mathcal{M}_{e,k}$, *and so the latter has the structure of a smooth manifold whose dimension is* $\iota_e - 2(g - k)$.

With e and k fixed, the residual set referred to here is denoted in what follows by $\mathcal{J}_{e,k}$.

Introduce $\vartheta_e$ to denote the set of pairs of the form (e´, k´) where e´ $\in H^2(X; \mathbb{Z})$ is a equal to $e_C$ for some pair (C, m) in a set $\Theta$ from $\mathcal{M}_e$ and where k´ is the genus of C's model curve. Proposition 3.1 and (2.1) guarantee that $\vartheta_e$ has but a finite collection of elements. Set $\mathcal{J}_e = \cap_{(e´,k´) \in \vartheta_e} \mathcal{J}_{e´,k´}$.

The next proposition is a corollary of sorts to Proposition 3.3. To set the stage for this proposition, fix a positive class $e \in H^2(X; \mathbb{Z})$ and then define a set $\mathcal{S}_e \subset H^2(X; \mathbb{Z})$ as follows: A class $\hat{e} \in \mathcal{S}_e$ if and only if

- $0 < \hat{e} \cdot [\omega] \le e \cdot [\omega]$.
- $\hat{e} \cdot \hat{e} = -1$.

(3.3)

This is a finite set. There exists a subset $\mathcal{S}_{e*} \subset \mathcal{S}_e$ with the following property: If J is $\omega$-tamed, then $\hat{e} \in \mathcal{S}_{e*}$ is Poincaré dual to the fundamental class of an embedded, J-holomorphic sphere. Note in this regard that a class so represented for one particular $\omega$-tamed J is represented in this manner by all. This follows by virtue of the fact that the cokernel of the relevant version of the operator $D_{(\cdot)}$ is necessarily trivial. In any event, say that a class $e \in H^2(X; \mathbb{Z})$ is *positive* if

- e *is positive, which is to say that* $e \cdot [\omega] > 0$.
- $e \cdot e > 0$.
- $e \cdot \hat{e} > 0$ *for all* $\hat{e} \in \mathcal{S}_{e*}$.

(3.4)

The class e is *strictly positive* if the inequality in the third bullet is a strict, thus if $e \cdot \hat{e} > 0$ for all $\hat{e} \in \mathcal{S}_{e*}$. A class e that obeys the first two bullets here will obey the third if there exists an $\omega$-tamed almost complex structure whose version of $\mathcal{M}_e$ contains an element $\Theta$ with the following property: If $(C, m) \in \Theta$ then $e_C \cdot e_C \geq 0$. This is because distinct J-holomorphic subvarieties have positive local intersection numbers. In any event, there exist strictly positive classes. Indeed, given that $[\omega]$ comes from $H^2(X; \mathbb{Q})$, then some large multiple of $[\omega]$ will be integral and as such, automatically strictly positive.

With the stage set, what follows is the promised proposition.

**Proposition 3.4**: *Fix a positive class* $e \in H^2(M; \mathbb{Z})$. *There is a residual set of $\omega$-tamed almost complex structures inside $\mathcal{J}_e$ with the following property*: *Suppose that J is from this set.*
- *If $\Theta \in \mathcal{M}_e$ and $(C, m) \in \Theta$, then $\iota_{e_C} \geq 0$.*
- $\sum_{(C,m) \in \Theta} m \iota_{e_C} \leq \iota_e$.
- *The preceding is an equality if and only if $\Theta = (C, 1)$ with $C \in \mathcal{M}_{e,g}$. If $\Theta$ is not of this form, then $\sum_{(C,m) \in \Theta} m \iota_{e_C} \leq \iota_e - 2$.*

The residual subset inside $\mathcal{J}_e$ referred to here is denoted by $\mathcal{J}_{e1}$.

*Proof of Proposition 3.4*: What with Proposition 2.5, the condition $\iota_{e_C} \geq 0$ follows from the fact that the index of any given $(C, m) \in \Theta$ version of the operator $\mathfrak{D}_C$ is no greater than $\iota_{e_C}$. To see about the second assertion, note first that

$$\iota_e = \sum_{(C,m) \in \Theta} (m \iota_{e_C} + (m^2 - m) e_C \cdot e_C) + \sum_{(C,m),(C',m') \in \Theta: C \neq C'} m m' e_C \cdot e_{C'}$$

(3.5)

Each contribution to the right most sum in (3.5) is non-negative by virtue of the fact that any given intersection between distinct, J-holomorphic subvarieties has positive local interesection number. The second assertion of the proposition follows if it is the case that $m = 1$ if $(C, m) \in \Theta$ and $e_C \cdot e_C < 0$. Note in this regard that $e_C \cdot e_C \geq -1$ in any case; and if $e_C \cdot e_C = -1$, the C is an embedded sphere. Indeed, this follows from the adjunction formula for J-holomorphic subvarieties. Thus, $e_C \in \mathcal{S}_{e*}$. With this last point in mind, suppose that $(C, m) \in \Theta$ and $e_C \in \mathcal{S}_{e*}$. Then

$$e \cdot e_C = -m + \sum_{(C',m') \in \Theta: C' \neq C} m' e_{C'} \cdot e_C.$$

(3.6)

This last expression is non-negative by virtue of the fact that e is strictly positive. With (3.6) in mind, split $\Theta = \Theta_+ \cup \Theta_-$ where $\Theta_-$ contains the pairs of the form $(C, m)$ with C a sphere that has $e_C \cdot e_C = -1$. Now use (3.6) to write

$$\iota_e = \sum_{(C,m)\in\Theta_+} (m\,\iota_{e_C} + (m^2 - m)\,e_C\cdot e_C + m\sum_{(C',m')\in\Theta_+ : C'\neq C} m'\,e_{C'}\cdot e_C) + \sum_{(C,m)\in\Theta_-} m(1 + e\cdot e_C).$$
(3.7)

This last expression implies the second assertion of the proposition.

**c) Moduli spaces with marked points and with singularities**

Assume in what follows that $e \in H^2(X; \mathbb{Z})$ as described in Proposition 3.4 and that the almost complex structure comes from $\mathcal{J}_{e1}$. Subsequent arguments require a refined version of $\mathcal{M}_{e,k}$ to parametrize the J-holomorphic subvarieties that pass through a given set of points in X. To make a precise definition, let d denote the number of points under consideration. Given a class $e \in H^2(X; \mathbb{Z})$ set $g = \frac{1}{2}(e\cdot e + c\cdot e) + 1$. Fix $k \in \{0, 1, \ldots, g\}$ and introduce $\mathcal{M}_{e,k,d} \subset \mathcal{M}_{e,k} \times (\times_d X)$ to denote the subspace that consists of pairs (C, q) with $C \subset \mathcal{M}_{e,k}$ and $q = (x_1, \ldots, x_d) \subset \times_d X$ such that $\{x_i\}_{1\leq i\leq d} \subset C$.

The next two propositions use the term *image variety* to describe certain subsets of a smooth manifold. Here is the definition: Let Y denote a manifold and let n denote a positive integer no greater than the dimension of Y. A subset $\mathcal{S} \subset Y$ is a codimension n image variety if it is closed, and if each point in Y has a neighborhood that intersects $\mathcal{S}$ as the image via a smooth, proper map of a manifold with a finite number of connected components, none with dimension greater than dim(Y) - n.

**Proposition 3.5**: *For any given $k \in \{0, 1, \ldots, g\}$, the space $\mathcal{M}_{e,k,d}$ is a smooth, codimension 2d image variety in $\mathcal{M}_{e,k} \times (\times_d X)$. It is a smooth submanifold if $k = g$.*

A proof of this is in Subsection Ad) of the appendix.

The rest of this subsection is concerned with the structure of the set of J-holomorphic subvarieties near a subvariety with singular points.

Set $d = \frac{1}{2}\iota_e$ and. For each $k \in \{0, 1, \ldots, g\}$, use $\pi_{d-2}$ to denote the map from $\mathcal{M}_{e,k,d}$ to $\times_{d-2} X$ that is obtained by composing $\pi_d$ with the map to $\times_{d-2} X$ which is defined by writing $\times_d X$ as $(\times_2 X) \times (\times_{d-2} X)$ and then projecting to the $\times_{d-2} X$ factor. Meanwhile, $\pi_\mathcal{M}$ will denote the map from $\mathcal{M}_{e,k,d}$ to $\mathcal{M}_{e,k}$ that is obtained by restricting the projection map from $\mathcal{M}_{e,k} \times (\times_d X)$ to $\mathcal{M}_{e,k,d}$.

**Proposition 3.6**: *There is a residual set in $\mathcal{J}_{e1}$ such that if the almost complex structure is from this set, then the following is true: There is residual set of points in $\times_{d-2} X$ such that if $k \in \{0, 1, \ldots, g-3\}$ and if $\mathfrak{w}$ is from this residual set, then $\pi_{d-2}^{-1}(\mathfrak{w}) \subset \mathcal{M}_{e,k,d}$ is empty. If $k \in \{g - 2, g - 1\}$, it is an image variety defined by a proper map from a manifold of dimension $(4 - 2(g-k))$. In the case $k = g$, it is a smooth 4-dimensional submanifold. Moreover, there exists $\kappa \geq 1$ such that*
- *At most $\kappa$ subvarieties in $\mathcal{M}_{e,g-1}$ from $\pi_\mathcal{M}(\pi_{d-2}^{-1}(\mathfrak{w}))$ are singular at an entry of $\mathfrak{w}$.*

- $\pi_{\mathcal{M}}(\pi_{d-2}^{-1}(\mathfrak{w})) \subset \mathcal{M}_{e,g-2}$ *is a finite set. If* $C \in \mathcal{M}_{e,g-2}$ *comes via* $\pi_{\mathcal{M}}(\pi_{d-2}^{-1}(\mathfrak{w}))$, *then C's singular points are not entries of* $\mathfrak{w}$.
- *The* $\pi_{\mathcal{M}}$ *image in* $\mathcal{M}_e$ *of* $\pi_{d-2}^{-1}(\mathfrak{w}) \subset (\mathcal{M}_{e,g,d} \cup \mathcal{M}_{e,g-1,d} \cup \mathcal{M}_{e,g-2,d})$ *is closed.*
- $\mathfrak{w}$ *is a regular value for the map* $\pi_{d-2}$ *on* $\mathcal{M}_{e,g,d-2}$.

This proposition is proved in Subsection Ae) of the appendix.

Use $\mathcal{J}_{e2} \subset \mathcal{J}_{e1}$ that is referred to by this proposition. A point $\mathfrak{w} \in \times_{d-2} X$ from the residual set in $\times_2 X$ given by Proposition 3.6 is said below to be a *regular point*. Note, by the way, that a regular point must have distinct entries.

Assume in what follows that the almost complex structure is from $\mathcal{J}_{e2}$. To set the stage for the next proposition, suppose that $\mathfrak{w} \in \times_{d-2} X$ is a regular point and let $\mathcal{M}^{\mathfrak{w}} \subset \mathcal{M}_{e,g}$ denote the subset of curves that contain all entries of $\mathfrak{w}$. This is to say that $\mathcal{M}^{\mathfrak{w}}$ is the $\pi_{\mathcal{M}}$-image of $\pi_{d-2}^{-1}(\mathfrak{w}) \subset \mathcal{M}_{e,g,d-2}$. It is consequence of the fact that $\mathfrak{w}$ is a regular point that $\mathcal{M}^{\mathfrak{w}}$ is a smooth, 4-dimensional submanifold of $\mathcal{M}_{e,g}$. Given $(x, \mathfrak{w}) \in X \times (\times_{d-2} X)$, let $\mathcal{M}^{(x,\mathfrak{w})} \subset \mathcal{M}^{\mathfrak{w}}$ denote the subset of curves that contain $x$ and all entries of $\mathfrak{w}$.

Use $\mathcal{M}^{\mathfrak{w}}{}_X \subset \mathcal{M}^{\mathfrak{w}} \times X$ to denote the subspace of pairs $(C, x)$ such that $x \in C$. Let $\pi^{\mathfrak{w}}{}_X: \mathcal{M}^{\mathfrak{w}}{}_X \to X$ denote the map that is induced by the projection from $\mathcal{M}^{\mathfrak{w}} \times X$ to $X$. Write the entries of $\mathfrak{w}$ as $(w_1, \ldots, w_{d-2})$. The set of critical points of $\pi^{\mathfrak{w}}{}_X$ is the union of $(\cup_{1 \le m \le d-2} (\mathcal{M}^{\mathfrak{w}} \times w_m)) \subset \mathcal{M}^{\mathfrak{w}} \times X$ with a set $\mathcal{Z}^{\mathfrak{w}}{}_X \subset \mathcal{M}^{\mathfrak{w}}{}_X$. Use $\mathcal{Z}^{\mathfrak{w}} \subset \mathcal{M}^{\mathfrak{w}}$ to denote the image via the map that is induced by the projection from $\mathcal{M}^{\mathfrak{w}} \times X$ to its $\mathcal{M}^{\mathfrak{w}}$ factor, and use $\mathcal{Z}^{(x,\mathfrak{w})} \subset \mathcal{M}^{\mathfrak{w}}$ to denote the set of curves in $\mathcal{Z}^{\mathfrak{w}}$ that contain $x$ and all entries of $\mathfrak{w}$, thus the intersection between $\mathcal{Z}^{\mathfrak{w}}$ and $\mathcal{M}^{(x,\mathfrak{w})}$.

The upcoming proposition also refers to the symmetric, non-negative function, $\mathfrak{d}$, on $\mathcal{M}_e \times \mathcal{M}_e$ that is define by the following rule: Let $\Theta, \Theta'$ denote elements in $\mathcal{M}_e$. Then

$$\mathfrak{d}(\Theta, \Theta') = \sup\nolimits_{z \in (\cup_{(C,m) \in \Theta} C)} \mathrm{dist}(z, \cup_{(C',m') \in \Theta'} C') + \sup\nolimits_{z' \in (\cup_{(C',m') \in \Theta'} C')} \mathrm{dist}(\cup_{(C,m) \in \Theta} C, z') \, . \tag{3.9}$$

The function $\mathfrak{d}$ is used to measure distances on $\mathcal{M}_e$. This is to say that if $C \in \mathcal{M}_{e,g}$, then the $\mathfrak{d}$-*distance* of $C$ to a given set $\mathcal{K} \subset \mathcal{M}_e$ denotes the infimum of the set of numbers $\{\mathfrak{d}((C, 1), \Theta) : \Theta \in \mathcal{K}\}$.

By way of a reminder from what is said at the end of Part 1 in Section 2c, a subvariety in $\mathcal{M}_{e,g-1}$ has precisely one singular point, and this is an immersion singularity.

**Proposition 3.7**: *If* $\mathfrak{w} \in \times_{d-2} X$ *is a regular point, then the space* $\mathcal{M}^{\mathfrak{w}}{}_X$ *is a smooth, 6-dimensional manifold. Moreover, there exists a residual subset in* $\mathcal{J}_{e2}$ *whose significance is as follows. Take the complex structure from this smaller set.*

- *There exists a residual set $\mathcal{X}^2 \subset \times_2(\times_{d-2} X)$ whose points have the following property: If $(\mathfrak{w}, \mathfrak{w}')$ come from $\mathcal{X}^2$, then they have no entry in common; and there exists $\kappa > 1$ such that each curve in $\mathcal{M}^{\mathfrak{w}}$ has $\mathfrak{d}$-distance at least $\kappa^{-1}$ from $\mathcal{M}^{\mathfrak{w}'}$.*
- *Let $x \in X$. There exists a residual set $\mathcal{X}_x \subset \times_{d-2} X$ whose points have the following property: If $\mathfrak{w} \in \mathcal{X}_x$, then there exists $\kappa > 1$ such that*
  a) *$x$ has distance greater than $\kappa^{-1}$ from each entry of $\mathfrak{w}$.*
  b) *The distance from $x$ to any point in any subvariety from $\mathcal{M}_{e,g-2}$ in $\pi_\mathcal{M}(\pi_{d-2}^{-1}(\mathfrak{w}))$ is greater than $\kappa^{-1}$.*
  c) *At most $\kappa$ elements from $\pi_\mathcal{M}(\pi_{d-2}^{-1}(\mathfrak{w})) \subset \mathcal{M}_{e,g-1}$ contain $x$. Moreover, each such curve is rigid in the following sense: Let $C \subset \mathcal{M}_{e,g-1}$ denote such a curve. Introduce the corresponding operator $D_C$ as defined in (2.12). The kernel of $D_C$ has no non-trivial elements that vanish at $x$ and at all entries of $\mathfrak{w}$.*
  d) *The distance from $x$ to the immersion singular point of any subvariety from $\mathcal{M}_{e,g-1}$ in $\pi_\mathcal{M}(\pi_{d-2}^{-1}(\mathfrak{w}))$ is greater than $\kappa^{-1}$.*
  e) *The space $\mathcal{Z}^{\mathfrak{w}}$ is a 3-dimensional image variety.*
  f) *The space $\mathcal{Z}^{(x,\mathfrak{w})}$ is a 1-dimensional image variety.*

Subsection Af) of the appendix has the proof of this last proposition. By the way, Lemma A.15 in this same section of the appendix asserts that there is a residual set of almost complex structures with the following property: If J comes from this set, then the set $\mathcal{X}_x$ can be chosen so that any $\mathfrak{w} \in \mathcal{X}_x$ version of $\mathcal{M}^{(x,\mathfrak{w})}$ is a smooth, 2-dimensional manifold.

Use $\mathcal{J}_{e3}$ in what follows to denote the residual set in $\mathcal{J}_{e2}$ that is described by this last proposition.

**d) Curves in $\mathcal{M}^{\mathfrak{w}}$ near an entry of $\mathfrak{w}$**

Assume in what follows that the almost complex structure is from $\mathcal{J}_{e2}$. Suppose that $\mathfrak{w} \in \times_{d-2} X$ is a regular point and let $w$ denote an entry of $\mathfrak{w}$. Fix an orthonormal frame for $T_{1,0}X|_w$ to identify the space of complex, 1-dimensional subspace in $T_{1,0}X|_w$ with $\mathbb{CP}^1$. By definition, all curves in $\mathcal{M}^{\mathfrak{w}}$ contain the point $w$. As all curves in $\mathcal{M}^{\mathfrak{w}}$ are smooth, each has a tangent plane at $w$, and the assignment of tangent plane to curve defines a smooth map

$$\phi_w: \mathcal{M}^{\mathfrak{w}} \to \mathbb{CP}^1.$$

(3.10)

Use $\mathcal{Y}^w \subset \mathcal{M}^{\mathfrak{w}}$ to denote the set of critical points of $\phi_w$.

**Proposition 3.8**: *There exists a residual subset in $\mathcal{J}_{e3}$ whose significance is as follows. Take the complex structure from this smaller set. Let $x \in X$. The residual set $X_x$ from Proposition 3.7 can be chosen so that the following is also true: Take $\mathfrak{w}$ from this set. If $w$ is an entry of $\mathfrak{w}$, then $x$ is not on a curve from $\mathcal{Y}^w$.*

This proposition is proved in Subsection g) of the Appendix.

Use $\mathcal{J}_{e4}$ to denote the residual set in $\mathcal{J}_{e3}$ that is described by Proposition 3.8.

## 4. Proofs of Proposition 1.1-1.3

This section supplies the proofs to the propositions in Section 1.

### a) Proof of Proposition 1.1

Introduce the set $\mathcal{J}_{e1}$ as described in Proposition 3.4. It is a consequence of Propositions 3.5 that the conclusions of Proposition 1.1 hold with $\Sigma = \Sigma_e$ if $\mathcal{J}_{e*}$ is a subset of $\mathcal{J}_{e1}$. To elaborate, Suppose that $\Theta \in \mathcal{M}_e - \mathcal{M}_{e,g}$. Let n denote the number of elements in $\Theta$ and let $\mathfrak{p}$ denote an ordering of $\Theta$. Introduce the n-tuple $\Xi_{\Theta,\mathfrak{p}} = ((e_1, k_1), \ldots, (e_n, k_n))$ where the $i \in \{1, \ldots, n\}$ version of $(e_i, k_i)$ are defined as follows: Use $\mathfrak{p}$ to order the pairs that comprise $\Theta$ from 1 to n. Use C to denote the subvariety from the i'th pair in $\Theta$. Then $e_i = e_C$ and $k_i =$ genus(C) where C is the subvariety from the i'th pair in $\Theta$ with its ordering defined by $\mathfrak{p}$. Next, introduce $\mathcal{X}$ denote the set of distinct elements in the set $\{\Xi_{\Theta,\mathfrak{p}}: \Theta \in \mathcal{M}_e - \mathcal{M}_{e,g}$ and $\mathfrak{p}$ is an ordering of $\Theta\}$. Proposition 3.4 guarantees that $\mathcal{X}$ is a finite set. Meanwhile, Proposition 3.3 finds a residual subset in $\mathcal{J}_{e1}$ such that if J is from this subset, then all of the spaces in $\{\mathcal{M}_{(e',k')}\}_{(e',k') \in \Xi \text{ and } \Xi \in \mathcal{X}}$ are smooth manifolds of the asserted dimensions. Take J from this set. Given $\Xi \subset \mathcal{X}$, let $\mathcal{M}_\Xi = \times_{(e',k') \in \Xi} \mathcal{M}_{e',k'}$, let $\mathcal{M}_{\Xi,d} \subset \mathcal{M}_\Xi \times (\times_d X)$ denote the subspace whose elements are of the form $((C_1, \ldots, C_n), (x_1, \ldots, x_d))$ with $\{x_m\}_{1 \leq m \leq d} \in \cup_{1 \leq i \leq n} C_i$. Here, n denotes the number of elements in $\Xi$. Proposition 3.5 implies that these spaces are all image varieties in that each is the image of map from a smooth manifold whose dimension is 2d more than that of $\mathcal{M}_\Xi$. It is a consequence of Proposition 3.4 that this number is at most 4d - 2 except in the case when $\mathcal{M}_\Xi = \mathcal{M}_{e,g,d}$, in which case it is 4d. This the case, the compactness result from Proposition 3.1 implies that the image in $\times_d X$ of the map from $(\times_{\Xi \in \mathcal{X}} \mathcal{M}_{\Xi,d}) - \mathcal{M}_{e,g,d}$ to $\times_d X$ that comes from the projection of $(\times_{\Xi \in \mathcal{X}} \mathcal{M}_\Xi) \times (\times_d X)$ is a measure zero set. Indeed, if $\mathfrak{y} \in \times_d X$ is not in this set, then Proposition 3.1 implies that there are no subvarieties from $\mathcal{M}_e - \mathcal{M}_{e,g}$ that contain all points of $\mathfrak{y}$. It also implies that $\pi_d^{-1}(\mathfrak{y}) \subset \mathcal{M}_{e,g,d}$ is compact. This last point implies that $\pi_d^{-1}(\mathfrak{y})$ is finite if $\mathfrak{y}$ is not a critical point of $\pi_d$.

### b) Proof of Proposition 1.2

Take the almost complex structure from the set $\mathcal{J}_{e1}$ as described in Proposition 3.4. The proofs of the assertions made in the first bullet has four steps.

*Step 1*: Let $S_{crit} \subset \times_d X$ denote the set of critical values of the map $\pi_d$. This set is measurable and it has measure zero. If $\pi_d(\mathcal{K}) \subset S_{crit}$, then $\pi_d(\mathcal{K})$ is measurable, and it too has measure zero.

*Step 2*: Suppose that $\pi_d(\mathcal{K})$ is not contained in $S_{crit}$. Introduce $n \in \{1, 2, \ldots\}$ to be the smallest integer with the property that $\pi_d^{-1}(\eta) \cap \mathcal{K}$ has at most n elements when $\eta \in \pi_d(\mathcal{K}) - (\pi_d(\mathcal{K}) \cap S_{crit})$. The set $\pi_d(\mathcal{K}) - (\pi_d(\mathcal{K}) \cap S_{crit})$ has a stratification as $\mathcal{U}_1 \subset \mathcal{U}_2 \subset \cdots \subset \mathcal{U}_n = \pi_d(\mathcal{K}) - (\pi_d(\mathcal{K}) \cap S_{crit})$ where $\mathcal{U}_k$ is such that $\eta \in \mathcal{U}_k$ if and only if the set $\pi_d^{-1}(\eta) \cap \mathcal{K}$ has at most k elements. Let $\eta \in \mathcal{U}_n - \mathcal{U}_{n-1}$ and suppose that $(C, \eta) \in \mathcal{K}$. As $\eta$ is not a critical point of $\pi_d$, and as $\mathcal{K}$ is open, there is a neighborhood of $(C, \eta)$ to which $\pi_d$ restricts as a diffeomorphism onto its image in $\times_d X$. This implies that $\mathcal{U}_n - \mathcal{U}_{n-1}$ is open, and thus measurable with positive measure.

Introduce $\mathcal{V}_{n-1}$ to denote the complement in $\mathcal{U}_{n-1}$ of the closure of $\mathcal{U}_n$. If $\mathcal{V}_{n-1} = \emptyset$, then $\mathcal{U}_{n-1}$ has measure zero.

*Step 3*:. Suppose that $\mathcal{V}_{n-1} \neq \emptyset$, and let n´ denote the largest integer k with the property that $\mathcal{V}_{n-1} \cap \mathcal{U}_k \neq \emptyset$. By definition, n´ < n. Let $\eta \in \mathcal{V}_{n-1} \cap (\mathcal{U}_{n'} - \mathcal{U}_{n'-1})$ and suppose that $(C, \eta) \in \mathcal{K}$. As $\eta$ is not a critical point of $\pi_d$, and as $\mathcal{K}$ is open, there is a neighborhood of $(C, \eta)$ to which $\pi_d$ restricts as a diffeomorphism onto its image in $\times_d X$. Moreover, there exists such a neighborhood whose image is disjoint from $\mathcal{U}_n$. This implies that $\mathcal{V}_{n-1} \cap (\mathcal{U}_{n'} - \mathcal{U}_{n'-1})$ is also open. Thus, it too is measurable with positive measure.

Introduce $\mathcal{V}_{n'-1}$ to denote the complement inside the set $\mathcal{V}_{n-1} \cap \mathcal{U}_{n'-1}$ of the closure of $\mathcal{V}_{n-1} \cap (\mathcal{U}_{n'} - \mathcal{U}_{n'-1})$. If $\mathcal{V}_{n'-1}$ is empty, then $\mathcal{U}_{n'-1}$ has measure zero.

*Step 4*: If $\mathcal{V}_{n'-1}$ is not empty, replace n with n´ in Step 3 to prove that there exists n´´ < n´ such that $\mathcal{V}_{n'-1} \cap (\mathcal{U}_{n''} - \mathcal{U}_{n''-1})$ is also open. Continuing with this looping algorithm will decompose $\pi_d(\mathcal{K})$ as a finite union of open sets and measure zero sets.

Turn next to the assertion made by the proposition's second and third bullets. Let n denote the integer from Step 1 above. It then follows from Lemma 2.2, that the assignment $|\phi_\eta(\upsilon)| \leq c_0 n \, [\omega] \cdot e \sup_X |\upsilon|$.

To see that the function $\phi_{(\cdot)}(\upsilon)$ is measurable, return to the proof of the first bullet. This step decomposes $\pi_d(\mathcal{K})$ into a finite union of open sets and measure zero sets. Let $\mathcal{V}$ denote one of the open sets. It is a consequence of the way $\mathcal{V}$ is defined that any given $\eta \in \mathcal{V}$ has a neighborhood with the following property: Let U denote the neighborhood in question. Then $\pi_d$ restricts to $\mathcal{K} \cap \pi_d^{-1}(U)$ as a covering map with some $n_U \leq n$ sheets. For example, if $\mathcal{V} = \mathcal{U}_n - \mathcal{U}_{n-1}$, then $n_U = n$. This implies that $\phi_{(\cdot)}(\upsilon)$ varies smoothly on U. Taking all such U into consideration shows that $\phi_{(\cdot)}(\upsilon)$ is measurable.

Given the bound on $\phi_{(\cdot)}(\upsilon)$, and given that it is measurable, it then follows that $\Phi_{\mathcal{K}}(\cdot)$ defines a current. It is closed, non-negative and type 1-1 because it is a weighted average of such currents.

The fact that $\Phi_{\mathcal{K}}$ is non-trivial when $\pi_d(\mathcal{K})$ has positive measure follows from the fact $\Phi_{\mathcal{K}}(\omega)$ is no less than $[\omega]\cdot e$ times the measure of $\pi_d(\mathcal{K})$, this because the integral of $\omega$ on any given curve C from $\mathcal{M}_{e,g}$ is equal to $[\omega]\cdot e$.

### c) Proof of Proposition 1.3: The definition of $\mathcal{K}$

Fix e as in Proposition 3.4, and set $d = \frac{1}{2}\iota_e$ and $g = \frac{1}{2}(e\cdot e + c\cdot e) + 1$. Introduce $\mathcal{J}_{e4}$ as defined in Proposition 3.8 and take the almost complex structure from this last set. Use this complex structure to define $\mathcal{M}_{e,g,d}$. The set $\mathcal{K}$ for use in Proposition 1.3 will be chosen from a family of sets supplied by the upcoming Definition 4.1. The definition of this family requires some preliminary stage setting; and this is done in a two part digression that follows directly. The first part explains how to construct a fiberwise finite open set in $\mathcal{M}_{e,g}$ from any given open set with compact closure. The second part establishes notation.

*Part 1*: This part of the digression would not be necessary were the following known: There is a residual subset of $\mathcal{J}_{e4}$ such that if J is in this set, and if $\mathcal{K} \subset \mathcal{M}_{e,g,d}$ has compact closure, then there is an upper bound to the number of points in $\pi_d^{-1}(\cdot) \cap \mathcal{K}$. As was explained to the author by S. Givental this assertion is true for the generic smooth map from $\mathcal{M}_{e,g,d}$ to $\times_d X$ as a consequence of what is said in [B]. (Y. Eliashberg and J. Mather separately pointed the author to arguments for this assertion). It is a good bet that a generic choice of J makes $\pi_d$ generic in the appropriate sense; but a proof that such is the case appeared less than straightforward.

In any event, what follows constructs a fiberwise finite, open subset of a given open set in $\mathcal{M}_{e,g,d}$ with compact closure. Let $\mathcal{O} \subset \mathcal{M}_{e,g,d}$ denote the given open set with compact closure. Introduce $\mathcal{S}_{\mathcal{O}} \subset \mathcal{M}_{e,g,d}$ to denote the set of critical points of $\pi_d$ in the closure of $\mathcal{O}$. The set $\pi_d^{-1}(\pi_d(\mathcal{S}_{\mathcal{O}})) \cap \mathcal{O}$ is a closed set, so its complement is open. Let $\mathcal{O}^0$ denote this complement. The map $\pi_d: \mathcal{O}^0 \to \times_d X$ has no critical points, so it is locally a covering map onto the complement of $\pi_d(\mathcal{S}_{\mathcal{O}})$ in $\pi_d(\mathcal{O})$. Fix a locally finite, open cover of $\pi_d(\mathcal{O}^0)$. Such a cover exists with the following two properties: First, there exists $K \geq 1$ such that any given point in $\pi_d(\mathcal{O}^0)$ is contained in at most K sets from this cover. Second, let U denote a set from this cover. Then $\pi_d^{-1}(U)$ interects $\mathcal{O}^0$ as a disjoint union of open sets that are mapped diffeomorphically onto U by $\pi_d$. A cover of $\pi_d(\mathcal{O}^0)$ of this sort is said in what follows to be *fiberwise finite*. Let $\mathfrak{U}$ now denote such a fiberwise finite cover of $\pi_d(\mathcal{O}^0)$. Assign to each open set $U \subset \mathfrak{U}$ one of its $\pi_d$-inverse image sets in $\mathcal{O}^0$. Use $\mathcal{O}^U \subset \mathcal{K}^0$ to denote the associated set. Set $\mathcal{K}$ to be $\cup_{U \in \mathfrak{U}} \mathcal{O}^U$. This is an open, fiberwise finite subset of $\mathcal{O}$ that is mapped by $\pi_d$ onto the complement in $\pi_d(\mathcal{O})$ of the closed, measure zero set $\pi_d(\mathcal{S}_\mathcal{O})$. A subset $\mathcal{K} \subset \mathcal{O}$ as just described is said to be a $\pi_d$-*finite core* of $\mathcal{O}$.

*Part 2*: Let $\Lambda \subset \times_{d-2} X$ denote a chosen finite set of distinct regular values. With $\Lambda$ chosen, reintroduce the sets $\{\mathcal{Z}^\mathfrak{w}\}_{\mathfrak{w} \in \Lambda}$ from Proposition 3.7. Let $\pi^\mathfrak{w}{}_\mathcal{M}: \mathcal{M}^\mathfrak{w}{}_X \to \mathcal{M}_{e.g}$ denote the composition of the map induced by the projection from $\mathcal{M}^\mathfrak{w} \times X$ to $\mathcal{M}^\mathfrak{w}$ followed by the tautological inclusion of $\mathcal{M}^\mathfrak{w}$ into $\mathcal{M}_{e.g}$. Use $\mathcal{Z}_\Lambda$ to denote the subset $\cup_{\mathfrak{w} \in \Lambda} \mathcal{Z}^\mathfrak{w}$ in $\mathcal{M}_{e.g}$. Meanwhile, introduce for each $\mathfrak{w} \in \Lambda$ and each entry $w$ of $\mathfrak{w}$, the set $\mathcal{Y}^w$ as defined in Proposition 3.8. Then use $\mathcal{Y}_\Lambda$ to denote $\cup_{\mathfrak{w} \in \Lambda}(\cup_{w: w \text{ is an entry of } \mathfrak{w}} \mathcal{Y}^w)$.

**Definition 4.1**: *Fix $r > 0$, $s > 0$. Define the set $\mathcal{O}_{r,s,\Lambda} \subset \mathcal{M}_{e.g.d}$ to be the set whose pairs have the form $(C, \mathfrak{y})$ where*
- *$C$ has $\mathfrak{d}$-distance greater than $r$ from $\mathcal{Z}_\Lambda$, from $\mathcal{Y}_\Lambda$, and from $\mathcal{M}_e - \mathcal{M}_{e.g}$.*
- *$\mathfrak{y} = (x_1, x_2, \mathfrak{z})$ with $(x_1, x_2) \in X \times X$ and $\mathfrak{z} \in \times_{d-2} X$ such that*
  a) *The distance from $\mathfrak{z}$ to some point in $\Lambda$ is less than $s$.*
  b) *The distance from both $x_1$ and $x_2$ to any entry of $\mathfrak{z}$ is greater than $rs$.*
  c) *The distance from $x_1$ to $x_2$ is greater than $rs$.*

*Fix a $\pi_d$-finite core of $\mathcal{O}_{r,s,\Lambda}$ and denote the latter by $\mathcal{K}_{r,s,\Lambda}$.*

With regards to the size of $s$, keep in mind that the entries of any given point from the set $\Lambda$ are pairwise disjoint as each point in $\Lambda$ is a regular point. This understood, there exists $s_0 > 0$ such that if $s < s_0$, then the set of points in $\times_{d-2} X$ with distance $100 s$ or less from a point in $\Lambda$ is a union of disjoint balls, with each centered on a point in $\Lambda$. This upper bound for $s$ is assumed implicitly for what follows.

As is argued momentarily, the conclusions of Proposition 1.3 are satisfied if $\mathcal{K}$ is taken to be a suitable $r$, $s$ and $\Lambda$ version of $\mathcal{K}_{r,s,\Lambda}$. What follows directly explains why $\mathcal{K}_{r,s,\Lambda}$ for small $r$ and $s$ satisfies the criteria in Proposition 1.2 so as to guarantee that the $\mathcal{K} = \mathcal{K}_{r,s,\Lambda}$ version of $\Phi_\mathcal{K}$ is a closed, non-negative type 1-1 current. To this end, note first that $\mathcal{K}_{r,s,\Lambda}$ is open. Meanwhile, $\mathcal{K}_{r,s,\Lambda}$ is fiberwise finite by construction. The arguments given in the upcoming Part d of this section establish that $\mathcal{K}_{r,s,\Lambda}$ is non-empty when $r$ and $s$ are small. Granted this, then small $r$ and $s$ versions of $\mathcal{K}_{r,s,\Lambda}$ meet all of the criteria for Proposition 1.2.

**d) Proof of Proposition 1.3: The lower bound**

Fix $x \in X$ and adapted coordinates $(z, w)$ centered at $x$. Choose $t > 0$, but such that the ball of radius t centered at the origin in $\mathbb{C}^2$ is well within the domain of these coordinates. Fix a smooth, non-increasing function $\chi: [0, \infty) \to [0, 1]$ with value 1 on the interval $[0, \frac{1}{4}]$ and value 0 on $[\frac{1}{2}, \infty)$. Use $\chi_t$ to denote the function $\chi(t^{-1}|\cdot|)$ on $\mathbb{C}^2$. The lower bound asserted by Proposition 1.3 follows if

$$c_0^{-1} s^{4(d-2)} t^4 \leq \Phi_{\mathcal{K}_{r,s,\Lambda}} (\chi_t i dz \wedge d\bar{z}).$$

(4.1)

The verification of (4.1) for suitable $r$, $s$ and $\Lambda$ requires the following preliminary result.

**Proposition 4.2**: *There exist a finite set $\Lambda \subset \times_{d-2} X$ of regular values and $\sigma > 10^4$ whose properties are herewith described. Here is the first property: Suppose that $\mathfrak{w}$, $\mathfrak{w}'$ are distinct elements in $\Lambda$ and let $w$ and $w'$ denote respective entries of $\mathfrak{w}$ and $\mathfrak{w}'$. Then $\mathrm{dist}(w, w') > \sigma^{-1}$. Moreover, any pair of curves $C \in \mathcal{M}^{\mathfrak{w}}$ and $C' \in \mathcal{M}^{\mathfrak{w}'}$ are such that $\mathfrak{d}(C, C') \geq \sigma^{-1}$. To state the remaining properties, let $x \in X$ denote any given point. There exists $\mathfrak{w} \in \Lambda$ and an entry $w$ of $\mathfrak{w}$ such that the following is true: Let $\mathfrak{B} \subset \times_{d-2} X$ denote the ball of radius $\sigma^{-2}$ centered at $\mathfrak{w}$ and use $B \subset X$ to denote the ball of radius $\sigma^{-3}$ centered at $x$. Fix adapted coordinates $(z, w)$ centered at $x$ to identify $B$ with a neighborhood of the origin in $\mathbb{C}^2$. A radius $\sigma^{-3}$ disk $O_x \subset \mathbb{CP}^1$ exists such that*

- *Each $\hat{\theta} \in O_x$ has a lift $\theta = (\theta_z, \theta_w)$ to $S^3$ with $|\theta_z| > \sigma^{-1}$.*
- *Fix $\mathfrak{z} \in \mathfrak{B}$, $\hat{\theta} \in O_x$ and points $x_1 \neq x_2 \in B$ with distance less than $\frac{1}{4}\sigma^{-3}$ from $x$ and such that $x_2 - x_1 \in \mathbb{C}^2 - 0$ projects to $\hat{\theta}$. Set $\mathfrak{y} = (x_1, x_2, \mathfrak{z}) \in \times_d X$. Then $\pi_d^{-1}(\mathfrak{y}) \in \mathcal{M}_{e.g}$ is non-empty, and each curve in $\pi_d^{-1}(\mathfrak{y})$ has the following three properties:*
  a) *The curve has $\mathfrak{d}$-distance greater than $\sigma^{-1}$ from $\mathcal{M}_e - \mathcal{M}_{e.g}$, $\mathcal{Z}^{\mathfrak{w}}$ and $\mathcal{Y}^{w}$.*
  b) *The curve intersects the ball of radius $\sigma^{-3}$ centered at $x$ as the image of a map from the disk about the origin in $\mathbb{C}$ to $\mathbb{C}^2$ of the form $u \to x_1 + \theta u + \mathfrak{r}$ with*
     i) *$|\mathfrak{r}| \leq \sigma|u|(|x_1| + |u|)$ and $|d\mathfrak{r}| \leq \sigma(|x_1| + |u|)$,*
     ii) *$\theta = (\theta_z, \theta_w) \in \mathbb{C}^2$ is a unit vector with image $\hat{\theta}$ in $\mathbb{CP}^1$.*

This proposition is proved momentarily.

The four steps that follow directly derive the left hand inequality in (4.1) for suitable $r$ and $s$ using the set $\Lambda$ from Proposition 4.2.

*Step 1*: Let $\sigma$ denote the constant given by Proposition 4.2. Take $s < \sigma^{-2}$ and take $r < \sigma^{-3}$ to define $\mathcal{K}_{r,s,\Lambda}$.

*Step 2*: With x ∈ X given for use in (4.1), fix $\mathfrak{w} \in \Lambda$ and an entry $w$ of $\mathfrak{w}$ such that the triple $(x, \mathfrak{w}, w)$ obeys the conclusions of the second bullet in Proposition 4.2. Define $\mathfrak{B} \subset \times_{d-2} X$, $B \subset X$ and $O_x \in \mathbb{CP}^1$ as in Proposition 4.2. The following geometric considerations are used in the subsequent steps.

- *If $\mathfrak{z} \in \mathfrak{B}$, then each point in B has distance greater than $s^2$ from each entry of $\mathfrak{z}$.*
- *If $x_1$ and $x_2$ are in B and if their distance apart is greater than $\frac{1}{64}\sigma_\varepsilon^{-2}$, then their distance apart is greater than $rs$.*
- *Fix $\mathfrak{z} \in \mathfrak{B}$, and also a pair $\{x_1, x_2\} \subset B$ and set $\mathfrak{y} = (x_1, x_2, \mathfrak{z}) \in \times_d X$. If $C \in \pi_d^{-1}(\mathfrak{y})$ then it has $\eth$-distance greater than $r$ from $\mathcal{Z}_\Lambda$, $\mathcal{Y}_\Lambda$ and $\mathcal{M}_e - \mathcal{M}_{e,g}$.*

(4.2)

These conclusions all follow directly from the second bullet of Proposition 4.2 given the upper bounds given in Step 1 for $r$ and $s$.

*Step 3*: Fix $\mathfrak{z} \in \mathfrak{B}$. Given $x_1 \in B$, write its $\mathbb{C}^2$ coordinates as $(z', w') \in \mathbb{C}^2$. Granted this notation, fix $x_1$ so that its coordinate $z'$ obeys $|z'| \leq \frac{1}{32}\sigma^{-2}$. Meanwhile, require that the coordinate $w'$ obey $|w'| \leq \frac{1}{64}t$. Fix a second point, $x_2 \in B$, with distance greater than $\frac{1}{64}\sigma^{-2}$ from $x_1$ which is described by the second bullet of Proposition 4.2 for some choice of $\hat{\theta} \in O_x$. Suppose that $\mathfrak{y} = (x_1, x_2, \mathfrak{z})$ is not a critical value of $\pi_d$. With (4.2) understood, it follows that the second bullet of Proposition 4.2 supplies $C \in \mathcal{M}_{e,g}$ such that $(C, \mathfrak{y} = (x_1, x_2, \mathfrak{z}))$ is in $\mathcal{K}_{r,s,\Lambda}$.

*Step 4*: What with Item b) of the second bullet of Proposition 4.2, the norm of the restriction of $|dz|$ to C's intersection with the ball of radius t centered on x is greater than $c_\varepsilon^{-1}$. Given the aforementioned lower bound for $|dz|$, and given what is said in Lemma 2.2, it follows that the integral of $\chi_t \, dz \wedge d\bar{z}$ over C must be greater than $c_\varepsilon^{-1} t^2$.

*Step 5*: The final remarks in Steps 3 and 4 have the following implications: The integral in (4.1) is no less than the product of four factors: First, the factor of $c_\varepsilon^{-1} t^2$ given by the final remark in Step 4. Second, with $\mathfrak{z} \in \mathfrak{B}$ fixed and $x_1 \in B$ fixed as in Step 3, the volume of the set of $x_2$ as described in Step 3. The latter volume is at least $c_\varepsilon^{-1}$. By way of an explanation, note that the set of $\hat{\theta}$ to choose from has diameter greater than $c_\varepsilon^{-1}$, and as the distance from $x_1$ can vary in an interval of length at least $c_\varepsilon^{-1}$ so there is an $c_\varepsilon^{-1}$ volume's worth of choices for the point $x_2$. Third, with $\mathfrak{z} \in \mathfrak{B}$ fixed, the volume of the set of $x_1$ considered in Step 3. The latter is no smaller than $c_\varepsilon^{-1} t^2$. Finally, the volume of $\mathfrak{B}$. The latter is at least $c_\varepsilon^{-1} s^{4(d-2)}$.

Multiply these four factors to get what is asserted by (4.1).

The proof of Proposition 4.2 requires an auxilliary lemma that asserts a somewhat weaker version of what is asserted in Proposition 4.2.

**Lemma 4.3**: *There exists a finite set $\Lambda_0 \subset \times_{d-2} X$ of regular points, and $\sigma_0 > 10^4$ whose significance is now explained. Fix a point $x \in X$. Then there is a point $\mathfrak{w} \in \Lambda_0$ with the following properties: Use $\mathfrak{B} \subset \times_{d-2} X$ to denote the ball of radius $\sigma_0^{-1}$ at $\mathfrak{w}$ and use $B \subset X$ to denote the ball of radius $\sigma_0^{-2}$ centered at $x$. Fix adapted coordinates $(z, w)$ centered at $x$ to identify $B$ with a neighborhood of the origin in $\mathbb{C}^2$. A radius $\sigma_0^{-1}$ disk $O \subset \mathbb{CP}^1$ exists such that*

- *Each $\hat{\theta} \in O$ has a lift $\theta = (\theta_z, \theta_w)$ to $S^3$ with $|\theta_z| > \sigma_0^{-1}$.*
- *Fix $\mathfrak{z} \in \mathfrak{B}$, $\hat{\theta} \in O$ and points $x_1 \neq x_2 \in B$ with distance less than $\frac{1}{4}\sigma_0^{-2}$ from $x$ and such that $x_2 - x_1 \in \mathbb{C}^2 - 0$ projects to $\hat{\theta}$. Set $\mathfrak{y} = (x_1, x_2, \mathfrak{z}) \in \times_d X$. Then $\pi_d^{-1}(\mathfrak{y})$ is non-empty, and each curve in $\pi_d^{-1}(\mathfrak{y})$ has the following two properties:*
  a) *The curve has $\mathfrak{d}$-distance greater than $\sigma_0^{-1}$ from $\mathcal{M}_e - \mathcal{M}_{e,g}$.*
  b) *The curve intersects the ball of radius $\sigma_0^{-2}$ centered at $x$ as the image of a map from the disk about the origin in $\mathbb{C}$ to $\mathbb{C}^2$ of the form $u \to x_1 + \theta u + \mathfrak{r}$ with*
    i) *$|\mathfrak{r}| \leq \sigma |u|(|x_1| + |u|)$ and $|d\mathfrak{r}| \leq \sigma(|x_1| + |u|)$,*
    ii) *$\theta = (\theta_z, \theta_w) \in \mathbb{C}^2$ is a unit vector with image $\hat{\theta}$ in $\mathbb{CP}^1$.*

*Proof of Lemma 4.3*: The first task is to define the set $\Lambda_0$. To this end, fix $x \in X$ and choose a point from the residual set $\mathcal{X}_x$ that is supplied by Propositions 3.7 and 3.8. Use $\mathfrak{w}_x$ to denote this point. As argued momentarily, there is an open neighborhood, $U_x \subset X$, of $x$ and a constant $\kappa_x > 0$ such that if $x' \in U_x$, then the following is true:

- $x'$ *has distance greater than $\kappa_x^{-1}$ from each entry of $\mathfrak{w}_x$.*
- *The distance from $x'$ to any point in any subvariety from $\mathcal{M}_{e,g-2}$ in $\pi_\mathcal{M}(\pi_{d-2}^{-1}(\mathfrak{w}_x))$ is greater than $\kappa_x^{-1}$.*
- *At most $\kappa_x$ elements from $\pi_\mathcal{M}(\pi_{d-2}^{-1}(\mathfrak{w}_x)) \subset \mathcal{M}_{e,g-1}$ contain $x'$.*
- *The distance from $x'$ to the immersion singular point of any subvariety from $\mathcal{M}_{e,g-1}$ in $\pi_\mathcal{M}(\pi_{d-2}^{-1}(\mathfrak{w}_x))$ is greater than $\kappa_x^{-1}$.*

(4.3)

Indeed, the existence of an open set and $\kappa_x$ such that the first, second, and fourth bullets hold follows directly from Items a), b) and d) of the second bullet in Proposition 3.7. As explained next, the third bullet in (4.3) follows from Item c) of the second bullet in Proposition 3.7. To see how this comes about, assume that no open set and $\kappa_x$ exist that make the third bullet hold. This requires a sequence $\{x_k\}_{k=12,...} \subset X$ and a divergent sequence $\{n_k\}_{k=1,2,...}$ with the following property: There are $n_k$ elements in $\mathcal{M}_{e,g-1}$ that contain $x_k$ and all entries of $\mathfrak{w}_x$. Given Proposition 3.1 and (2.14), and also Item b) of the second bullet in Proposition 3.7, there must exist a curve $C \in \mathcal{M}_{e,g-1}$ that contains x and all entries of $\mathfrak{w}$; and whose version of kernel($D_C$) has a non-trivial element that vanishes at x and all entries of $\mathfrak{w}$. But, such an event is impossible given Item c) of the second bullet of Proposition 3.7.

To finish the task of defining $\Lambda_0$, note that the collection $\{U_x\}_{x \in X}$ is an open cover of X. Take a finite subcover, and set $\Lambda_0$ to denote the corresponding finite subset of $\{\mathfrak{w}_x\}_{x \in X}$. Set $\kappa_0$ to denote the supremum of the corresponding finite subset of the numbers $\{\kappa_x\}_{x \in X}$.

With $\Lambda_0$ now defined, consider next the constant $\sigma_0$. Assume that no such $\sigma_0$ exists with the asserted properties so as to derive some nonsense. This assumption requires a divergent, increasing sequence $\{s_k\}_{k=1,2,\ldots} \subset [1, \infty)$ with the following significance: Given the index k, there exist points $x_k \in X$, $\mathfrak{w}_k \in \Lambda$, and adapted coordinates $(z_k, w_k)$ centered on $x_k$ for which there is no disk $O_{x_k} \in \mathbb{CP}^1$ as described in the statement of the proposition. To elaborate, suppose that $\hat{O} \in \mathbb{CP}^1$ is any given disk of radius greater than $s_k^{-1}$ with the property that each $\hat{\theta} \in \hat{O}$ has a lift to $S^3 \subset \mathbb{C}^2$ as $(\theta_z, \theta_w)$ with $|\theta_z| > s_k^{-1}$. Then the following two conditions must hold:

<u>Condition 1</u>: *There is a point $\mathfrak{z}_k \in \times_{d-2} X$ with distance less than $s_k^{-1}$ from $\mathfrak{w}_k$.*

<u>Condition 2</u>: *There are points $x_{1k} \neq x_{2k} \in \mathbb{C}^2$ with distance less than $\frac{1}{2} s_k^{-2}$ from the origin such that $x_{1k} - x_{2k}$ projects to a point $\hat{\theta}_k \in \hat{O}$. In addition, one or more of the next three assertions is true.*

- $\pi_d^{-1}(x_{1k}, x_{2k}, \mathfrak{z}_k) = \emptyset$,
- *There is a curve in $\pi_d^{-1}(x_{1k}, x_{2k}, \mathfrak{z}_k)$ with $\mathfrak{d}$-distance less than $s_k^{-1}$ from $\mathcal{M}_e - \mathcal{M}_{e,g}$.*
- *There is a curve in $\pi_d^{-1}(x_{1k}, x_{2k}, \mathfrak{z}_k)$ whose intersection with the ball of radius $2 s_k^{-2}$ centered at $x_k$ is not the image of a map from $\mathbb{C}$ to $\mathbb{C}^2$ of the form $u \to x_{1k} + \theta u + \mathfrak{r}$ where $|\mathfrak{r}| \leq \frac{7}{8} s_k |u|(|x_{1k}| + |u|)$ and $|d\mathfrak{r}| \leq \frac{7}{8} s_k (|x_{1k}| + |u|)$,*

The derivation of nonsense from this data has five steps.

*Step 1*: This first step asserts that the n = 0 version of what is said in (2.2) is stable with respect to deformations of the pseudoholomorphic subvariety. The details are given by the lemma that follows. Note that the lemma does not require that $J \in \mathcal{J}_{e3}$ or that the class e obey the conditions from Proposition 3.4.

**Lemma 4.4**: *Suppose that e is a given class in $H^2(X; \mathbb{Z})$ and that J is an $\omega$-tamed almost complex structure. Fix $x \in X$ and suppose that $C \subset X$ is an irreducible, J-holomorphic subvariety that contains x and whose fundamental class is Poincaré dual to e. Fix an* *adapted coordinate chart centered at* x *so as to identify a neighborhood of* x *in* X *with a ball about the origin in* $\mathbb{C}^2$. *Suppose that* R > 1 *and that* C *appears in the radius* $R^{-1}$ *ball about the origin in* $\mathbb{C}^2$ *as the image of a map from* $\mathbb{C}$ *to* $\mathbb{C}^2$ *that has the form*

$$u \to \theta u + \mathfrak{r}$$

*where* $|\mathfrak{r}| < R|u|^2$ *and* $|d\mathfrak{r}| < R|u|$, *and such that* $\theta \in \mathbb{C}^2$ *has norm* 1. *Fix* $\varepsilon > 0$ *and there is a neighborhood of* (C, 1) *in* $\mathcal{M}_e$ *whose elements have the form* (C´, 1) *where* C´ *intersects the ball of radius* $\frac{1}{2} R^{-1}$ *about the origin, and this intersection is the image of a map from* $\mathbb{C}$ *to* $\mathbb{C}^2$ *of the form* $u \to x´ + \theta´ u + \mathfrak{r}´$ *with*

- $|x´| < \varepsilon$,
- $|\mathfrak{r}´| < R|u|(|x´| + |u|)$ *and* $|d\mathfrak{r}´| < R(|x´| + |u|)$,
- $\theta´ \in \mathbb{C}^2$ *is a unit vector with* $|\theta´ - \theta| < \varepsilon$.

***Proof of Lemma 4.4*** The assertion follows from what is said in [Wo], [Ye] or [McS] about limits of sequence of pseudoholomorphic maps.

*Step 2*: One can assume without loss that $\{x_k\}_{k=1,2,\ldots}$ converges. This the case, introduce $x \in X$ in what follows to denote the limit. Fix adapted coordinates (z, w) at x. The sequence of loci $\{w_k = 0\}_{k \gg 1}$ appear in this coordinate chart as a sequence of submanifolds in $\mathbb{C}^2$ that are very nearly complex lines. The corresponding sequence of points in $\mathbb{CP}^1$ has a convergent subsequence. Pass to a subsequence and relable consecutively from 1 so that these points converge. If necessary, rotate the (z, w) coordinates so that the limit is the w = 0 complex line. The sequence $\{\mathfrak{w}_k\}_{k=1,2\ldots}$ can also be assumed to converge. The limit, $\mathfrak{w}$, must be a point in $\Lambda_0$ such that the $(x, \mathfrak{w}, \kappa_0)$ obey the assertions of all bullets in (4.3) with equality allowed where (4.3) has an inequality.

*Step 3*: According to the second and third bullet in (4.3), there are at most $c_0$ singular subvarieties through x and all entries of $\mathfrak{w}$. The singularities of each such subvariety avoids x. This last conclusion follows from the second and fourth bullets of (4.3). Since the singularities of these subvarieties avoid x, they all have unique tangent planes at x. The collection of these tangent planes define a set $I_x \subset \mathbb{CP}^1$ of at most $c_0$ points. To summarize: If (C, 1) is *any* element in $\mathcal{M}_e$ such that C contains x and all entries of $\mathfrak{w}$, then C has a well defined tangent plane in $T_{1,0}X|_x$. This last conclusion has the following implication: If (C, 1) is as just described, then C must lie in $\mathcal{M}_{e,g}$ if its tangent plane at x defines a point in $\mathbb{CP}^1 - I_x$.

As $I_x$ is finite, there is a disk $O_x \subset \mathbb{CP}^1$ of radius $c_0^{-1}$ with the following property: Any given point in $O_x$ has distance greater than $c_0^{-1}$ from $I_x$; and any such point lifts to a point $\theta = (\theta_z, \theta_w) \in \mathbb{C}^2$ with $|\theta| = 1$ and $|\theta_z| > c_0^{-1}$. Let $O_x' \subset O_x$ denote the concentric disk with $\frac{1}{2}$ the radius of $O_x$.

*Step 4*: Fix a point $\hat{\theta} \in O_x'$ and a lift, $\theta$, to a point in $\mathbb{C}^2$ with norm 1. For each index k, fix points $x_{1k} \neq x_{2k}$ with distance less than $\frac{3}{4} s_k^{-2}$ from $p_k$ and such that $x_{1k} - x_{2k}$ projects to the point $\hat{\theta}$. According to Proposition 2.1, there is a set $\Theta_k \in \mathcal{M}_e$ that contains $x_{1k}, x_{2k}$ and all entries of $\mathfrak{z}_k$. The sequence $\{\Theta_k\}_{k=1,2,...}$ has a subsequence that converges in $\mathcal{M}_e$. Let $\Theta$ denote any limit of this sequence, and let $C = \cup_{(C',m') \in \Theta} C'$. This subvariety C must contain x and all entries of $\mathfrak{w}$. As $\mathfrak{w}$ is x-regular, the limit $\Theta$ must have the form (C, 1) where C is a subvariety which is non-singular at x. The fact that C is non-singular at x implies that each $C \in \mathcal{M}_{e,g}$. To explain these last remarks, note that (2.2) asserts that C intersects a ball of radius $c_0^{-1}$ centered at x as the image of a map from a disk about the origin in $\mathbb{C}$ to $\mathbb{C}^2$ that has the form $u \to \theta u + \mathfrak{r}(u)$ where $|\mathfrak{r}| \leq c_0 |u|^2$ and where $|d\mathfrak{r}| \leq c_0 |u|$. Since $\theta$ is not in $I_x$, the subvariety C must be in $\mathcal{M}_{e,g}$. This implies that each large k version of $\Theta_k$ has the form $(C_k, 1)$ with $C_k \in \mathcal{M}_{e,g}$. Given that $\{C_k\}_{k=1,2,...}$ converges to C and $C \in \mathcal{M}_{e,g}$, this in turn implies that each such $C_k$ has ∂-distance at least $c_0^{-1}$ from $\mathcal{M}_e - \mathcal{M}_{e,g}$.

*Step 5*: Lemma 4.4 asserts that each large k version of $C_k$ must intersect the ball of radius $c_0^{-1}$ about x as the image of a map from a disk about the origin in $\mathbb{C}$ to $\mathbb{C}^2$ that has the form $u \to x_{1k} + \theta u + \mathfrak{r}_k$ where $|\mathfrak{r}_k| \leq c_0 |u|(|p_k - x_{1k}| + |u|)$ and $|d\mathfrak{r}_k| \leq c_0 (|p_k - x_{1k}| + |u|)$. Given that such a $C_k$ exists for any $\hat{\theta} \in O_x'$ and any $x_{1k}, x_{2k}$ as described in Step 4, this conclusion is not compatible with the demands made by Condition 2.

With Lemma 4.3 in hand, turn now to the

***Proof of Proposition 4.2***: The proof has four steps.

<u>Step 1</u>: Fix $x \in X$ and let $\mathfrak{w}_{x0}$ denote a point from Lemma 4.3's set $\Lambda_0$ for which the conclusions of the lemma hold for the pair $(x, \mathfrak{w}_x)$. Introduce the set $\mathcal{X}_x \in X$ as described in Propositions 3.7 and 3.8. Fix a point $\mathfrak{w}_x \in \mathcal{X}_x$ that is very close to $\mathfrak{w}_{x0}$, and in particular, well inside Lemma 4.3's ball $\mathfrak{B}$. Let $w$ denote the entry of $\mathfrak{w}_x$ that is given by Proposition 3.8.

Fix an adapted coordinate chart (z, w) centered at x so as to identify a neighborhood of x with a ball about the origin in $\mathbb{C}^2$. Let $\mathcal{M}^{(x,\mathfrak{w}_x)}$ denote the space of curves in $\mathcal{M}_{e,g}$ that contain x and all entries of $\mathfrak{w}_x$. Any given curve $C \in \mathcal{M}^{(x,\mathfrak{w}_x)}$ has a well defined tangent plane at x. The identification between the neighborhood of x in X and $\mathbb{C}^2$ supplied by the adapted coordinates identifies this plane with a plane through the origin in $\mathbb{C}^2$, and the latter defines a point in $\mathbb{CP}^1$. This assignment of curve to point in $\mathbb{CP}^1$ defines a map $\phi_x: \mathcal{M}^{(x,\mathfrak{w}_x)} \to \mathbb{CP}^1$.

Step 2: Fix $\varepsilon$, and let $\mathcal{M}^{(x,\mathfrak{w}_x),\varepsilon} \subset \mathcal{M}^{(x,\mathfrak{w}_x)}$ denote the subset of curves with $\eth$-distance at least $\varepsilon$ from $\mathcal{M}_e - \mathcal{M}_{e,g}$. Given Item f) of Proposition 3.7, what follows is a consequence of the fact that $\mathcal{M}^{(x,\mathfrak{w}_x),\varepsilon}$ has compact closure in $\mathcal{M}_{e,g}$: If $\varepsilon$ is suitably generic, then $\phi_x(\mathcal{Z}^\mathfrak{w} \cap \mathcal{M}^{(x,\mathfrak{w}_x),\varepsilon})$ is a compact, 1-dimensional image variety in $\mathbb{CP}^1$. In particular, this is also the case with the latter's intersection with the closure of Lemma 4.3's disk $O \subset \mathbb{CP}^1$.

What follows is a consequence of this observation. There exists a constant $c_1 \geq 10^6$ that depends only on x, $\mathfrak{w}_x$ and $\varepsilon$; and there exists a disk $\hat{O} \subset O$ with radius $c_1^{-1}$ such that all points in $\hat{O}$ have distance $c_1^{-1}$ or more from $\phi_x(\mathcal{Z}^{\mathfrak{w}_x} \cap \mathcal{M}^{(x,\mathfrak{w}_x),\varepsilon})$, and distance $4c_1^{-1}$ or less from the point in $\mathbb{CP}^1$ defined by $(1,0) \in \mathbb{C}^2$.

This last conclusion has the following consequence: Any curve $C \in \mathcal{M}^{(x,\mathfrak{w}_x),\varepsilon}$ with $\phi_x(C) \in \hat{O}$ has $\eth$-distance at greater than $c_2^{-1}$ from $\mathcal{Z}^{\mathfrak{w}_x}$ where $c_2 > 10^6$ also depends only on x, $\mathfrak{w}_x$ and $\varepsilon$.

What with Proposition 3.8, the fact that $\mathcal{M}^{(x,\mathfrak{w}_x),\varepsilon}$ is compact implies that any curve C from this set has $\eth$-distance at least $c_3^{-1}$ from $\mathcal{Y}^w$. Here, $c_3 > 10^6$ also depends only on x, $\mathfrak{w}_x$ and $\varepsilon$.

Step 3: Let $B \subset X$ denote the set supplied by Lemma 4.3. What follows is a consequence of the last two conclusions from Step 2 with the fact that $\mathcal{M}^{(x,\mathfrak{w}_x),\varepsilon}$ is compact. There exist concentric balls $B_x \subset B$ and $\mathfrak{B}_x \subset \mathfrak{B}$ and a concentric disk $O_x \subset \hat{O}$ with the following consequence: Let $(x', \mathfrak{z}')$ denote a point in $B_x \times \mathfrak{B}_x$. Let $C \subset \mathcal{M}_{e,g}$ denote a curve with $\eth$-distance at least $2\varepsilon$ from $\mathcal{M}_e - \mathcal{M}_{e,g}$ that contains $x'$ and each entry of $\mathfrak{z}$. Suppose, in addition that the tangent plane to C at $x'$ defines a point in $O_x$. Then C has $\eth$ distance at least $\frac{1}{2} c_2^{-1}$ from $\mathcal{Z}^{\mathfrak{w}_x}$ and at least $\frac{1}{2} c_3^{-1}$ from $\mathcal{Y}^w$. With regards to C's tangent plane, the latter defines a vector in $T_{1,0}X|_{x'}$; but given that $x'$ is close to x, such a vector is very nearly in the span of the coordinate vectors $\{\frac{\partial}{\partial z}, \frac{\partial}{\partial w}\}$ and thus defines a point in $\mathbb{CP}^1$.

Step 4: The collection $\{B_x\}_{x \in X}$ defines an open cover of X, and so there is a finite subcover. Let $U \subset X$ denote the finite set that labels this subcover. Define $\Lambda'$ to be the corresponding finite set from $\{\mathfrak{w}_x\}_{x \in U}$. The set $\Lambda'$ obeys all but the first requirement set by Proposition 4.2. In particular, it may be the case that there are distinct pairs from $\Lambda'$ that do not lie in Proposition 3.7's set $\mathcal{X}^2$. However, given that each $x \in X$ version of $\mathcal{X}_x$ is residual, and given that $\mathcal{X}^2$ is residual, it follows from what is said in the previous steps that there exists a set $\{\mathfrak{w}_x'\}_{x \in U}$ that does obey all of Proposition 4.2's requirements, and is such that each $x \in U$ version of $\mathfrak{w}'_x$ is as close as desired to the corresponding $\mathfrak{w}_x$. This last conclusion follows from two observations. Here is the first: If $\delta > 0$ has been fixed in advance, and if $\mathfrak{w}'_x$ is very close to $\mathfrak{w}_x$, then the all curves in the $\mathfrak{w}'_x$ version of $\mathcal{Z}^\mathfrak{w}$ with $\eth$-distance at least $\frac{1}{2}\varepsilon$ from $\mathcal{M}_e - \mathcal{M}_{e,g}$ will have $\eth$-distance less than $\delta$ from $\mathcal{Z}^{\mathfrak{w}_x}$. By the same token, if $\mathfrak{w}'_x$ is very close to $\mathfrak{w}_x$, then the former has an entry, $w'$, that is very close to $w$; and any curve in $\mathcal{Y}^{w'} \subset \mathcal{M}^{\mathfrak{w}'_x}$ with $\eth$-distance at least $\frac{1}{2}\varepsilon$ from $\mathcal{M}_e - \mathcal{M}_{e,g}$ will have $\eth$-distance less than $\delta$ from $\mathcal{Y}^w$. To state the second observation, suppose that $x, x'$ are distinct pairs in $U \times U$. Let $\mathcal{X}^{(x,x')} \subset \times_U X$ denote the subset of points with the property that the pair of entries labeled by $x$ and $x'$ sit in $\mathcal{X}^2$. There is no constraint on the points in the other entries. This set is residual. As a consequence, the subset in $\cap_{(x,x') \in U \times U} \mathcal{X}^{(x,x')}$ with distinct entries is a residual subset in $\times_U X$ with the following property: Any pair of entries from any point in this set defines a point in $\mathcal{X}^2$.

Take $\Lambda$ to be such a set $\{\mathfrak{w}_x'\}_{x \in U}$.

### e) Proof of Proposition 1.3: A start on the upper bound

To obtain the desired the upper bound, introduce the notation used in (4.1). The upper bound asserted by Proposition 1.3 follows if

$$\Phi_{\mathcal{K}_{r,s,\Lambda}} (\chi_t i \, dz \wedge d\bar{z}) < c_{r,s} t^4$$

(4.4)

holds with $c_{r,s}$ here denoting a constant that depends on $r$ and $s$, but not on $t$ nor on $x$. As is explained in this and the next subsection, the bound in (4.4) holds if $r$ is first chosen small, and then $s$ is chosen with an $r$-dependent upper bound. Let $\sigma$ denote the constant from Proposition 4.2. The constraints $s < \sigma^{-2}$ and $r < \sigma^{-3}$ used to prove (4.1) are also assumed implicitly.

The upcoming lemmas are used to define this $r$-dependent upper bound. These lemmas together define an $(x, t)$-independent constant $\kappa_{**r} > 1$; and the proof of (4.4) requires that

$$s < \kappa_{**r}^{-2}.$$

(4.5)

This bound on s is the minimum of a various upper bounds that are given in Lemmas 4.7-4.14.

To prove (4.4), introduce $B_t(x)$ to denote the ball of radius t centered at x. The value of $\Phi(\chi_t i\, dz \wedge d\bar{z})$ is no greater than

$$c_0 \int_{\eta \in \pi_d(\mathcal{K}_{r,s,\Lambda})} \left(\sum_{C \in \pi_d^{-1}(\eta):\, C \cap B_t(x) \neq \emptyset} \int_C \chi_t \omega\right) \quad .$$

(4.6)

where it is understood that $\pi_d$ in this and subsequent formulae has domain $\mathcal{K}_{r,s,\Lambda}$. What with Lemma 2.2, this last expression is itself no greater than

$$c_0 t^2 \int_{\eta \in \pi_d(\mathcal{K}_{r,s,\Lambda})} \left(\sum_{C \in \pi_d^{-1}(\eta):\, C \cap B_t(x) \neq \emptyset} 1\right) \quad .$$

(4.7)

To continue, associate to each point $\mathfrak{w} \in \Lambda$, the ball $\mathfrak{B}_\mathfrak{w} \subset \times_{d-2} X$ of radius s. Associate to each $\mathfrak{z} \in \mathfrak{B}_\mathfrak{w}$ the set $X^2_\mathfrak{z} \subset X \times X$ of pairs $(x_1, x_2)$ such that $\eta = (x_1, x_2, \mathfrak{z}) \in \pi_d(\mathcal{K}_{r,s})$. Granted this notation, write (4.7) as

$$c_0 t^2 \sum_{\mathfrak{w} \in \Lambda} \int_{\mathfrak{z} \in \mathfrak{B}_\mathfrak{w}} \left(\int_{(x_1,x_2) \in X^2_\mathfrak{z}} \left(\sum_{C \in \pi_d^{-1}(x_1,x_2,\mathfrak{z}):\, C \cap B_t(x) \neq \emptyset} 1\right)\right) \quad .$$

(4.8)

There are six parts to the discussion that follows in this subsection Part 6 uses what is said in Parts 1-5 to derive an upper bound for the contribution to (4.8) from the points $\mathfrak{w} \in \Lambda$ whose entries all have distance $4s$ or more from x. There is at most one point in $\Lambda$ that does not have this property. The next subsection bounds the contribution from a point $\mathfrak{w} \in \Lambda$ with an entry having distance less than $4s$ from x. This subsection and the next use $c_{r,s}$ to denote an (x, t)-independent constant greater than 1, but dependent on r and s. Its value can be assumed to increase between subsequent appearances.

*Part 1*: To set the stage for Lemma 4.5, recall from the discussion just prior to Propostition 3.7 the definition of the space $\mathcal{M}^\mathfrak{w}$ with $\mathfrak{w}$ a given point $\mathfrak{w} \in \times_{d-2} X$. By way of reminder, this is the space of curves in $\mathcal{M}_{e,g}$ that contain all entries of $\mathfrak{w}$. As such, it is a smooth submanifold of $\mathcal{M}_{e,g,d-2}$ when $\mathfrak{w}$ is a regular point; and thus when $\mathfrak{w} \in \Lambda$. Given $\varepsilon > 0$, let $\mathcal{M}^{\mathfrak{w},\varepsilon} \subset \mathcal{M}^\mathfrak{w}$ denote the subspace of curves with $\partial$-distance greater than $\varepsilon$ from $\mathcal{M}_e - \mathcal{M}_{e,g}$. Lemma 4.5 also refers to the subspace $\mathcal{M}^{\mathfrak{w},\varepsilon}{}_X \subset \mathcal{M}^\mathfrak{w} \times X$ that consists of the pairs (C, x) with $C \in \mathcal{M}^{\mathfrak{w},\varepsilon}$ and $x \in C$. Use $\pi^\mathfrak{w}{}_\mathcal{M}$ and $\pi^\mathfrak{w}{}_X$ to denote the respective projection induced maps from this space to $\mathcal{M}^\mathfrak{w}$ and to X. Lemma 4.5 introduces one final bit of notation; this being $\pi_{\mathcal{M};d-2}$ to denote the map from $\mathcal{M}_{e,g,d-1}$ to $\mathcal{M}_{e,g,d-2}$ that is induced by the projection map from $\mathcal{M}_{e,g} \times (\times_{d-1} X)$ to $\mathcal{M}_{e,g} \times (\times_{d-2} X)$ that comes by writing the $\times_{d-1} X$ factor as $X \times (\times_{d-2} X)$ and then forgetting factor of X.

**Lemma 4.5**: *Given $r > 0$ but small, there exists $\kappa_r > 10^{10}$ with the following significance: Fix $\mathfrak{w} \in \Lambda$ and let $\mathfrak{B} \subset \times_{d-2} X$ denote the ball of radius $\kappa_r^{-2}$ centered on $\mathfrak{w}$.*

- *If $\mathfrak{z} \in \mathfrak{B}$, then $\mathcal{M}^{\mathfrak{z},r/4}$ is a smooth submanifold in $\mathcal{M}_{e,g}$.*
- *If $\mathfrak{z} \in \mathfrak{B}$ then $\mathcal{M}^{\mathfrak{z},r/4}{}_X$ is a smooth submanifold of $\mathcal{M}^{\mathfrak{z},r/4} \times X$. Moreover, if C is in $\mathcal{M}^{\mathfrak{z},r/2}$ and if it has $\mathfrak{d}$-distance $\frac{1}{16} r$ or more from $\pi^{\mathfrak{w}}{}_{\mathcal{M}}(\mathcal{Z}^{\mathfrak{w}})$, then C has distance $\frac{1}{32} r$ or more from the $\pi^{\mathfrak{z}}{}_{\mathcal{M}}$ image of any critical point of the map $\pi^{\mathfrak{z}}{}_X : \mathcal{M}^{\mathfrak{z},r/2}{}_X \to X$.*
- *There exists an embedding $\psi_X : \mathcal{M}^{\mathfrak{w},r/4}{}_X \times \mathfrak{B} \to \mathcal{M}_{e,g,d-1}$ onto an open set that*
  - a) *maps any given $\mathfrak{z} \in \mathfrak{B}$ version of $\mathcal{M}^{\mathfrak{w},r/2}{}_X \times \mathfrak{z}$ diffeomorphically onto an open set in $\pi_{d-2}^{-1}(\mathfrak{z})$ that contains $\mathcal{M}^{\mathfrak{z},r/2}{}_X$;*
  - b) *restricts to $\mathcal{M}^{\mathfrak{w},r/4}{}_X \times \mathfrak{w}$ as the inclusion map into the fiber of $\pi_{d-2}(\mathfrak{w})$;*
  - c) *is such that the composition $\psi = \pi_{\mathcal{M},d-2} \circ \psi_X$ is an embedding from $\mathcal{M}^{\mathfrak{w},r/4} \times \mathfrak{B}$ into $\mathcal{M}_{e,g,d-2}$.*

*Proof of Lemma 4.5*: Given that $\mathfrak{w}$ is a regular point, the implicit function theorem can be used in a straightforward way to obtain the constant $\kappa_r$ and the map $\psi_X$. The first two bulleted items follow directly from the existence of $\psi_X$.

The subsequent parts assume implicitly that $s$ obeys the bound $s < \kappa_r^{-2}$ with $\kappa_r$ as described in Lemma 4.5.

*Part 2*: Suppose that $\mathfrak{z} \in \times_{d-2} X$ has distance $s$ or less from a point in $\Lambda$. Lemma 4.5 guarantees that $\mathcal{M}^{\mathfrak{z},r/2}$ is a smooth, 4-dimensional submanifold of $\mathcal{M}_{e,g}$. Let C denote a curve in this submanifold, and let $N \to C$ denote C's normal bundle. Introduce the operator $D_C$ as given in (2.12). The tangent space to $\mathcal{M}^{\mathfrak{z}}$ at a curve C can be identified with the vector space $\ker_{C,\mathfrak{z}} \subset C^\infty(C; N)$ that consists of the sections in the kernel of $D_C$ that vanish at each entry of $\mathfrak{z}$. Use the $L^2$ inner product on $C^\infty(C; N)$ to define inner products and the norm on $\ker_{C,\mathfrak{z}}$. The norm is denoted by $\|\cdot\|_2$ in what follows. Note that $\|\cdot\|_2$ dominates any given $C^k$ norm on the elements of $\ker_{C,\mathfrak{z}}$.

**Lemma 4.6**: *Given $r > 0$, there exists $\kappa > 1$ with the following significance: Fix a point $\mathfrak{z} \in \times_{d-2} X$ with distance less than $s < \kappa_r^{-2}$ from a point in $\Lambda$. Let C denote a given curve from $\mathcal{M}^{\mathfrak{z},r/2}$. There is a diffeomorphism from the radius $\kappa^{-2}$ ball in $\ker_{C,\mathfrak{z}}$ onto an open set in $\mathcal{M}^{\mathfrak{z},r/4}$ that contains the set of curves in $\mathcal{M}^{\mathfrak{z},r/2}$ with $\mathfrak{d}$-distance less then $\kappa^{-3}$ from C. Moreover, there exists a diffeomorphism of the following sort: Introduce the map $\exp_C$* *that appears in (2.13). The map in question sends a given small normed vector $\eta \in \ker_{C,\mathfrak{z}}$ to $\exp_C(\eta + \phi_{C,\mathfrak{z}}(\eta))$ where $\phi_{C,\mathfrak{z}}: \ker_{C,\mathfrak{z}} \to C^\infty(C; N)$ is such that*

- *$\phi_{C,\mathfrak{z}}(\cdot) = 0$ at all entries of $\mathfrak{z}$.*
- *$\varsigma = \eta + \phi_{C,\mathfrak{z}}(\eta)$ obeys (2.14).*
- *Any given $C^k$ norm of $\phi_{C,\mathfrak{z}}(\eta)$ is bounded by $c_{r,k} \|\eta\|_2^2$ where $c_{r,k}$ depends only on k and r; and in particular, not on $\eta, C$, nor $\mathfrak{z}$.*

***Proof of Lemma 4.6***: Given that $\mathcal{M}^{3,r/2}$ is a smooth manifold, there exists $R_{C,\mathfrak{z}} > 1$ such that a map of the sort described by the first two bullets of the lemma exists with domain the ball of radius $R_{C,\mathfrak{z}}^{-2}$ in $\ker_{C,\mathfrak{z}}$ and with $\|\phi_{C,\mathfrak{z}}(\eta)\|_\infty \le R_{C,\mathfrak{z}} \|\eta\|_2^2$. Likewise, the k'th order derivatives are bounded by a contant $R_{C,\mathfrak{z},k}$. The fact that these constants can be chose to depend only on $r$, and for the latter, k, follows directly when Lemma 4.5's map $\psi_X$ is used to parametrize $\mathcal{M}^3$.

*Part 3*: Introduce again the subspace $\mathcal{M}^{3,r/2}{}_X \subset \mathcal{M}^{3,r/2} \times X$ to denote the subspace of pairs $(C, x)$ such that $x \in C$. Lemma 4.5 guarantees that this is a smooth, 6-dimensional submanifold of $\mathcal{M}^{3,r/2} \times X$. A pair $(C, x) \in \mathcal{M}^{3,r/2}{}_X$ is a critical point of the map $\pi^3_X$ to X if the subspace

$$\ker_{C,\mathfrak{z},x} = \{\eta \in \ker_{C,\mathfrak{z}} : \eta(x) = 0\}$$

(4.9)

has dimension greater than 2. This follows by virtue of the fact that the tangent space to $(C, x) \in \mathcal{M}^3{}_X$ is the vector space

$$\{(\eta, v) \in \ker_{C,\mathfrak{z}} \times TX|_x : \eta(x) - \Pi v = 0\},$$

(4.10)

where $\Pi$ here denotes the orthogonal projection from $TX|_C$ to the normal bundle N. Let $(\ker_{C,\mathfrak{z},x})^\perp \subset \ker_{C,\mathfrak{z}}$ denote the orthogonal complement to $\ker_{C,\mathfrak{z},x}$.

The next lemma concerns the points in the complement of the critical locus of $\pi^3_X$. The lemma refers to the distance between a given point in $\mathcal{M}^{3,r/2}{}_X$ and the critical locus of $\pi^3_X$. This distance is defined using the sum of the $\eth$-distance from $\mathcal{M}^{3,r/2}$ and the Riemannian distance from X.

**Lemma 4.7**: *Given $r > 0$ but small and given $\delta > 0$, there is a constant $\kappa > 1$ with the following significance: Suppose that $\mathfrak{z} \in \pi_{d-2}X$ has distance $s < \kappa_r^{-2}$ or less from a point in $\Lambda$. Let $(C, x) \in \mathcal{M}^{3,r/2}{}_X$ denote a pair with distance $\delta$ or more from the critical locus of $\pi^3_X$. Then $\|\eta\|_2 + \sup_C |\eta| \le \kappa |\eta(x)|$ if $\eta \in (\ker_{C,\mathfrak{z},x})^\perp$.*

*Proof of Lemma 4.7*: This follows from three observations: First, the bound in question holds with $\kappa$ replaced by some (C, x) and $\mathfrak{z}$ dependent contant since any two continuous norms on a finite dimensional vector space are equivalent. The latter constant can be chosen to be independent of (C, x) because the critical locus of $\pi^3_X$ on $\mathcal{M}^{3,r/2}{}_X$ is compact. More to the point, there is a uniform bound as (C, z) varies since it stays uniformly far from the critical locus. The fact that it can be chosen to be independent of $\mathfrak{z}$ is proved using the map $\psi$ from Lemma 4.5 to parametrize $\mathcal{M}^{3,r/2}{}_X$.

The critical points of $\pi^3_X$ are of two sorts. The first sort consists of the points of the form (C, x) where x is an entry of $\mathfrak{z}$. Use $\mathcal{Z}^{3,r/2} \subset \mathcal{M}^{3,r/2}{}_X$ to denote the complement in the critical locus of $\pi^3_X$ of the latter set.

*Part 4*: Suppose that $x \in X$ has distance at least $4s$ from all entries of $\mathfrak{w}$. Let $\mathcal{C}_{\mathfrak{z},x} \subset \mathcal{M}^{3,r/2}$ denote the subspace of curves with the following two properties:

- *The curve lies in the image via $\pi^3_\mathcal{M}$ of $(\pi^3_X)^{-1}(x) \subset \mathcal{M}^{3,r/2}{}_X$.*
- *The curve has $\mathfrak{d}$-distance greater than $\frac{1}{16}r$ from $\pi^\mathfrak{w}_\mathcal{M}(\mathcal{Z}^\mathfrak{w})$.*

(4.11)

The definition of $\mathcal{C}_{\mathfrak{z},x}$ is designed expressedly to keep all $C \in \mathcal{C}_{\mathfrak{z},x}$ versions of (C, x) away from the critical locus of $\pi^3_X$.

The next lemma states a consequence of this last observation.

**Lemma 4.8**: *Fix $r > 0$. There is a constant $\kappa_{*r} > \kappa_r$ such that if $s < \kappa_{*r}^{-2}$ then there exists another constant, $\kappa_{r,s} > 1$, with the following significance: Fix $\mathfrak{w} \in \Lambda$ and suppose that $\mathfrak{z} \in \times_{d-2} X$ has distance less than $s$ from $\mathfrak{w}$ and that $x \in X$ has distance at least $4s$ from all entries of $\mathfrak{w}$..*
- *The space $\mathcal{C}_{\mathfrak{z},x}$ is a smooth, dimension 2 submanifold in $\mathcal{M}^{3,r/2}$.*
- *Fix $\tau \in (0, \kappa_{r,s}^{-2})$. The submanifold $\mathcal{C}_{\mathfrak{z},x}$ is contained in a union of $\kappa_{r,s}\tau^{-2}$ balls in $\mathcal{M}^{3,r/4}$ with $\mathfrak{d}$-radius $\tau$.*

*Proof of Lemma 4.8*: The assertion that $\mathcal{C}_{\mathfrak{z},x}$ is a submanifold follows using the inverse function theorem from the fact that it lacks critical points of $\pi^3_X$. The assertion made by the second bullet follows directly from what is said by the final bullet in Lemma 4.5 using the Vitali covering lemma.

*Part 5*: Let $\kappa_{*r}$ and $\kappa_{r,s}$ denote the constants supplied by Lemma 4.8. Fix $\mathfrak{w} \in \Lambda$, take $s < \kappa_{*r}^{-2}$ and suppose again that $\mathfrak{z} \in \times_{d-2} X$ has distance $s$ or less $\mathfrak{w}$. Given t positive but less than the minimum of $s$ and $\kappa_{r,s}^{-2}$, let $\mathcal{B}_t(x; \mathfrak{z}) \subset \mathcal{M}^{\mathfrak{z},r/2}$ denote the set of curves of the following sort:

- *The curve contains all entries of $\mathfrak{z}$.*
- *The curve has $\mathfrak{d}$-distance at least $\frac{3}{4} r$ from $\mathcal{M}_e - \mathcal{M}_{e,g}$.*
- *The curve has $\mathfrak{d}$-distance at least $\frac{1}{4} r$ from $\pi^{\mathfrak{z}}_{\mathcal{M}}(\mathcal{Z}^{\mathfrak{z},r/2})$.*
- *The curve intersect the ball of radius $t$ centered at $x$.*

(4.12)

This set is the image in $\mathcal{M}^{\mathfrak{z},r/2}$ via the projection map from $\mathcal{M}^{\mathfrak{z},r/2} \times X$ of the set of pairs $(C', x') \in \mathcal{M}^{\mathfrak{z},r/2}{}_X$ that obey

- *$C'$ has $\mathfrak{d}$-distance at least $\frac{1}{4} r$ from $\pi^{\mathfrak{z}}_{\mathcal{M}}(\mathcal{Z}^{\mathfrak{z},r/2})$ and $\frac{3}{4} r$ from $\mathcal{M}_e - \mathcal{M}_{e,g}$.*
- *$x' \in B_t(x)$.*

(4.13)

Suppose that $x$ has distance at least $4s$ from all entries of $\mathfrak{w}$. Any pair $(C', x')$ that appears in (4.13) has distance at least $\frac{1}{4} r$ from the critical locus of $\pi^{\mathfrak{z}}_X$. The upcoming Lemma 4.9 asserts a consequence.

To prepare for the lemma, recall that each $C \in \mathcal{C}_{\mathfrak{z},x}$ has the associated vector space $(\ker_{C,\mathfrak{z},x})^\perp$. This space has dimension 2 because $(C, x)$ is not a critical point of $\pi^{\mathfrak{z}}_x$. Moreover, the collection of these spaces define a rank 2 vector bundle over $\mathcal{C}_{\mathfrak{z},x}$. Use $\mathcal{N}_{\mathfrak{z},x}$ to denote this bundle. It follows from Lemma 4.7 that the assignment to a pair $(C, \eta)$ with $C \in \mathcal{C}_{\mathfrak{z},x}$ and $\eta \in (\ker_{C,\mathfrak{z},x})^\perp$ of $|\eta(x)|$ defines a norm on $(\ker_{C,\mathfrak{z},x})^\perp$; and that these norms define a fiber norm on the vector bundle $\mathcal{N}_{\mathfrak{z},x}$. Disk subbundles of fixed radius in $\mathcal{N}_{\mathfrak{z},x}$ are defined using this fiber norm.

Lemma 4.9 identifies a fixed radius disk subbundle in $\mathcal{N}_{\mathfrak{z},x}$ with a neighborhood in $\mathcal{M}^{\mathfrak{z},r/4}$ of $\mathcal{C}_{\mathfrak{z},x}$. The identification involves the map $\exp_C$ from (2.13). The lemma refers also to the constant $\kappa_{*r}$ from Lemma 4.8

**Lemma 4.9**: *Fix $r > 0$ and then $s < \kappa_{*r}^{-2}$. There exists $\kappa \geq 1$ that depends only on $r$ and $s$ and has the following significance: Fix $\mathfrak{w} \in \Lambda$ and suppose that $x \in X$ has distance $t < \kappa^{-2}$ from all entries of $\mathfrak{w}$. Fix a point $\mathfrak{z} \in \times_{d-2} X$ with distance less than $s$ from $\mathfrak{w}$. There is an embedding from the radius $\kappa t$ disk bundle in $\mathcal{N}_{\mathfrak{z},x}$ onto a neighborhood in $\mathcal{M}^{\mathfrak{z},r/4}$ of $\mathcal{C}_{\mathfrak{z},x}$ with the properties listed below. Let $\lambda_{\mathfrak{z},x}$ denote this embedding. Then*

- *$\lambda_{\mathfrak{z},x}$ maps the zero section to $\mathcal{C}_{\mathfrak{z},x}$ as the identity map.*
- *The set $\mathcal{B}_t(x; \mathfrak{z})$ is an open set in the image of $\lambda_{\mathfrak{z},x}$ with compact closure in this image.*
- *Let $C \in \mathcal{C}_{\mathfrak{z},x}$ and $\eta \in (\ker_{C,\mathfrak{z},x})^\perp = \mathcal{V}_{\mathfrak{z},x}|_C$. The image of $\eta$ via $\lambda_{\mathfrak{z},x}$ can be written as*

$\exp_C(\eta + \phi_{C,\mathfrak{z},x}(\eta))$ *where* $\phi_{C,\mathfrak{z},x}(\cdot)$ *maps the ball* $\{\eta \in (\ker_{C,\mathfrak{z},x})^\perp : |\eta(x)| < 2\kappa^{-1}\}$ *smoothly into* $C^\infty(C; N)$. *Moreover, if* $\eta$ *is in its domain, then any given* $C^k$ *norm of* $\phi_{C,\mathfrak{z},x}(\eta)$ *is bounded by a multiple of* $|\eta(x)|^2$ *that depends only on* k, ε, *and* r.

*Proof of Lemma 4.9*: Given what is said in Lemmas 4.5-4.8, these conclusions all follow in a straightforward manner using the implicit function theorem.

*Part 6*: Fix $r > 0$ and $s < \kappa_* r^{-2}$ with the latter constant from Lemmas 4.8 and 4.9. Fix a point $\mathfrak{w} \in \Lambda$ and suppose that $x \in X$ has distance at least $4s$ from all entries of $\mathfrak{w}$. Let κ denote the constant supplied by Lemma 4.9 and assume that $t < \kappa^{-2}$. What follows derives an upper bound for the contribution to the integral in (4.8) from the subset in $\times_d X$ whose points can be written as $(x_1, x_2, \mathfrak{z})$ with $(x_1, x_2) \in \times_2 X$ and with $\mathfrak{z} \in \times_{d-2} X$ a point with distance less than $s$ from $\mathfrak{w}$.

The derivation has four steps.

Step 1: If $(x_1, x_2, \mathfrak{z})$ is as just described, then $(x_1, x_2) \in X \times X$ lie on a curve $C' \subset \mathcal{B}_t(x; \mathfrak{z})$. Lemma 4.9 supplies the following data: First, a curve $C \in \mathcal{C}_{\mathfrak{z},x}$ and a vector $\eta \in (\ker_{C,\mathfrak{z},x})^\perp$ that has norm $|\eta(x)| \leq \kappa t$. Second, points $p_1$ and $p_2 \in C$ such that $x_1 = \lambda_{C,\mathfrak{z},x}(\eta(p_1))$ and $x_2 = \lambda_{C,\mathfrak{z},x}(\eta(p_2))$. What with Lemma 4.7, this implies in particular that both $x_1$ and $x_2$ are constrained so as to lie in a tubular neighborhood in X of C whose radius is bounded by $c_{r,s} t$ ..

Step 2: According to Lemma 4.8 that there exists a set $\mathcal{U} \subset \mathcal{C}_{\mathfrak{z},x}$ containing at most $c_{r,s} t^{-2}$ curves such that any curve in $\mathcal{C}$ has ∂-distance at most $c_{r,s} t$ from a curve in $\mathcal{U}$. This implies that any pair $(x_1, x_2)$ as described in Step 1 must lie in a radius $c_{r,s} t$ tubular neighborhood of some curve from $\mathcal{U}$.

Step 3: The volume in X of the radius $c_{r,s} t$ tubular neighborhood of a curve from $\mathcal{U}$ is bounded by $c_{r,s} t^2$. Thus, the volume of the set of pairs $(x_1, x_2)$ in $X \times X$ that both lie in this tubular neighborhood is bounded by $c_{r,s} t^4$.

Step 4: It follows from what is said in Step 2 and 3 that the volume in $X \times X$ of the set of pairs that lie on some curve from $\mathcal{B}_t(x; \mathfrak{z})$ is bounded by the product of two factors: The first is the upper bound, $c_{r,s} t^4$, for the volume of the relevant radius tubular neighborhood in $X \times X$ of $C \times C$ when C is a curve from $\mathcal{U}$. The second is the upper bound $c_{r,s} t^{-2}$ to the number of curves in $\mathcal{U}$. Thus, the volume in $X \times X$ of the set of pairs that lie on some curve in $\mathcal{B}_t(x; \mathfrak{z})$ is bounded by $c_{r,s} t^2$.

Step 5: Since the volume in $\times_{d-2} X$ of the set of points $\mathfrak{z}$ under consideration is bounded by $c_0 s^{4(d-2)}$, it follows that the contribution to (4.8) from the subset in $\times_d X$ as described at the outset is no greater than $c_{r,s} t^2$ times (4.8)'s explicit $t^2$ factor; thus by $c_{r,s} t^4$.

**f) Proof of Proposition 1.3: Finishing the upper bound**

It remains still to bound (4.8) for points $x \in X$ with distance $4s$ or less from some entry of an element from $\Lambda$. The seven parts of this subsection derive a suitable bound. Note that all parts of this subsection use the following convention: Given $r > 0$, what is denoted by $c_r$ is a constant whose value is greater than $10^{10}$ and depends only on $r$. Its value can increase between successive appearances.

To set the stage for what follows, fix $\mathfrak{w} \in \Lambda$ and a given entry, $w$, of $\mathfrak{w}$.

*Part 1*: Given a $\mathfrak{z} \in \times_{d-2} X$ and an entry $z$ of $\mathfrak{z}$, reintroduce from Proposition 3.8 the map $\phi_z$ and its critical locus $\mathcal{Y}^z$ in the smooth part of $\mathcal{M}^{\mathfrak{z}}$. The following lemma is a corollary to Lemma 4.5.

**Lemma 4.10**: *Given $r > 0$ but small, the constant $\kappa_r$ that appears in Lemma 4.5 can be chosen so that the following additional conclusions hold: Fix $\mathfrak{w} \in \Lambda$ and let $\mathcal{B} \subset \times_{d-2} X$ again denote the ball of radius $\kappa_r^{-2}$ centered on $\mathfrak{w}$. Fix an entry $w$ of $\mathfrak{w}$. Let $\mathfrak{z} \in \mathcal{B}$ and let $z$ denote the entry of $\mathfrak{z}$ with distance less than $s < \kappa_r^{-2}$ from $w$. Then $\mathcal{M}^{\mathfrak{z}, r/2}$ consists of smooth points. Moreover, if $C \in \mathcal{M}^{\mathfrak{z}, r/2}$ has $\mathfrak{d}$-distance $\frac{1}{16} r$ or more from $\mathcal{Y}^w$ then it has $\mathfrak{d}$-distance greater than $\frac{1}{32} r$ from $\mathcal{Y}^z$.*

*Proof of Lemma 4.10*: This is a consequence of what is said about the map $\psi$ from Item c) of the third bullet in Lemma 4.5.

*Part 2*: Let $\mathfrak{z}$ be as described in Lemma 4.10 and let $C \in \mathcal{M}^{\mathfrak{z}}$ denote a curve that is not a critical point of $\phi_z$. This being the case, the differential of $\phi_z$ at $C$ maps $\ker_{C,\mathfrak{z}}$ surjectively to the tangent space $T\mathbb{CP}^1|_{\phi_z(C)}$. Use $\ker_{C,\mathfrak{z},\phi} \subset \ker_{C,\mathfrak{z}}$ to denote the kernel of the differential at $C$ of $\phi_z$, and use $(\ker_{C,\mathfrak{z},\phi})^\perp$ to denote the orthogonal complement in $\ker_{C,\mathfrak{z}}$ to $\ker_{C,\mathfrak{z},\phi}$.

The lemma that follows uses $\phi_{z*}$ to denote the differential at $C$ of $\phi_z$.

**Lemma 4.11**: *Given $r > 0$ but small, the constant $\kappa_r$ in Lemma 4.10 can be chosen so that the following is also true: Suppose that $\mathfrak{z}$ is as described in Lemma 4.10, and* suppose that $C \in \mathcal{M}^{3,r/2}$ has distance $\mathfrak{d}$-distance $\frac{1}{16} r$ or more from $\mathcal{Y}^w$. If $\eta \in (\ker_{C,\mathfrak{z},\phi})^\perp$ then $\|\eta\|_2 + \sup_C |\eta| \le \kappa_r |\phi_{z*}\eta|$.

*Proof of Lemma 4.11*: The argument for this differs only cosmetically from the argument used to prove Lemma 4.7.

*Part 3*: The contents of this part of the argument are summarized by the upcoming Lemma 4.12. This lemma refers to the constant $\kappa_r$ from Lemma 4.11.

**Lemma 4.12**: *Given $r > 0$, the constant $\kappa_r$ in Lemma 4.11 can be chosen so that given also $\Delta > 0$, there exists $\kappa > 1$ whose significance is as follows: Suppose that $s < \kappa_r^{-2}$ and that $\mathfrak{z} \in \times_{d-2} X$ is a point with distance at most $s$ from $\mathfrak{w}$. Given $\hat{\theta} \in \mathbb{CP}^1$, there exists a set $\mathcal{Q}_{\hat{\theta}} \subset \mathcal{M}^{3,r/4}$ of at most $\kappa$ curves such that*

- *Each curve in $\mathcal{Q}_{\hat{\theta}}$ is mapped by $\phi_z$ to $\hat{\theta}$.*
- *Each curve in $\mathcal{Q}_{\hat{\theta}}$ has $\mathfrak{d}$-distance at least $\frac{1}{32} r$ from $\mathcal{Y}^w$.*
- *Each curve in $\mathcal{M}^{3,r/2}$ that is*
  a) *mapped by $\phi_z$ to the disk of radius $\kappa^{-1}$ in $\mathbb{CP}^1$ centered at $\hat{\theta}$ and*
  b) *has $\mathfrak{d}$-distance at least $\frac{1}{32} r$ from $\mathcal{Y}^w$*
  
  *has $\mathfrak{d}$-distance less than $\Delta$ from some curve in $\mathcal{Q}_{\hat{\theta}}$.*

*Proof of Lemma 4.12*: Consider first this assertion for the case when $\mathfrak{z} = \mathfrak{w}$. If only $\mathfrak{z} = \mathfrak{w}$ is considered, then the only issue is that of $\kappa$'s existence. As explained next, such a constant exists for any value of $r$ as a consequence of the fact that the set of curves in $\mathcal{M}^{\mathfrak{w},r/2}$ with $\mathfrak{d}$-distance at least $\frac{1}{16} r$ from the critical point set of $\varphi_w$ is compact. To start the explanation, let $\mathcal{N} = \mathcal{N}(r)$ denote this same compact set. The fact that $\mathcal{N}$ is compact implies the following: Give $\varepsilon > 0$, there is a finite collection, $Q_\varepsilon \subset \mathcal{N}$, such that each point in $\mathcal{N}$ has $\mathfrak{d}$-distance $\varepsilon$ or less from a curve in $Q_\varepsilon$, and such that each point in $\phi_w(\mathcal{N})$ has distance $\varepsilon$ or less from a point in $\phi_w(Q_\varepsilon)$.

To continue, let $\hat{\theta} \in \mathbb{CP}^1$. The fact that $\mathcal{N}$ is compact and $\varphi_w$ is smooth has the next assertion as a consequence: Suppose that $\varepsilon < c_r^{-2}$. Let $D \subset \mathbb{CP}^1$ denote the disk of radius $c_r \varepsilon$ centered on $\hat{\theta}$. Introduce $Q_{\varepsilon,D} \subset Q_\varepsilon$ to denote the inverse image of $D$ via $\phi_w$. Then any point in $\mathcal{N}$ with $\phi_w$ image in the disk of radius $\varepsilon$ centered on $\hat{\theta}$ has $\mathfrak{d}$-distance at most $\varepsilon$ from a point in $Q_{\varepsilon,D}$. Meanwhile, the fact that each point in $Q_\varepsilon$ is uniformly far from a critical point of $\phi_w$ has the following consequence: Assume that $\varepsilon < c_r^{-2}$. Then each point in $Q_{\varepsilon,D}$ has $\partial$-distance less than $c_r \varepsilon$ from a point in $\mathcal{M}^{3,r/4}$ which is mapped by $\phi_w$ to $\hat{\theta}$ and which has $\partial$-distance greater than $\frac{3}{64} r$ from $\mathcal{Y}^w$. This understood, replace each point in $Q_{\varepsilon,D}$ by a point of the sort just described and let $Q_{\varepsilon,\hat{\theta}}$ denote the resulting set. Each point in $\mathcal{N}$ that is mapped by $\varphi_w$ to the disk of radius $\varepsilon$ centered on $\hat{\theta}$ has $\partial$-distance less than $c_r \varepsilon$ from some point in $Q_{\varepsilon,\hat{\theta}}$. This understood, fix some $\varepsilon < c_r^{-2}\Delta$ and set $\mathcal{Q}_{\hat{\theta}}$ to denote the corresponding set $Q_{\varepsilon,\hat{\theta}}$. Take $\kappa$ to be the maximum of $\varepsilon^{-1}$ and the number of elements in $Q_\varepsilon$ to obtain the conclusions of the lemma for the case $\mathfrak{z} = \mathfrak{w}$.

Granted the existence of $\kappa$ for the case $\mathfrak{z} = \mathfrak{w}$, the existence of Lemma 4.5's map $\psi$ implies that there are values for $\kappa_r$ and $\kappa$ that make the lemma hold when $\mathfrak{z}$ has distance $\kappa_r^{-2}$ or less from $\mathfrak{w}$. Here one uses the fact that points in the $\mathcal{N}(\frac{1}{1024} r)$ are uniformly far from the critical points of $\phi_w$ to see that the diffeomorphism $\psi$ can be modified so that the following is also true: The restriction to this set intertwines $\phi_z$ with $\phi_w$ if $\mathfrak{z}$ has distance less than $c_r^{-1}$ from $\mathfrak{w}$.

*Part 4*: The contents of this part of the argument are summarized by the next lemma. The lemma refers to the $r$-dependent constant $\kappa_r$ from Lemma 4.12.

**Lemma 4.13**: *Given $r > 0$ and given $s < \kappa_r^{-2}$, there exist a constant $\kappa_* > \kappa_r$ so that the following is true: Let $\mathfrak{z} \in \times_{d-2} X$ denote a point with distance less than $s$ from $\mathfrak{w}$ and let $x' \in X$ denote a point with distance at least $rs$ from all entries of $\mathfrak{z}$. Fix a point $\hat{\theta} \in \mathbb{CP}^1$ and fix $\delta < \kappa_*^{-2}$. Introduce $\mathcal{N}_{\mathfrak{z},x',\hat{\theta},\delta} \subset \mathcal{M}^{3,r/2}$ to denote the set of curves that have $\partial$-distance at least $\frac{1}{16} r$ from $\mathcal{Y}^z$, contain the point $x'$, and are mapped by $\phi_z$ to the disk of radius $\delta$ centered on the point $\hat{\theta}$. The volume in $X$ of the points that lie on curves in $\mathcal{N}_{\mathfrak{z},x',\hat{\theta},\delta}$ is bounded by $\kappa_* \delta^2$.*

***Proof of Lemma 4.13***: Let $\kappa$ denote the constant that appears in Lemma 4.6, and then set $\Delta = \frac{1}{2} \kappa^{-4}$. Invoke Lemma 4.12 to obtain a set $\mathcal{Q}_{\hat{\theta}}$ with the properties described by the lemma. Let C denote a curve from $\mathcal{Q}_{\hat{\theta}}$. Lemma 4.6 describes a diffeomorphism from the radius $\kappa^{-2}$ ball in $\ker_{C,\mathfrak{z}}$ onto an open set in $\mathcal{M}^{3,r/4}$ that contains the subset of curves in $\mathcal{M}^{3,r/2}$ with $\partial$-distance less than $\kappa^{-3}$ from C. Let $\mathcal{O}_C$ denote the latter set of curves and let $\mathcal{O}_{C,\delta} \subset \mathcal{O}_C$ denote the subset of curves that are mapped by by $\varphi_z$ to the radius $\delta$ disk centered on $\hat{\theta}$. Let $O_{C,\delta} \subset X$ denote the subset of points that lie on some curve from $\mathcal{O}_{C,\delta}$. It follows from Lemmas 4.6 and Lemma 4.11 that the volume of $O_{C,\delta}$ is no greater than $c_r \delta^2$.

If $\delta \leq c_r^{-1}$, then Lemma 4.12 implies the that $\mathcal{N}_{\mathfrak{z},x',\hat\theta,\delta} \subset (\cup_{C \in \mathcal{Q}_{\hat\theta}} \mathcal{O}_{C,\delta})$. Given that there are at most $c_r$ elements in $\mathcal{Q}_{\hat\theta}$, it follows from the conclusion of the previous paragraph that the volume of the set of points that lie on some curve from $\mathcal{N}_{\mathfrak{z},x',\hat\theta,\delta}$ is no greater than $c_r \delta^2$.

*Part 5*: The next lemma adds some to what is asserted by Lemma 4.4.

**Lemma 4.14**: *Fix $\varepsilon > 0$ and there is a constant, $\kappa_\varepsilon > 10^{10}$ with the following significance: Fix a point $p \in X$ and an adapted coordinate chart centered at p so as to identify a neighborhood of p with a ball about the origin in $\mathbb{C}^2$. Let $C \in \mathcal{M}_{e,g}$ denote a curve with $\mathfrak{d}$-distance $\varepsilon$ or more from $\mathcal{M}_e - \mathcal{M}_{e,g}$ that contains p. Then C intersects the ball of radius $\kappa_\varepsilon^{-2}$ centered p as the image of a map from a disk in $\mathbb{C}$ about the origin to $\mathbb{C}^2$ that has the form $u \to \theta u + \mathfrak{r}(u)$ where $\theta \in \mathbb{C}^2$ has norm 1 and where $\mathfrak{r}(u)$ is such that $|\mathfrak{r}(u)| \leq \kappa_\varepsilon |u|^2$ and $|d\mathfrak{r}|_v| \leq \kappa_\varepsilon |u|$.*

***Proof of Lemma 4.14***: This follows from Lemma 4.4 given that the set of curves under consideration has compact closure in $\mathcal{M}_{e,g}$.

Given $r > 0$, let $\kappa_r$ denote the constant from Lemma 4.12. It is assumed in what follows that $s$ is less than the minima of $\kappa_r^{-2}$ and $(\kappa_{r/4})^{-4}$ with $\kappa_{(\cdot)}$ the constant from Lemma 4.14.

*Part 6*: With (4.8) in mind, suppose that $x \in X$ has distance $4s$ or less from the entry $w$ of $\mathfrak{w}$. It is sufficient to bound the $\mathfrak{w} \in \Lambda$ contribution to (4.8) when $t < rs^2$. To this end, it proves useful to separate the integral over $\mathfrak{B}_\mathfrak{w}$ into two parts, that where the entry of $\mathfrak{z}$ near $w$ has distance on the order of t or less from x, and that where the distance from this entry to x has distance much greater than t from x.

To make a quantitative statement with regards to the first part, fix for the moment $R > 10^4$ so as to consider the contribution to the $\mathfrak{w} \in \Lambda$ term in (4.8) from the subset of points $\mathfrak{z} \in \mathfrak{B}_\mathfrak{w}$ whose entry, $z$, near $w$ is such that $\text{dist}(z, x) \leq Rt$. A suitable value for R is specified in the next part of this subsection.

Given that the number of curves in $\pi_d^{-1}(x_1, x_2, \mathfrak{z})$ has a bound that is independent of $x_1, x_2$ and $\mathfrak{z}$, it follows that the sum inside the $\mathfrak{w} \in \Lambda$ integral in (4.8) is bounded by a fixed, t independent constant, $c_0$. Meanwhile, the integral over the allowed set of pairs $(x_1, x_2)$ is no greater than the square of the volume of X, thus by $c_0$ also. Finally, the integral over the region of $\mathfrak{B}_\mathfrak{w}$ of interest is no greater than $c_0 s^{4(d-3)} R^4 t^4$. Thus, the contribution to the $\mathfrak{w} \in \Lambda$ term in (4.8) from the subset of points $\mathfrak{z} \in \mathcal{B}_{\mathfrak{w}}$ with $\text{dist}(z, x) < Rt$ is no greater than $c_0 \, s^{4(d-3)} R^4 t^6$.

*Part 7*: This part bounds the contribution to the $\mathfrak{w} \in \Lambda$ term from the subset of points $\mathfrak{z} \in \mathcal{B}_{\mathfrak{w}}$ whose entry near $w$ has distance much greater than t from x. As before, use $z$ to denote the entry of $\mathfrak{z}$ in question.

To start the analysis, introduce adapted coordinates $(z, w)$ centered at $z$ such that their identification of a neighborhood of $z$ with a ball about the origin in $\mathbb{C}^2$ makes x the point $(\mathfrak{d}, 0)$ with $\mathfrak{d} > 0$. Keep in mind that $\mathfrak{d} > c_0^{-1} \text{dist}(z, x)$ and so $\mathfrak{d} > c_0 Rt$. What follows is a key observation that is used in the analysis:

*If $\mathfrak{d} > 10^4 t$, then the set of complex 1-dimensional subspaces of $\mathbb{C}^2$ that intersect the ball of radius $t > 0$ centered at x is contained in a disk in $\mathbb{CP}^1$ of radius less than $c_0 t/\mathfrak{d}$ with center the image of $(1, 0)$ in $\mathbb{CP}^1$.*

(4.14)

Assume in what follows that R is such that $\mathfrak{d} > 10^4 t$.

Let $x_1$ and $x_2$ denote points in X that lie on a curve $C \in \mathcal{M}^{\mathfrak{z}, r/2}$ that intersects the ball of radius t centered at x. Given that these points appear in (4.8), their distance apart is at least $rs$. Keeping in mind that $t < rs^2$, this implies that at least one of them has distance $\Delta \geq t s^{-1} \gg t$ from x. No generality is lost by assuming that this is the point $x_1$. This point also has distance at least $rs$ from $z$.

With $\mathfrak{z}$ and $x_1$ fixed, it follows from (4.14) and Lemma 4.14 that the set of curves that contribute to any $x_2 \in X$ version of (4.8) lie in the set $\mathcal{N}_{\mathfrak{z}, x_1, \hat{\theta}, \delta}$ with $\hat{\theta}$ here denoting the image of the point $(1, 0)$ and with $\delta < c_r t/\mathfrak{d}$. This understood, then Lemma 4.13 has the following consequence: If $t/\mathfrak{d} \leq c_r^{-2}$, then the volume of the set of points $x_2 \in X$ that can appear in (4.8) with $x_1$ and $\mathfrak{z}$ fixed as above is no greater than $c_r \, t^2/\mathfrak{d}^2$. Meanwhile, with $\mathfrak{z}$ fixed, the volume of the set of points $x_1$ that can appear in (4.8) is bounded by the volume of X.

Granted the constraint $\mathfrak{d} > c_r^2 t$, granted the preceding conclusions, and given that $\mathfrak{d} \geq c_0^{-1} \text{dist}(z, x)$, the contribution to the integral in (4.8) from triples $(x_1, x_2, \mathfrak{z})$ with $\mathfrak{z} \in \mathcal{B}_{\mathfrak{w}}$ such that $\text{dist}(z, x) > c_r t$ is no greater than

$$c_r \, t^2 \, s^{4(d-3)} \int_{\text{dist}(\cdot, w) < s} \frac{1}{\text{dist}(\cdot, x)^2} \leq c_r t^2 \, s^{4d-10} .$$

(4.15)

Given what is said in Part 6, this implies that $\mathfrak{w} \in \Lambda$ contribution to (4.8) is no greater than $c_r s^{4d-10} t^4$.

**Appendix: Generic almost complex structures**

The purpose of this appendix is to supply the proofs to the various propositions in Section 3.

As a preview for what is to follow, note that the standard proofs of genericity assertions about almost complex structures invoke at some point the Sard-Smale theorem [Sm]. This is a generalization of Sard's theorem that applies to Fredholm maps between Banach spaces. In particular, its application requires introducing some sort of Banach space of almost complex structures on X. It is an unfortunate fact that the space of smooth almost complex structures on X is not a Banach space. The standard approach replaces the latter with a space of complex structures with bounded derivatives to some finite, but suitably large order. The Sard-Smale theorem is then invoked using such a space to obtain a residual set of complex structures with finite differentiability. This residual set of finitely differentiable almost complex structures is seen to be a countable intersection of open sets, each with open and dense intersection with the subspace of smooth, almost complex structures. As such, its intersection with the set of smooth almost complex structures is residual in the latter space. This is the usual route taken to obtain the desired residual set of smooth almost complex structures, and it is the route used here.

**a) Moduli spaces of pseudoholomorphic maps**

This part of the appendix sets the stage for what follows by setting up the machinery that is used in the subsequent parts of the appendix to prove the various genericity assertions that are made in Section 3. In particular, this first part of the appendix talks about moduli spaces of pseudoholomorphic maps as opposed to moduli spaces of pseudoholomorphic subvarieties. Most of what is done here consists of variations of what can be found with perhaps different notation in many sources (see, e.g. [McD]). Even so, the detailed discussion is warranted by virtue of the fact that pseudoholomorphic subvarieties are of interest here, and the treatment in the literature focuses for the most part on pseudoholomorphic maps. The discussion here also uses a 'Morrey space' as the Banach space for maps into X. The norm involves only the $L^2$ properties of the differential of a map; and as a consequence, norm bounds do not require estimates for higher order derivatives.

The discussion that follows has five parts culminating in Proposition A.5.

*Part 1*: Fix a smooth, oriented surface of genus k, here denoted by $\Sigma$, and fix a class $e \in H^2(X; \mathbb{Z})$. Choose a 'reference' complex structure, $j_0$, on $\Sigma$ and a compatible area form so as to define a Riemannian metric. Norms and covariant derivatives are defined using this metric. Select a reference smooth metric on X for the same purpose.

Fix $\upsilon \in (0, 1)$. The choice of such a constant determines a norm on the space of smooth maps from $\Sigma$ to X. This norm is denoted in what follows by $\|\cdot\|_*$; and it is defined by the rule

$$\|\varphi\|_*^2 = \sup_{p \in \Sigma} \sup_{d \in (0,1]} d^{-\upsilon} \int_{\text{dist}(p,\cdot)<d} |\varphi_*|^2$$

(A.1)

Here, $\varphi_*$ denotes the differential of $\varphi$. The next lemma describes the basic facts needed concerning this norm.

**Lemma A.1**: *There exists $\kappa > 1$ and given $\upsilon \in (0, 1)$, there exists $\kappa_\upsilon > 1$; and these have the following significance: Let $\varphi$ denote a smooth map from $\Sigma$ into X. Then the $L^2$ norm of $\varphi_*$ is bounded by $\kappa \|\varphi\|_*$ and the norm of $\varphi$ in the exponent $\frac{1}{2}\upsilon$ Holder space is bounded by $\kappa_\upsilon \|\varphi\|_*$.*

***Proof of Lemma A.1***: The $L^2$ bound on $\varphi_*$ follows directly from the definition. The Holder norm bound follows from Theorem 3.5.2 in [Mo].

Given maps $\varphi, \varphi'$ from $\Sigma$ to X, define their $*$-distance to be

$$D_*(\varphi, \varphi') = \sup_{p \in \Sigma} \text{dist}_X(\varphi(p), \varphi'(p)) + \sup_{p \in \Sigma} \sup_{d \in (0,1]} d^{-\upsilon} \int_{\text{dist}(p,\cdot)<d} |\varphi_* - \varphi'_*|^2.$$

(A.2)

Introduce $L_*(\Sigma_e, X)$ to denote the metric completion using $D_*$ of the space of smooth maps from $\Sigma$ to X that pushes forward the fundamental class of $\Sigma$ as the Poincaré dual of e. It is a consequence of Lemma A.1 that maps in this space are Holder continuous with exponent $\frac{1}{2}\upsilon$. As a consequence, $L_*(\Sigma_e, X)$ has the structure of a smooth Banach manifold. Let $\mathcal{E}_* \subset L_*(\Sigma_e, X)$ denote the subspace of maps that are somewhere 1-1. This is to say that there is at least one point in the image that has single point in its pre-image. This is an open set in $L_*(\Sigma_e, X)$. Let $\mathcal{E} \subset \mathcal{E}_*$ denote the subset of smooth maps.

Let $\varphi \in \mathcal{E}_*$. Since $\varphi$ is continuous, it pulls back vector bundles over X so as to give a vector bundle over $\Sigma$. In particular, $\varphi^*TX$ is a vector bundle over $\Sigma$ with exponent $\frac{1}{2}\upsilon$ Holder continuous transition functions. Moreover, the first derivatives of these transition functions are square integrable on their domain of definition, and are such that their integral over any disk of radius $d > 0$ is bounded by $c_0 d^{\upsilon/2}\|\varphi\|_*$. This understood, the notion of a Sobolev class $L^2_1$ section of $\varphi^*TX$ makes good sense. More to the point, there is an infinite dimensional space of sections of $\varphi^*TX$ on which the defined by

$$\|v\|_{1,2;\upsilon}^2 = \sup_{p\in\Sigma} \sup_{d\in(0,1]} d^{-\upsilon} \int_{\text{dist}(p,\cdot)<d} (|\nabla v|^2 + |v|^2)$$

(A.3)

is finite. The completion of the latter space is denoted in what follows by $L^2_*(\Sigma; \varphi^*TX)$. The assignment to any given map $\varphi \in L_*(\Sigma_e, \varphi^*TM)$ of the corresponding Banach space $L_*(\Sigma; \varphi^*TX)$ defines a smooth vector bundle over $\mathcal{E}_*$; this its tangent bundle.

Let $\iota_e = e \cdot e - c \cdot e$ and fix an integer $m > 10^{100}\iota_e$. Let $\mathcal{J}$ denote the Fréchet space of smooth, $\omega$-compatible almost complex structures on X and let $\mathcal{J}^m$ denote the corresponding Banach space of m-times differentiable, almost complex structures. When $J \in \mathcal{J}^m$, use $T_{J(1,0)}X$ to denote the (1, 0) part of $TX_\mathbb{C}$ as defined by J. The assignment to a pair $(\varphi, J) \in \mathcal{E}_* \times \mathcal{J}^m$ of $L_*(\Sigma, \varphi^*T_{J(1,0)}X)$ for $k \in \{0, 1, 2\}$ defines an m-1 times differentiable vector bundle over the Banach space $\mathcal{E}_* \times \mathcal{J}^m$.

*Part 2*: Given non-negative integers p and q, use $T_0^{p,q}\Sigma$ to denote the space of differentials of type (p, q) on $\Sigma$ as defined by the reference complex structure $j_0$. Suppose that the genus, k, of $\Sigma$ is greater than 1. Teichmuller theory asserts the existence of 3k-3 dimensional complex vector subspaces in $C^\infty(\Sigma, \text{Hom}(T_0^{1,0}\Sigma; T_0^{0,1}\Sigma))$ with the following property: Let V denote the subspace, and let j denote any given almost complex structure on $\Sigma$ whose derivatives to order m-2 are square integrable. There is an m-3 times differentiable homeomorphism of $\Sigma$ that pulls back j as a complex structure whose forms of type (1, 0) appear in $T^*\Sigma_\mathbb{C}$ as the graph of an element from V with pointwise norm less than 1. Conversely, any element from V with pointwise norm less than 1 defines a complex structure on $\Sigma$, this the complex structure whose forms of type (1, 0) appear in $T^*\Sigma_\mathbb{C}$ as the graph of the given element. In the case k = 1, there exists a 1-complex dimensional subspace of this sort.

To say a sentence more about the constraints on V, introduce $T_{0(1,0)}\Sigma$ to denote the holomorphic tangent bundle of $\Sigma$ and identify in the usual way $\text{Hom}(T_0^{1,0}\Sigma; T_0^{0,1}\Sigma)$ with the bundle $T_{0(1,0)}\Sigma \otimes T_0^{0,1}\Sigma$. A vector subspace $V \subset C^\infty(\Sigma, \text{Hom}(T_0^{1,0}; T_0^{0,1}))$ will suffice if projection to the cokernel of $\bar\partial : C^\infty(\Sigma, T_{0(1,0)}) \to C^\infty(\Sigma, T_{0(1,0)} \otimes T_0^{1,0})$ restricts as an isomorphism to V.

Fix such a subspace when $k \geq 1$ and let $V_1$ denote the set of elements in this subspace with pointwise norm less than 1. View $V_1$ as a finite dimensional Banach manifold with tangent space modeled on the vector space V with the norm given by the supremum norm for sections of $\text{Hom}(T_0^{1,0}\Sigma; T_0^{0,1}\Sigma)$.

As noted above, an element in $V_1$ defines a complex structure on $\Sigma$. If $j \in V_1$ denotes a given element, the corresponding almost complex structure on $T\Sigma$ is also denoted by j in what follows. Meanwhile, the corresponding space of differentials of type (p, q) is denoted by $T_j^{(p,q)}\Sigma$ and the (1, 0) tangent space of $\Sigma$ is denoted by $T_{j(1,0)}$. Note that an element in the tangent space to $V_1$ at a given $j \in V_1$ can be viewed in a canonical way as a section over $\Sigma$ of $T_{j(1,0)} \otimes T_j^{0,1}$. This view is used at times without accompanying remarks.

In the case of genus zero, set $V_1$ to be the zero dimensional Banach space.

*Part 3*: Let $(\varphi, J, j)$ denote a given element in $\mathcal{E}_* \times \mathcal{J}^m \times V_1$. Associated to the triple $(\varphi, J, j)$ is the Banach space of sections of $\varphi^* T_{J(0,1)} X \otimes T_j^{0,1} \Sigma \to \Sigma$ with norm whose square is given by the version of (A.1) that uses the section in question instead of $\varphi_*$. This norm is denoted by $\|\cdot\|_{2;\upsilon}$. The association of this Banach space to any given triple $(\varphi, J, j)$ defines an m-1 times differentiable vector bundle over $\mathcal{E}_* \times \mathcal{J}^m \times V_1$. This bundle is denoted below by $\mathcal{V}$.

What follows describes a canonical, m-1 times differentiable section of $\mathcal{V}$. It is denoted by $\eth$, and it is defined by the following rule: If $(\varphi, J, j)$ is such that $\varphi$ is smooth, it is the element in $C^m(\Sigma; \varphi^* T_{J(0,1)} X \otimes T_j^{0,1} \Sigma)$ that is obtained from the section

$$v \to (J|_\varphi + i) \varphi_*(j - i) v .$$
(A.4)

of $\mathrm{Hom}(T\Sigma_{\mathbb{C}}, \varphi^* T_{J(1,0)} X)$ by restriction to $T_{j(1,0)} \Sigma$. Here, $\varphi_*$ denotes the differential of $\varphi$. Note that its definition is such that $\varphi_*$ intertwines $j$ and $J$ if and only if $\eth = 0$ at $(\varphi, J, j)$.

Let $\mathcal{P}^m$ denote the zero locus of $\eth$. The implicit function theorem asserts that $\mathcal{P}^m$ is a m-1 times differentiable Banach manifold on a neighborhood of any element where $\eth$ vanishes transversely. It is left for the reader to verify that such is the case along the whole of $\mathcal{P}^m$.

The implicit function theorem identifies the tangent spaces to $\mathcal{P}^m$ at a given point $(\varphi, J, j)$. This the space of sections $(\varphi', J', j')$ of $\varphi^* T_{J(1,0)} X \oplus \mathrm{End}(TX) \oplus V$ with $\varphi'$ in $L_*(\Sigma; \varphi^* TX)$ and $J'$ in $C^m$ that obey an equation that has the schematic form

$$(J|_\varphi + i) \nabla \varphi'(j - i) + (\nabla_\varphi J) \cdot \varphi_*(j - i) + J' \varphi_*(j - i) + (J|_\varphi + i) \varphi_* j' = 0 .$$
(A.5)

What follows is a consequence of the inverse function theorem: Given $(\varphi, J, j) \in \mathcal{P}^m$, there exists $\delta > 0$ and an m-1 times differentiable homeomorphism with m-1 times differentiable inverse from the radius $\delta$ ball in vector space of solutions to (A.5) onto a neighborhood of $(\varphi, J, j)$ in $\mathcal{P}^m$ that maps the origin to $(\varphi, J, j)$, and whose differential at the origin is the identity map.

To say more about this homeomorphism, and to connect with what is said in Part 1 of Section 2, let $(\varphi, J, j)$ denote a given element in $\mathcal{P}^m$. A look at (A.5) reveals that the vector space of elements $(\varphi', J' = 0, j' = 0)$ in $T\mathcal{P}^m|_J$ is isomorphic to the subspace of sections of $L_*(\Sigma; \varphi^* T_{J(1,0)} X)$ that obey a certain first order differential equation given by a Fredholm differential operator whose domain is the $L^2_1$ completion of the space of sections over $\Sigma$ of the bundle $\varphi^* T_{J(1,0)} X$ and whose range is the $L^2$ completion of the space of sections of $(\varphi^* T_{J(1,0)} X) \otimes T_j^{0,1} \Sigma$. To elaborate, let $\varphi'_{(1,0)}$ denote the $(1, 0)$ part of $\varphi'$.

When written in terms of $\eta$, this is the operator that is denoted by $\mathfrak{D}$ in (2.5) with it understood that the complex structure on $C_0$ is that defined by j. In general, an element $(\varphi', J' = 0, j')$ is in $T\mathcal{P}^m|_J$ if and only if

$$\mathfrak{D}\varphi'_{(1,0)} + (J|_\varphi + i)\varphi_* j' = 0.$$

(A.6)

The notation used here is admittedly sloppy because $j'$ in (A.6) is to be viewed as a section $T_{j(1,0)} \otimes T_j^{0,1}$. Let $\mathfrak{D}_*$ denote the operator from $L^2_1(\Sigma, \varphi^* T_{J(1,0)} X) \oplus V$ to the Hilbert space $L^2(\Sigma, \varphi^* T_{J(1,0)} X \otimes T_j^{1,0}\Sigma)$ that is defined by the left hand side of (A.6). This is a Fredholm operator. introduce $\mathfrak{D}_*$ to denote the latter operator and let $\ker(\mathfrak{D}_*)$ and $\mathrm{coker}(\mathfrak{D}_*)$ denote its respective kernel and cokernel. Note in this regard that this vector space $\mathrm{coker}(\mathfrak{D}_*)$ is canonically isomorphic to the cokernel of the operator $\mathfrak{D}_C$ in (2.9).

The implicit function theorem applied to (A.5) finds a smooth map from a ball about the origin $\ker(\mathfrak{D}_*) \times T\mathcal{J}^m|_J$ to $\mathrm{coker}(\mathfrak{D}_*)$ that maps the point $(0, 0)$ to 0, with the following significance: Let $\mathfrak{B}$ denote the afore-mentioned ball and let $f$ denote the map. Then 0 is a regular value for $f$ and $f^{-1}(0)$ is homeomorphic to a neighborhood of $(\varphi, J, j)$ in $\mathcal{P}^m$ via an m-1 times differentiable map with m-1 times differentiable inverse. This homeomorphism is the restriction to $f^{-1}(0)$ of a map, $\Psi$, from $\mathfrak{B}$ into a neighborhood of $(\varphi, J, j)$ in $\mathcal{E}_* \times \mathcal{J}^m \times V_1$. To say slightly more about $f$, note that it has a Taylor's expansion of the form

$$f(\mathfrak{v}, J') = f_{0,1} J' + f_{1,0} \mathfrak{v} + \cdots$$

(A.7)

where $f_{0,1}$ and $f_{1,0}$ are linear maps and the unwritten terms are higher order. A particularly important point for what follows is that $f_{0,1}$ is surjective. Note also that the tangent space to $\mathcal{P}^m$ at $(\varphi, J, j)$ is isomorphic to the kernel in $\ker(\mathfrak{D}_*) \times T\mathcal{J}^m|_J$ of the linear map $(\mathfrak{v}, J') \to f_{0,1} J' + f_{1,0} \mathfrak{v}$.

To say a bit more about $\Psi$, identify a neighborhood of J in $\mathcal{J}^m$ with a a ball about the origin in the space of m-times differentiable sections of $\mathrm{Hom}(T_J^{1,0} X; T_J^{1,0} X)$. Then $\Psi$ can be viewed so as to send any given pair $(\mathfrak{v}, J') \in \mathfrak{B}$ to an element that can be written as

$$\Psi(\mathfrak{v}, J') = (\varphi(\mathfrak{v}, J'), J', j(\mathfrak{v}, j)).$$

(A.8)

The following lemma is a useful result concerning the elements in $\mathcal{P}^m$.

**Lemma A.2**: *Suppose that $J \in \mathcal{J}^m$ has continuous derivatives to order $m' \geq m$ or is infinitely differentiable. If $(\varphi, J, j) \in \mathcal{P}^m$, then $\varphi$ has continuous derivatives up through order $m'$ - 2 or is infinitely differentiable as the case may be.*

***Proof of Lemma A.2***: This follows using standard elliptic regularity techniques as can be found, for example, in Chapter 6 of [Mo].

This lemma implies that $\mathcal{P}^m$ has a second Banach space topology, this coming from its inclusion in $C^{m-2}(\Sigma; X) \times \mathcal{J}^m \times V_1$. The upcoming Lemma A.3 asserts that these two topologies are equivalent.

Let $\mathcal{P} \subset \mathcal{P}^m$ denote the subset of triples $(\varphi, J, j)$ such that $\varphi$ and $J$ are smooth. Introduce $\pi_{\mathcal{J}} \colon \mathcal{P}^m \to \mathcal{J}^m$ to denote the restriction of the projection map from $\mathcal{E}_* \times \mathcal{J}^m \times V_1$. Lemma A.2 implies that $\mathcal{P} = \pi_{\mathcal{J}}^{-1}(\mathcal{J})$. The space $\mathcal{P}$ has a Frêchet topology from the $C^\infty$ topology on the space of smooth maps from $\Sigma$ to $X$ and on $\mathcal{J}$. It also has the topology coming from the inclusion in $\mathcal{P}^m$. A part of the Lemma A.3 asserts that these two topologies are also equivalent.

**Lemma A.3**: *The Banach space topology on $\mathcal{P}^m$ coming from its inclusion in the space $C^{m-2}(\Sigma; X) \times \mathcal{J}^m \times V_1$ is equivalent to that coming from its inclusion in $\mathcal{E}_* \times \mathcal{J}^m \times V_1$. By the same token, the $C^\infty$-Frêchet space topology on $\mathcal{P}$ is equivalent to the topology coming from its inclusion as a subspace of $\mathcal{E}_* \times \mathcal{J} \times V_1$. The $C^\infty$-Frêchet space topology endows $\mathcal{P}$ with the structure of a smooth, Frêchet manifold.*

***Proof of Lemma A.2***: The fact that the two topologies on $\mathcal{P}^m$ are equivalent follows from the fact that the elliptic regularity techniques from Chapter 6 in [Mo] control the norms of the derivatives a map from an element in $\mathcal{P}^m$ in terms of the $L_*$ norm of the map and the norms of the derivatives of the almost complex structure. These same norm bounds imply that the $C^\infty$ Frêchet topology on $\mathcal{P}$ is the same as the topology that comes from its inclusion in $\mathcal{E}_* \times \mathcal{J} \times V_1$. The structure on $\mathcal{P}$ of a smooth Frêchet manifold comes via the implicit function theorem in the same way that the latter endows $\mathcal{P}^m$ with the structure of an m-1 times differentiable Banach manifold.

*Part 4*: It is a consequence of (A.5) that the differential of the map $\pi_{\mathcal{J}}$ from $\mathcal{P}^m$ to $\mathcal{J}^m$ has at all points a Fredholm differential. Set $g = \frac{1}{2}(e \cdot e + c \cdot e) + 1$. The index of the differential of $\pi_{\mathcal{J}}$ is $\iota_e - (g - k)$ if $k \geq 1$. In the case $k = 1$, the index is $\iota_e - g + 2$, and in the case $k = 0$, it is $\iota_e - g + 3$. As m - 2 is greater than this index, the Sard-Smale theorem [S] finds a residual set in $\mathcal{J}^m$ of regular values. If J comes from this set, then $\pi_{\mathcal{J}}^{-1}(J)$ is a manifold of differentiablility class m-1 and dimension equal to the index.

Let $d_e = \iota_e - (g - k)$. As argued next, the inverse image of a regular value of $\pi_{\mathcal{J}}$ is empty if $d_e < 0$. This is an automatic consequence of the definition of a regular value in the case when $k > 1$. In the cases $k = 1$ and $k = 0$, the inverse image of a regular value has free group action, the torus $S^1 \times S^1$ if $k = 1$, and the group $Sl(2; \mathbb{C})$ in the case $k = 0$. These actions are free by virtue of the fact that the maps from $\Sigma$ are almost everywhere one-to-one. If $d_e < 0$, then the inverse images in these cases have dimension less than the dimension of the relevant Lie group. As a consequence, these spaces must also be empty when $d_e < 0$.

Let $(\varphi, J, j)$ denote a given element in $\mathcal{P}^m$ and introduce the operator $\mathfrak{D}_*$ that is depicted in (A.6). Identify the tangent space to $(\varphi, J, j)$ with the kernel in $\ker(\mathfrak{D}_*) \times T\mathcal{J}^m|_J$ of the linear part, $(\mathfrak{v}, J') \to f_{0,1}J' + f_{1,0}\mathfrak{v}$, of the right hand side of (A.7). Given the form of $\Psi'$ in (A.8), it follows that $(\varphi, J, j)$ is a critical point of $\pi_{\mathcal{J}}$ if and only if the operator $\mathfrak{D}_*$ has non-trivial cokernel.

*Part 5*: When $j \in V_1$, use $\|j\|_\infty$ to denote its $L^\infty$ norm when viewed as a section of $\mathrm{Hom}(T_0^{1,0}, T_0^{0,1})$. Given $R \geq 1$, introduce as notation $\mathcal{P}^m_R$ to denote the subspace of triples $(\varphi, j, J)$ in $\mathcal{P}^m$ with $\|\varphi\|_* + \|j\|_\infty < R$. This is an open set. In addition, $\mathcal{P}^m_R \subset \mathcal{P}^m_{R+1}$ and $\cup_{R \in \{1,2,\ldots\}} \mathcal{P}^m_R = \mathcal{P}^m$. These sets are the focus of the next lemma.

**Lemma A.4**: *Fix $R \in \{1, 2,\ldots\}$. Let $J \in \mathcal{P}^m$ denote a regular value for the restriction of $\pi_{\mathcal{J}}$ to $\mathcal{P}^m_{R+1}$. Then there is an open neighborhood of $J$ in $\mathcal{P}^m$ whose elements are regular values for the restriction of $\pi_{\mathcal{J}}$ to $\mathcal{P}^m_R$.*

This lemma is proved momentarily.

Given that the set of positive integers is countable, what follows is a consequence of Lemma A.4.

**Proposition A.5**: *The map $\pi_{\mathcal{J}}$ between the Fréchet manifolds $\mathcal{P}$ and $\mathcal{J}$ has a residual set of regular values in $\mathcal{J}$.*

*Proof of Proposition A.5*: For any given $R$, the set of regular values for $\pi_{\mathcal{J}}$'s restriction to $\mathcal{P}^m_{R+1}$ is a residual set, so in particular it is dense. This and Lemma A.4 imply that there is an open and dense set of regular values for $\pi_{\mathcal{J}}$'s restriction to $\mathcal{P}^m_R$. This being the case, there is a $C^m$ open and dense set of regular values for $\pi_{\mathcal{J}}$'s restriction to $\mathcal{P} \cap \mathcal{P}^m_R$. This set is the intersection between $\mathcal{P}$ and the afore-mentioned open dense set in $\mathcal{J}^m$. A $C^m$ open set in $\mathcal{P}$ is, by definiton, a Fréchet open set. A straightforward argument using mollifiers shows that a $C^m$ open and dense set in $\mathcal{P}$ is also open and dense in the $C^\infty$ Fréchet topology. It follows as a consequence that there is a Fréchet open and dense set in $\mathcal{P}$ of regular values for the restriction of $\pi_{\mathcal{J}}$ to $(\mathcal{P} \cap \mathcal{P}^m_R)$. The intersection of the countable set of integer $R$ versions of these open and dense sets is the desired residual set.

*Proof of Lemma A.4*: Suppose that $J \in \mathcal{J}^m$ is a regular value for the restriction of $\pi_{\mathcal{J}}$ to $\mathcal{P}^m_{R+1}$. Let $(\varphi, J, j)$ denote a point in $\pi_{\mathcal{J}}^{-1}(J) \cap \mathcal{P}^m_R$ and let $\mathfrak{D}_*$ denote the operator from $L^2_1(\Sigma, \varphi^*T_{J(1,0)}X) \oplus V$ to $L^2(\Sigma, \varphi^*T_{J(1,0)}X \otimes T_j^{1,0}\Sigma)$ that is defined by the left hand side of (A.6). This operator has trivial cokernel, and so there exists $\delta > 0$ such that if $(\varphi', J', j')$ is in $\mathcal{P}^m$ and $D_*(\varphi, \varphi') + \|J' - J\|_{C^2} + \|j' - j\|_\infty < \delta$, then the $(\varphi', J', j')$ version of $\mathfrak{D}_*$ also has trivial cokernel. Note in this regard that this assertion depends only on the $C^2$ norm of $J' - J$, not the full $C^m$ norm. Also keep in mind that the forgetful map from $C^m$ to $C^2$ is a compact mapping when $n > 1$. Lemma A.2 implies that $\pi_{\mathcal{J}}^{-1}(J) \cap \mathcal{P}^m_R$ has compact closure in $\pi_{\mathcal{J}}^{-1}(\mathcal{P}^m_{R+1})$. Granted these last two points, there exists a positive constant $\delta$ such that the $(\varphi', J', j') \in \mathcal{P}^m_R$ version of $\mathfrak{D}_*$ has trivial cokernel if $J' \in \mathcal{J}^m$ has $C^m$ distance less than $\delta$ from $J$.

### b) The topologies on $\mathcal{M}_{e,k}$ and $\mathcal{P}^m$.

This part starts to explain the relation between $\mathcal{P}$ and the space $\mathcal{M}_{e,k}$. To set the notation, fix for the moment an almost complex structure $J \in \mathcal{J}^m$. As in the case when $J$ is smooth, a closed set $C \subset X$ with finite, non-zero 2-dimensional Hausdorff measure is said to be a $J$-holomorphic subvariety if it has no isolated points; and if the complement of a finite set of points in $C$ is a $m-2$ times differentiable submanifold with $J$-invariant tangent space. Introduce $\mathcal{M}_{e,k}|_J$ to denote the space of irreducible, $J$-holomorphic subvarieties whose model curve has genus $k$. If $(\varphi, J, j)$ is an element in $\mathcal{P}^m$, then $\varphi(\Sigma) \in \mathcal{M}_{e,k}|_J$. Conversely, if $C \in \mathcal{M}_{e,k}|_J$, then there exists a map $\varphi: \Sigma \to X$ with continuous derivatives through order $m-2$ and an element $j \in V_1$ such that $(\varphi, J, j) \in \mathcal{P}^m$, and such that $C = \varphi(\Sigma)$. In particular, $\Sigma$ with its complex structure defined by $j$ is a model curve for $C$.

The space $\mathcal{M}_{e,k}|_J$ has a topology given by (3.2). As is argued next, this topology is essentially equivalent to the topology on $\pi_{\mathcal{J}}^{-1}(J)$ that comes by viewing the latter as a subset of $\mathcal{P}^m$. The proposition that follows makes this precise.

**Proposition A.6**: *Fix a point $(\varphi, J, j) \in \mathcal{P}^m$ and let $C = \varphi(\Sigma)$ denote the corresponding subvariety in $\mathcal{M}_{e,k}|_J$. Given $\varepsilon > 0$, there exists a neighborhood, $\mathcal{N}$, of $(\varphi, J, j)$ in $\mathcal{P}^m$ with the following property: If $(\varphi', J', j') \in \mathcal{N}$, then $\varphi'(\Sigma)$ has $\mathfrak{d}$-distance less than $\varepsilon$ from $C$. Conversely, given a neighborhood, $\mathcal{N}$, of $(\varphi, J, j)$ in $\mathcal{P}^m$, there exists $\varepsilon > 0$ with the following property: If $J' \in \pi_{\mathcal{J}}(\mathcal{N})$ has $C^m$-distance less than $\varepsilon$ from $J$, and if $C' \in \mathcal{M}_{e,k}|_{J'}$ has $\mathfrak{d}$-distance less than $\varepsilon$ from $C$, then there exists $(\varphi', j') \in \mathcal{E}_* \times V_1$ with $\varphi'$ $m-2$ times differentiable, and which is such that $(\varphi', J', j') \in \mathcal{N}$ and $\varphi'(\Sigma) = C'$.*

*Proof of Proposition A.6*: The first assertion follows from the continuity of the evaluation map from $\Sigma \times \mathcal{P}^m$ to $X$, this the map that sends any given pair $(p, \varphi)$ to $\varphi(p)$. The proof of the converse assertion has ten parts.

*Part 1*: The arguments are shorter if $J$ is assumed to be a smooth almost complex structure. To see that it is sufficient to prove only the latter case, remark that if $J$ is only m times differentiable, there are elements in $\mathcal{P}^m$ as close as desired to $(\varphi, J, j)$ with smooth almost complex structures. As explained momentarily, this follows from (A.7) and (A.8). Meanwhile, if $(\varphi'', J'', j'')$ is sufficiently close in $\mathcal{P}^m$ to $(\varphi, J, j)$, then the subvariety $C'' = \varphi''(\Sigma)$ will lie everywhere very close to C. This understood, it enough to prove the converse statement of the proposition with $(\varphi, J, j)$ replaced by an element in $\mathcal{P}^m$ very close to it but with smooth almost complex structure. The fact that there are elements in $\mathcal{P}^m$ with smooth almost complex structure near any given element $(\varphi, J, j)$ can be seen as follows: Identify a neighborhood of J in $\mathcal{J}^m$ with a neighborhood of the origin in the space of m-times differentiable sections of $\mathrm{Hom}(T_J^{1,0}X; T_J^{0,1})$. Let $\mathcal{C}$ denote the affine subspace in this space of sections that represent smooth almost complex structures. This affine space is dense in the space of all m-times differentiable sections. Restrict the map $f$ in (A.7) to the subspace of its domain consisting of points of the form $(0, J')$ with $J' \in \mathcal{C}$. Given that $\mathcal{C}$ is dense, and given that $f_{0,1}$ in (A.7) is surjective, it follows that there is a finite dimensional linear subspace $\mathcal{C}_0 \subset \mathcal{C}$ of dimension equal to that of $\mathrm{coker}(\mathfrak{D}_*)$, whose origin is as close as desired to the zero section of $\mathrm{Hom}(T_J^{1,0}X; T_J^{0,1})$ and with the following property: The map $f_{0,1}$ restricts to $\mathcal{C}_0$ as an isomorphism with inverse bounded by $c_0$ with $c_0$ here independent of the distance between $\mathcal{C}_0$ and the zero section of $\mathrm{Hom}(T_J^{1,0}X; T_J^{0,1})$ if this distance is less that $c_0^{-1}$. Granted that such is the case, it follows by degree theory that there are elements of the form $(0, J') \in f^{-1}(0)$ with $J' \in \mathcal{C}_0$.

The remaining parts of the proof assume implicitly that J is smooth.

*Part 2*: Fix $\delta > 0$ but very small; in particular much less than the distance between any two singular points of C. Let $B_\delta \subset X$ denote the union of balls of radius $\delta$ about these singular points. The intersection of C with the interior of $X - B_\delta$ is a properly embedded J-holomorphic submanifold. If $\delta$ is suitably generic, then C will have transversal intersection with the boundaries of the closures of $B_\delta$ and also $B_{k\delta}$ for $k = 2, 3$, and 4. Moreover, C's intersection with $B_{k\delta} - B_{(k-1)\delta}$ for each such n will be a disjoint union of embedded annuli. Assume that such is the case.

Let $C_\delta$ denote $C \cap (X - B_\delta)$ and likewise define $C_{k\delta}$. Note that the integral of $\omega$ over $C_{4\delta}$ is less than $[\omega] \cdot e$ by at most $c_0 \delta^2$. Take $\delta$ so that this is less than $10^{-100} [\omega] \cdot e$. Note also the following: Given $\delta_0 \in (0, \delta_1)$, then $\varphi^{-1}(B_{4\delta})$ is contained in the union of disks of radius $\delta_1$ about the $\varphi$-inverse images of the singular points of C if $\delta$ is sufficiently small. Let $D(\delta_0) \subset \Sigma$ denote this union of disks.

Let $N \to C_\delta$ denote the normal bundle to $C_\delta$. There exists $\varepsilon_1 > 0$ such that if $J'$ has $C^m$ distance less than $\varepsilon_1$ from J, then the constructions in Section 5d in [T3] can be used to find $\delta' > 0$ and an embedding, $\exp_C$, from the radius $\delta'$ subdisk bundle in N onto a neighborhood of $C_\delta$ in $X - U_{\delta/2}$ with the following properties:

- *It restricts to the zero section as the identity map.*

- *The image of each fiber disk is J´-pseudoholomorphic.*

(A.9)

Let $N_{\delta'} \subset N$ denote this subdisk bundle, and use $\exp_C$ to identify it with its image in X.

*Part 3*: Suppose that $C´ \subset X$ is J´-pseudoholomorphic. If C´ intersects $X-B_\delta$ so as to lie entirely in $N_\delta$, then it will intersect each fiber disk with each intersection having positive intersection number. Moreover, the number of such intersections is the same for each such disk. It follows as a consequence that C´ has precisely one intersection with each fiber disk. Indeed, if such isn't the case, then the integral of ω over $C´ \cap B_\delta$ will be larger than [ω]·e and this is not possible since ω is positive on each tangent plane of C´. Granted that C´ intersects each fiber disk once, it intersect $N_{\delta'}|_{C_{2\delta}}$ as the graph of a section of $N_{\delta'}$. Let η denote this section. Since φ maps $\Sigma - \varphi^{-1}(B_{2\delta})$ in a 1-1 fashion to $C-B_\delta$, the composition $\eta \circ \varphi$ is a section over $\Sigma - B_{2\delta}$ of φ*N with pointwise norm less than δ. Moreover, it must obey the inhomogeneous version of (2.14); this an equation of the form

$$D_C \eta + \mathfrak{r}_1 \cdot \partial \eta + \mathfrak{r}_2 = \mathfrak{e}$$

(A.10)

where $\mathfrak{r}_1$ and $\mathfrak{r}_2$ are as in (2.14) and where the norm of $\mathfrak{e}$ is bounded by $c_0 \varepsilon$ if J´ is ε-close to J in the $C^m$ topology.

Granted what was just said, the assertion that follows can be proved using the elliptic regularity arguments of the sort found in Chapter 6 of [Mo]: If C´ has ∂-distance $\varepsilon < 10^{-4}\delta´$ from C, then η defines an m -2 times differentiable section of φ*N over $\Sigma - B_{3\delta}$ with $C^{m-2}$ norm bounded by $c_\delta \varepsilon$. Here, and in what follows, $c_\delta > 1$ is a constant that depends on δ and also C, but nothing else of relevance. Its value can be assumed to increase between consecutive appearances.

The fact that $C´ \cap (X-B_{3\delta})$ is the image of $\eta \circ \psi$ implies that the singularities of C´ are contained in $B_{3\delta}$. This understood, the pull back via $\eta \circ \psi$ of the restriction of J to $T(C´\cap (X-B_{3\delta}))$ defines a complex structure, $j´_\delta$, on $\Sigma - \varphi^{-1}(B_{3\delta})$. Given that the $C^{m-2}$ norm of η is bounded by $c_\delta \varepsilon$, this complex structure differs from j by an endomorphism of $T(\Sigma - \varphi^{-1}(B_{3\delta}))$ whose $C^{m-3}$ norm is bounded by $c_\delta \varepsilon$.

*Part 4*: The subvariety C´ can, in any event, be written as φ´(Σ) where (φ´, J´, j´) $\in \mathcal{P}^m$. Let $\Sigma´_\delta$ denote $(\varphi´)^{-1}(C´ \cap (X-B_{3\delta/2}))$. The composition $(\eta \circ \varphi)^{-1} \circ \varphi´$ is an m-2 times differentiable homeomorphism between $\Sigma´_\delta$ and its image inside $\Sigma - \varphi^{-1}(B_{3\delta})$. Let $\psi_\delta$ denote this map. If $\varepsilon < c_\delta^{-1}$, then the image of $\psi_\delta$ contains $\Sigma - \varphi^{-1}(B_{4\delta})$ and thus $\Sigma - D(\delta_0)$. Note that the complement in Σ of $\psi_\delta^{-1}(\Sigma - D(\delta_0))$ must be a disjoint union of embedded disks that contains the points mapped by φ´ to the singularities of C´. Moreover, each disk in this union is mapped by φ´ to a component of $B_{4\delta}$. Let D´($\delta_0$) denote this union of disks. The map $\psi_\delta$ sets up a 1-1 correspondence between the disks that comprise $D(\delta_0)$ and those that comprise D´($\delta_0$) by identifying the boundary circles of these two sets. Keep in mind that the maps φ and φ´ send partnered components of D($\delta_0$) and D´($\delta_0$) to the same component of $B_{4\delta}$.

The task ahead is to extend the map η∘φ over each component of D($\delta_0$) in a controlled fashion. This is done in this part and in Parts 5-7 to come. The following notation is used here and in Parts 5-7. The symbol D denotes a component of D($\delta_0$), and D´ ⊂ D´($\delta_0$) the partner of D.

Suppose first that φ|$_D$ is an embedding. In this case, the exponential map $\exp_C$ can be extended across D and likewise the section η. This done, then φ´(D´) can be written as $\exp_C$(η∘φ). Granted this, no generality is lost by assuming in what follows that the differential of φ vanishes at the center of each component of D($\delta_0$).

It is necessary to find a new parametrization of φ(D) and likewise φ´(D´). The desired parametrization of D is described by the next lemma. This parametrization involves a new choice of coordinates for a neighborhood of the singular point of φ(D).

**Lemma A.7**: *Given J and C, there exist $\kappa > 10^4$ and an integer n > 1 with the following properties: Suppose that $\Delta < \kappa^{-1}$ and that J´ has $C^m$ distance less than $\Delta$ from J. There are m-1 times differentiable, complex coordinates (z, w) centered at φ(0) such that*
- *The coordinates are valid where both $|z| < \kappa^{-1}$ and $|w| < \kappa^{-1}$.*
- *The w = 0 disk is J´-pseudoholomorphic, as are the constant z disks.*
- *The bundle $T_{J´}^{(1,0)}$ is spanned by the 1-forms*

$$dz - \alpha \, d\bar{z} \quad and \quad dw - \gamma \, d\bar{z}$$

  *where α and γ are m-1 times differentiable functions that vanish on the w = 0 locus and otherwise obey $|\alpha| \le c_0^{-1}$, $|\gamma| \le c_0^{-1}$, $|d\alpha| \le c_0$ and $|d\gamma| \le c_0$.*
- *The subvariety φ(D) intersects the $|z| < \kappa^{-2}$ part of this coordinate chart as an m-1 times differentiable map from a disk in $\mathbb{C}$ about the origin of radius bounded by $\kappa^{-2/n}$ that has the form*

$$u \to (z = u^n(1 + \mathfrak{r}_z), w = u^n \mathfrak{r}_w)$$

  *where $|\mathfrak{r}_z| \le \Delta$ and $|d\mathfrak{r}_z| \le \kappa\Delta$; while $|\mathfrak{r}_w| \le \Delta$ and $|d\mathfrak{r}_w| \le \kappa\Delta$.*

*Proof of Lemma A.7*: Coordinates of this sort are obtained using the arguments in Section 5d of [T3]. To elaborate, let p denote the singular point in φ(D). Fix an orthonormal frame for $T_{J(1,0)}|_p$ so as to identify the space of 1-dimensional subspaces through the origin in $T_{J(1,0)}|_p$ with $\mathbb{CP}^1$. What is said in Section 5d of [T3] can used to associate to each $\theta \in \mathbb{CP}^1$ an embedded, pseudoholomorphic disk with center p and tangent plane $\theta$ at p. Moreover, each such disk will intersect the closed, radius $c_0^{-1}$ ball centered at p as disk and intersect the boundary transverally. Finally, any two $\theta \neq \theta'$ disks will intersect only at p. Each such disk has a finite number of intersections with $\varphi(D)$, each with positive intersection number. Use n to denote the minimum of these intersection numbers. Fix $\theta \in \mathbb{CP}^1$ whose disk has intersection number n.

Use what is said in Section 5d of [T3] to find complex coordinates, (z, w) centered on p, defined for $|z| < c_0^{-1}$ and $|w| < c_0^{-1}$ and such that the following is true: First, the bundle $T_J^{1,0}$ is spanned by $dz - a d\bar{z}$ and $dw - g d\bar{z}$ where $a$ and $g$ vanish on the $w = 0$ locus, and are such that their derivatives are bounded in norm by $c_0$. Second, $dw|_p$ vanishes on the complex line $\theta$. Note that it follows from the form of this basis for $T_J^{1,0}$ that the $w = 0$ disk is J-pseudoholomorphic, as are the constant z disks. Finally, each such constant z disk has intersection number n with $\varphi(D)$ and each of these intersections will have intersection number 1.

Granted the preceding, and granted that $dz - a d\bar{z}$ is a form of type (1, 0) it follows that the function $z_\varphi = z(\varphi(\cdot))$ must obey the equation

$$\bar{\partial} z_\varphi - a \bar{\partial} \bar{z}_\varphi = 0 \ .$$

(A.11)

Because $a = 0$ where $w = 0$, there exists a holomorphic coordinate, u, that is defined on a neighborhood of the center of D with $u = 0$ this same center point, and is such that $z_\varphi(u) = u^n(1 + \mathfrak{r}_z)$ where $|\mathfrak{r}_z| \le c_0^{-1}|u|$ and $|d\mathfrak{r}_z| \le c_0$. The proof that this is the case is an exercise using Taylor's expansion given that (A.11) implies via Aronzajn's unique continuation principle [A] that the function $u \to z_\varphi(u)$ can not vanish to infinite order at any point.

Meanwhile, let $w_\varphi = w(\varphi(\cdot))$. Then $|w_\varphi|(u)$ must be bounded by $c_0|u|^n$ because of the minimality of n. Furthermore, given that $dw - g d\bar{z}$ is of type (1, 0), the function $w_\varphi$ obeys the analog of (A.6) with $a$ replaced by $g$ and with only the left most $z_\varphi$ replaced by $w_\varphi$. Given that $\gamma$ also vanishes at the origin, a Taylor's theorem argument shows that $w = c(u^{n+k} + \mathfrak{r}_w)$ where $c \in \mathbb{C}-0$, k is a non-negative integer and $|\mathfrak{r}_w| \le c_0|u|^{n+k+1}$ and $|d\mathfrak{r}_w| \le c_0|u|^{n+k}$. Note that k must be finite or else $\varphi$ will map D as a multiple cover onto the $w = 0$ locus, and this is forbidden because $\varphi$ is assumed to be almost everywhere 1-1. (Keep in mind here that distinct J-holomorphic subvarieties can not be tangent to infinite order.) If $k = 0$ it is necessary to take a different choice for $\theta$. To explain, the linear change of coordinates to $z' = z$ and $w' = w - cz$ makes $w' = c'(u^{n+k} + e')$ where $k > 0$ and where $c' \neq 0$ and $|e'| \le c_0|u|^{n+k+1}$. Even so, the $w' = 0$ disk will no longer be J-holomorphic. However, the constructions in Section 5d of [T3] can be used to construct a different version of $w'$, this given by $w' = w - cz + \mathfrak{v}(z)$ with $|\mathfrak{v}| \leq c_0 |z|^2$ so that the $w' = 0$ disk is J-holomorphic.

Given that $J'$ is sufficiently close in the $C^{2n}$ norm to $J$, the constructions in Section 5d of [T3] can be used to perturb the coordinates to get coordinates $(z, w)$ which are as described by the lemma.

*Part 5*: This part begins the task of constructing the desired coordinates for $\varphi'(D')$. The basic construction is, to some extent, similar to what was done in the previous part of the proof for $\varphi(D)$. There is, however, the added complication that $\varphi'(D')$ may have more than one singular point.

To start, let $D_* \subset \mathbb{C}$ denote a disk about the origin that is contained in the disk from the fourth bullet of Lemma A.7. The latter disk is where the coordinate u from this fourth bullet is defined. Use $\rho$ in what follows to denote the radius of $D_*$. No generality is lost by assuming $\rho < 10^{-4} \delta_0$ and $\rho < \delta$. The freedom to choose $\rho$ conveniently and in particular small will be exploited in the subsequent discussions. These discussions use $\kappa_\rho$ to denote a constant that is greater than $10^4$ and that depends on $\rho$, C and J, but not on $J'$ nor $C'$. The value of $\kappa_\rho$ can be assumed to increase between subsequent appearances. On the subject of parameters, the discussion that follows will also involve a parameter r, this with values in $(0, \frac{1}{8} \rho)$. The discussion uses $\varepsilon_r$ to denote a positive constant that depends on r, C and J, but not on either $C'$ or $J'$. This constant $\varepsilon_r$ can be assumed to *decrease* between successive appearances.

In what follows, $D_*$ is identified with its image in D via the coordinate map. Let $\lambda: D_* \to \mathbb{C}^2$ denote the map depicted in this same fourth bullet of Lemma A.7, and use $z_\lambda$: $D_* \to \mathbb{C}$ to denote the composition of first $\lambda$ and then the coordinate function z. Introduce $U \subset \mathbb{C}$ to denote $z_\lambda(D_*)$. The fourth bullet of Lemma A.7 asserts that U contains the disk about the origin of radius $\rho^n(1 - c^{-1})$ and is contained in the concentric disk of radius $\rho^n(1 + c^{-1})$ where $c > 10^4$.

Consider now $(z \circ \varphi')^{-1}(U)$. If $\varepsilon < c_0^{-1}$, then $(z \circ \varphi')^{-1}(U)$ is an embedded disk in $D'$ with m-3 times differentiable boundary. This follows from the fact that the section $\eta$ is defined on $D - D_*$ if $\varepsilon < \kappa_\rho^{-1}$. The Riemann mapping theorem finds a m-3 times differentiable homeomorphism $\sigma: D_* \to (z \circ \varphi')^{-1}(U)$ so that the composition, $\varphi' \circ \sigma$, is a $J'$-holomorphic map from $D_*$ to $\varphi'(D')$. Let $\lambda'$ denote the latter map and introduce $z_{\lambda'}$ to denote $z \circ \lambda'$. This $z_{\lambda'}$, like $z_\lambda$, is also a function on $D_*$ with image equal to U. By composing $\lambda'$ with a suitable Möbius transformation of $D_*$, a version of $z_{\lambda'}$ can be had that maps the origin in $D_*$ to the origin in U. This condition on $\lambda'$ and $z_{\lambda'}$ is assumed in what follows as it plays an important role.

To say more about $z_{\lambda'}$, note that given $r \in (0, 8^{-1}\rho)$, there exists $\varepsilon_r$ with the following significance: If $\varepsilon < \varepsilon_r$, then the critical values of $z_{\lambda'}$ must lie in the disk of radius $r^n$ centered on the origin. This again follows from the fact that the coordinate z on D´ where $|z| > r$ looks very much like that on D for the reason that this part of D´ is parametrized by $\eta \circ \varphi$ if $\varepsilon$ is small, which is to say less than $\varepsilon_r$. For this same reason, the function $z_{\lambda'}$ on the annulus in $D_*$ where $|z_{\lambda'}| > r$ has a single valued n'th root which gives a bonafide coordinate on this annulus. The latter coordinate plays a role also in what follows.

Lemma A.7 also supplies the coordinate w for its neighborhood of $\varphi(0)$. Introduce $w_{\lambda'}: D_* \to \mathbb{C}$ to denote the composition $w \circ \lambda'$. Use $\alpha_{\lambda'}$ and $\gamma_{\lambda'}$ to denote the respective $\mathbb{C}$-valued functions on $D_*$ that come by composing the respective functions $\alpha$ and $\gamma$ from the third point of Lemma A.7 with the map $\lambda'$. Given the basis for $T_{J'}^{1,0}$ in from the third bullet of Lemma A.7, it follows that the pair $(z_{\lambda'}, w_{\lambda'})$ must obey the equations

$$\bar\partial z_{\lambda'} - \alpha_{\lambda'} \bar\partial \bar z_{\lambda'} = 0 \quad and \quad \bar\partial w_{\lambda'} - \gamma_{\lambda'} \bar\partial \bar z_{\lambda'} = 0$$

(A.12)

on $D_*$. Note that both $\alpha_{\lambda'}$ and $\gamma_{\lambda'}$ are bounded in absolute value by $c_0 |w_{\lambda'}|$.

This last equation leads to the following lemma concerning the size of the derivatives of $z_{\lambda'}$ and $w_{\lambda'}$.

**Lemma A.8**: *There exists $\kappa > 1$ with the following significance: Suppose that J´ has $C^{2n}$ distance less than $\kappa^{-1}$ from J and that $\rho < \kappa^{-2}$. Suppose in addition that $\varepsilon < \rho^n$. Then the functions $z_{\lambda'}$ and $w_{\lambda'}$ on the radius $\frac{7}{8}\rho$ disk in $D_*$ centered on the origin obey*

$$|z_{\lambda'}| + |w_{\lambda'}| + \rho(|dz_{\lambda'}| + |dw_{\lambda'}|) \leq \kappa \rho^n .$$

*Proof of Lemma A.8*: The bound on the absolute values of $z_{\lambda'}$ and $w_{\lambda'}$ follow from the fact that $z_\lambda$ and $w_\lambda$ obey similar bounds. To see about the bounds on the higher derivatives, reintroduce the functions $\alpha$ and $\gamma$ from the third bullet in Lemma A.7. The functions $\alpha_{\lambda'}$ and $\gamma_{\lambda'}$ are obtained from the latter by evaluating them at $z = z_{\lambda'}$ and $w = w_{\lambda'}$. This understood, and given the bounds on the norms of $z_{\lambda'}$ and $w_{\lambda'}$, the bounds on the derivatives follow from (A.12) using the elliptic bootstrap arguments from Chapter 6 in [Mo].

*Part 6*: This part constitutes a digression for two results about quasi-conformal maps that will be used in conjunction with (A.12) to say more about $z_\lambda$. The first lemma refers to the *hyperbolic radius* of a disk in $D_*$. This is the radius as measured by the hyperbolic metric on the disk $D_*$, this the metric obtained from the unit radius disk in $\mathbb{C}$ with line element $(1 - |u|^2)^{-1} |du|$ via pull-back by a holomorphic diffeomorphism from $D_*$ that maps the origin to the origin. The hyperbolic radius has the advantage of being invariant under the action of the group, Sl(2; $\mathbb{R}$), of holomorphic diffeomorphisms of $D_*$.

**Lemma A.9**: *Given* $d \in (0, 1)$, *there exists a constant* $\kappa > 1$ *whose significance is given in what follows. Fix* $r < \frac{1}{8}\rho$ *and let* $U_r \subset U$ *denote the disk of radius* $r^n$ *centered on the origin. Suppose that* $x: D_* \to U$ *is a twice differentiable map with the following three properties. First, $x$ restricts to the boundary of $D_*$ as a degree n-1 covering map of the boundary of U. Second, x is such that* $|\bar{\partial} x| \leq (1 - d)|\partial x|$ *at all points in $D_*$. Finally, the critical values of x lie in $U_r$. Then $x^{-1}(U_r)$ is contained in a disk with hyperbolic radius less than* $\kappa (r/\rho)^{d/d-2} \rho$.

*Proof of Lemma A.9*: Let $A_r = x^{-1}(U - U_r)$. This is a topological cylinder with $C^1$ inner boundary and outer boundary that of $D_*$. Meanwhile, the Riemann mapping theorem finds $s \in (0, 1)$ with the following significance: Let $\hat{U} \subset \mathbb{C}$ denote the annulus centered at the origin with outer radius 1 and inner radius s. The Riemann mapping theorem finds a holomorphic, degree n covering map, $g: \hat{U} \to U - U_r$. Note that $s \sim \rho^{-1} r$ because the boundary of U is nearly a round circle of radius $\rho^n$. It follows from the path lifting property of covering spaces that there is a lift, $\hat{x}: A_r \to \hat{U}$ such that $x = g \circ \hat{x}$. The map $\hat{x}$ is 1-1 and so a uniformly, quasi-conformal homeomorphism from the domain $A_r$ to the domain $\hat{U}$. Indeed, given that g is holomorphic, the chain rule finds $|\bar{\partial} \hat{x}| \leq (1 - d)|\partial \hat{x}|$.

The Riemann mapping theorem also finds $s' \in (0, 1)$ and a holomorphic homeomorphism to $A_r$ from the cylinder $\hat{A} \subset \mathbb{C}$ with outer radius 1 and inner radius $s'$. It is a consequence of Theorem 7.1 in [LV] that

$$\frac{d}{2-d} |\log(s)| \leq |\log(s')| \leq \frac{2-d}{d} |\log(s)|.$$

(A.13)

The number $\frac{1}{2\pi} |\log(s')|$ said to be the *modulus* of $A_r$. This understood, and given that $r/\rho$ << 1, then what is claimed by the lemma follows from the upper bound (A.13) using the inequalities given in Theorem 2.4 in [McM] that relate the hyperbolic diameter to the modulus.

The next lemma states the second of what is needed from quasi-conformal mapping theory. The lemma uses $\underline{D}_*$ to denote the closure in $\mathbb{C}$ of the disk $D_*$.

**Lemma A.10**: *Fix* $\varepsilon > 0$ *and* $0 < \rho' < \frac{1}{8}\rho$. *There exists constants* $\kappa_1, \kappa_2 > 100$ *with the following significance: Suppose that* $A \subset \underline{D}_*$ *is a closed annulus with outer boundary* *that of $\underline{D}_*$ and with $C^1$ inner where $|u| < \frac{1}{2}\rho'$. Suppose in addition that $u: A \to \underline{D}_*$ is a twice differentiable, 1-1 embedding whose image is the annulus in $\underline{D}_*$ where $\kappa_1^{-1}\rho' \le |u|$, and is such that $|\bar\partial u| < \kappa_2^{-1}|\partial u|$. This function $u$ must restrict to the annulus in $\underline{D}_*$ where $\rho' \le |u| \le (1 - \varepsilon)\rho$ so as to have the form $u(u) = \theta u + \tau$ with $|\theta| = 1$ and with $|\tau| < \varepsilon$.*

***Proof of Lemma A.10***: Suppose that no such constants exist. There would in this case exist numbers $\varepsilon$ and $\rho'$ as in the statement of the lemma, plus sequences $\{A_k \subset \underline{D}_*\}_{k=1,2,...}$ and $\{u_k\}_{k=1,2,...}$ of the following sort: First, each $A_k$ is an annulus of the sort under consideration. Meanwhile, each $u_k: A_k \to \underline{D}_*$ is a map as described in the statement of the lemma, and in particular one that sends $A_k$ to the annulus in $\underline{D}_*$ where $k^{-1}\rho' \le |u|$ and is such that $|\bar\partial u_k| \le k^{-1}|\partial u_k|$. It follows from Lemma A.9 and what is said in Chapter II.5 of [LV] that the sequence $\{u_k\}$ will have a subsequence (hence renumbered consecutively from 1) that converges uniformly in the $C^0$ and $L^2_{1;loc}$ topologies on compact subsets of $D_* - \{0\}$ to a 1-1 holomorphic self map of $D_*$ that sends 0 to 0. The latter must have the form $u \to \theta u$ for $\theta \in \mathbb{C}$ having norm 1. It follows as a consequence that $|u_k - \theta u| < \varepsilon$ for all k sufficiently large on the annulus where $\rho' < |u| < (1 - \varepsilon)\rho$.

*Part 7*: Since $z_{\lambda'}$ obeys the equation on the left hand side of (A.12) with $|\alpha_{\lambda'}| \ll 1$, when $\rho < c_\delta^{-1}$ and $\varepsilon < \varepsilon_r$, Lemma A.9 can be applied to see that the critical points of $z_{\lambda'}$ lie inside a disk in $D_*$ centered at the origin with radius no greater than $c_0 (r/\rho)^{3/4}\rho$. This said, fix $r < c_0^{-1}\rho$ to guarantee that these critical points lie inside the disk $D_{**} \subset D_*$ with radius $10^{-4}\rho$. View $z_{\lambda'}$ as a $\mathbb{C}$-valued function on $D_* - D_{**}$, and it follows that the latter has an n'th root that gives a bonafide m - 2 times differentiable coordinate for $D_* - D_*$. Use $u_{\lambda'}$ to denote this coordinate. It is a consequence of the left-hand equation in (A.12) that the map $u \to u_{\lambda'}(u)$ obeys

$$\bar\partial u_{\lambda'} - \alpha_{\lambda'} \left(\frac{\bar u_{\lambda'}}{u_{\lambda'}}\right)^{n-1} \bar\partial \bar u_{\lambda'} = 0$$

(A.14)

This equation is used in conjunction with Lemmas A.8-A.10 to prove

**Lemma A.11**: *Given $\delta$, $\rho$, and also $\varepsilon_1 > 0$, there exists $\kappa_2 > \kappa_1 > 100$ with the following signicance: Suppose that $J'$ has $C^m$ distance less than $\kappa_1^{-1}$ from $J$, that $\rho < \kappa_1^{-1}$, and that $C'$ has $\mathfrak{d}$-distance less than $\kappa_2^{-1}$ from $C$. Then $u_{\lambda'}(u)$ has the form $u_{\lambda'}(u) = \theta u + \tau$ on the annulus where $10^{-4}\rho < |u| < \frac{7}{8}\rho$ with $\theta \in \mathbb{C}$ having norm 1 and with $|\tau| < \varepsilon_1$ and $|d\tau| < \varepsilon_1$.*

***Proof of Lemma A.11***:  It is a consequence of what is said in Lemma A.8 that the function $\alpha_{\lambda'}$ in (A.9) has absolute value no greater than $\rho^n$ when $\varepsilon < \rho^n$.  As noted earlier, given $r > 0$, there exists $\varepsilon_r > 0$ such that if $\varepsilon < \varepsilon_r$, then $z_{\lambda'}$ has none of its critical values where $r^n < |z_{\lambda'}| \le \rho^n$.  This and what was just said about $\alpha_{\lambda'}$ imply via Lemma A.9 that $z_{\lambda'}$ has no critical points on the part of $D_*$ where $c_0(r/\rho)^{1/c_0} \rho < |u| \le \rho$.  Given this last fact and the aforementioned bound on $|\alpha_{\lambda'}|$, Lemma A.10 supplies constants $c_1, c_2 > 1$ such that if $\rho < c_1^{-1}$ and $\varepsilon < c_2^{-1}$, then $u_{\lambda'}$ can be written on the annulus where $10^{-5}\rho < |u| < \frac{15}{16}\rho$ as $u_{\lambda'}(u) = \theta u + \mathfrak{r}$ where $\theta \in \mathbb{C}$ has norm 1 and $|\mathfrak{r}| < \varepsilon_1$.  This understood, it remains yet to bound $|d\mathfrak{r}_1|$.  To do so, use (A.14) to see that $\mathfrak{r}$ obeys an equation of the form

$$\overline{\partial}\,\mathfrak{r} - \alpha_{\lambda'}\,\Big(\tfrac{\overline{z}_{\lambda'}}{z_{\lambda'}}\Big)^{1-1/n}\,\overline{\partial}\,\overline{\mathfrak{r}} - \alpha_{\lambda'}\,\Big(\tfrac{\overline{z}_{\lambda'}}{z_{\lambda'}}\Big)^{1-1/n}\,\overline{\theta} = 0 \,.$$

(A.15)

Let $\alpha$ denote the function that appears in the third bullet of Lemma A.7.  Given that $\alpha_{\lambda'} = \alpha(z_{\lambda'}, w_{\lambda'})$, given the first derivative bounds from Lemma A.8, and given that $|\mathfrak{r}| < \varepsilon_1$, the elliptic boot-strapping techniques from Chapter 6 in [Mo] applied to (A.15) give the desired bound on $|d\mathfrak{r}|$.

    *Part 8*:  Compose that map $u \to (z_{\lambda'}, w_{\lambda'})$ with a rigid rotation of the disk so as to make the constant $\theta$ from Lemma A.11 equal to 1.  Use the same notation, $z_{\lambda'}$ and $w_{\lambda'}$ for the result of this composition.

    If $\rho < c_0^{-1}$, if $r < c_0^{-1}\rho$ then there is a constant, $\varepsilon_r > 0$ such that if $\varepsilon < \varepsilon_r$, two parametrizations of the $2^n 10^{-4n}\rho^n < |z| < (\tfrac{7}{8})^n$ part of $D'$ are available.  The first is via the map

$$u \to (z = z_\lambda(u) + \eta_z(u),\ w = w_\lambda(u) + \eta_w(u))\,.$$

(A.16)

where $\eta_z$, $\eta_w$ and their derivatives to order $m - 3$ are bounded by $c_\rho \varepsilon$.  The second is via the map

$$u \to (z = z_{\lambda'}(u),\ w = w_{\lambda'}(u))$$

(A.17)

The latter has the advantage that it extends over the whole of $D_*$.  The plan is to use this extension to extend (A.16) over $D_*$ too.

    The first step in this task uses Lemmas A.7 and A.11 to bound the distance between their respective images.  In particular, these lemmas assert that the z-coordinates of the images by these respective maps can be written as

$$u \to z_\lambda(u) + \eta_z(u) = u^n + \mathfrak{e} \quad and \quad u \to z_{\lambda'}(u) = u^n + \mathfrak{e}'$$

(A.18)

where $|\mathfrak{e}| + |\mathfrak{e}'|$ can be assumed bounded by $10^{-10n}\rho^n$, and $|d\mathfrak{e}| + |d\mathfrak{e}'|$ can be assumed to be bounded by $10^{-10n}\rho^{n-1}$ if $\rho < c_\delta^{-1}$ and $\varepsilon < c_0^{-1}\rho^n$.

Granted (A.18), the parametrization of the $10^{-3n} < |z| < (\tfrac{3}{4})^n$ part of $D'$ given by the map $u \to (z_{\lambda'}(u), w_{\lambda'}(u))$ can be composed by a rigid rotation, $u \to \lambda u$ with $\lambda^n = 1$ so that

$$z_\lambda(u) + \eta_z(u) = z_{\lambda'}(v(u)) \quad and \quad w_\lambda(u) + \eta_w(u) = w_{\lambda'}(v(u))$$
(A.19)

where $u \to v(u)$ is a map that has the form $v(u) = u + \mathfrak{z}(u)$ with $|\mathfrak{z}| < 10^{-10}\rho$ and $|d\mathfrak{z}| < 10^{-10}$. As a consequence, the parametrization given by (A.16) can be extended by writing it as

$$u \to (z = z_\lambda(u + \beta_\rho \mathfrak{z}(u)), w = w_\lambda(u + \beta_\rho \mathfrak{z}(u)))$$
(A.20)

where $u \to \beta_\rho(u)$ is given by $\beta(\rho^{-1}|u|)$ where $\beta: [0, 1] \to [0, 1]$ is a smooth function that equals 1 on $[\tfrac{1}{2}, 1]$, equals 0 on $[0, \tfrac{1}{4}]$ and is such that $|d\beta| < 8$.

Use $\varphi^\diamond: \Sigma \to C'$ to denote the parametrization that results by constructing just such an extension for each component of $D(\delta_0)$.

*Part 9*: It is a consequence of what is said in Part 3 that the restriction of the map $\varphi^\diamond$ to $\Sigma - (\cup_{D \subset D(\delta_0)} D_*)$ has $L_*$ distance at most $\kappa_\rho \varepsilon$ from $\varphi$'s restriction to this same domain if $\varepsilon < \kappa_\rho^{-2}$. Meanwhile, if $\rho < c_\delta^{-1}$ and if $\varepsilon < \kappa_\rho^{-2}$, it follows from Lemmas A.7 and A.8 that the restriction of $\varphi^\diamond$ to where $|u| \leq 2\rho$ in any disk $D \subset D(\delta_0)$ has $L_*$ distance at most $c_0 \rho^3$ from $\varphi$'s restriction to this same domain. It follows as a consequence that $\varphi^\diamond$ has $L_*$ distance bounded by $c_0 \rho^3$ from $\varphi$ on the whole of $\Sigma$ if $\rho < c_\delta^{-1}$ and if $\varepsilon < \kappa_\rho^{-1}$.

It remains yet to investigate the distance from the complex structure, $j$, to the almost complex structure on $\Sigma$ that makes $\varphi^\diamond$ a $J'$-holomorphic map. Let $j^\diamond$ denote the latter. View both $j$ and $j^\diamond$ as sections over $\Sigma$ of $End(T\Sigma)$. Viewed this way, Part 3 says that the $C^{m-2}$ norm of $j - j^\diamond$ is bounded by $\kappa_\rho^{-1}\varepsilon$ on $\Sigma - (\cup_{D \subset D(\delta_0)} D_*)$.

Now let $D \subset D(\delta_0)$ so as to consider $j - j^\diamond$ on the $|u| \leq 2\rho$ part of $D$. Consider first the behavior of $j$. The fourth item in Lemma A.7 implies the following: Given $\Delta > 0$, there exists $R_\Delta > 1$ such that if $J'$ has $C^{2m}$ distance less than $1/R_\Delta$ from $J$, then

$$|J(\varphi_* \tfrac{\partial}{\partial u}) - \varphi_* \tfrac{\partial}{\partial u}|/|\varphi_* \tfrac{\partial}{\partial u}| \leq \Delta.$$
(A.21)

Meanwhile, (A.12) and Lemma A.8 imply that there exists $R_\Delta' > 1$ such that if $\rho < 1/R_\Delta'$, then (A.21) also holds with $J'$ replacing $J$ and with $\varphi^\diamond$ replacing $\varphi$. These two versions of (A.21) with what was said in the preceding paragraph imply that $|j - j^\diamond| \leq c_0 \Delta$.

*Part 10*: Introduce the reference almost complex structure $j_0$ from Part 2 of Section a) of this Appendix and the corresponding decomposition $T\Sigma_\mathbb{C} = T_0^{1,0}\Sigma \oplus T_0^{0,1}\Sigma$ into forms of type (1, 0) and (0, 1) as defined by $j_0$. This done, view j and also $j^\diamond$ as homomorphism from $T_0^{1,0}\Sigma$ to $T_0^{0,1}\Sigma$ with everywhere norm less than 1. The homomorphism j comes from the disk $V_1$, but this is not necessarily the case for $j^\diamond$. The next lemma is used to deal with this event. This lemma refers to the $D_*$ distance between self maps of $\Sigma$. This distance is defined as in (A.2) using $\Sigma$ in lieu of X.

**Lemma A.12**: *Given an almost complex structure, j, from $V_1$, and given $\tau > 0$, there exists $\kappa_\tau > 1$ such that the following is true: If $j^\#$ is an m-3 times differentiable almost complex structure on $\Sigma$ with $|j - j^\#| < \kappa_\tau^{-1}$, then there exists a homeomorphism of $\Sigma$ which is m-3 times differentiable with m-3 times differentiable inverse, with $D_*$ distance at most $\tau$ from the identity, and that pulls back $j^\#$ to a complex structure given by an element in $V_1$ with distance at most $\tau$ from j.*

*Proof of Lemma A.12*: Let exp: $T\Sigma \to \Sigma$ denote the metric's exponential map. The desired homeomorphism will have the form exp(v) where v is a suitably chosen section of $T\Sigma$. To elaborate, a map from $\Sigma$ to itself of this sort will pull back $j^\#$ to give an element from $V_1$ if the (1, 0) part of v obeys an equation with the schematic form

$$\bar{\partial} v_{1,0} + (1 - \Pi_V)\mathfrak{a} = 0$$

(A.22)

where the notation is as follows: First, $\Pi_V$ is the $L^2$ orthogonal projection to V. Second, $\mathfrak{a}$ is a section of $T^{0(1,0)}\Sigma \otimes T^{0(0,1)}\Sigma$ that is determined apriori by v, j and $j^\#$; it obeys

- $|\mathfrak{a}| \leq (1 - c_0^{-1})|\partial v| + c_0 \sup_\Sigma |j - j^\#|$
- $|\mathfrak{a}[v] - \mathfrak{a}[v']| \leq (1 - c_0^{-1})|\partial(v - v')| + c_0|v - v'||\nabla v|$.

(A.23)

Granted (A.23), a straightforward construction using the contraction mapping theorem in conjunction what is said in Theorems 3.5.2 and 5.4.2 in [Mo], and in Chapter 6 in [Mo], finds a Holder continuous section, v, of $T\Sigma$ such that exp(v) gives the desired homeomorphism.

Given $\tau > 0$ this lemma can be applied with $j^\# = j^\diamond$ when $\Delta$ from Part 9 is less than $c_0^{-1}\kappa_\tau$. Let $\psi^\tau$ denote the resulting homeomorphism and let $j^\tau$ denote the resulting almost complex structure. Then $(\varphi^\diamond \circ \psi^\tau, J', j^\tau) \in \mathcal{P}^m$ is such that $\varphi^\diamond \circ \psi^\tau$ maps $\Sigma$ to $C'$ and has $D_*$ distance at most $c_0(\tau + \rho^2)$ from $\varphi$. Meanwhile, and $j^\tau - j$ has norm less than $\tau$.

This last conclusion implies what is asserted in Proposition A.6.

### c) The manifold structure on $\mathcal{M}_{e,k}$

This part of the appendix focuses on the manifold structure for the space $\mathcal{M}_{e,k}$. The discussion contains proofs for Proposition 3.2 and 3.3.

To start, fix $J \in \mathcal{J}^m$ and introduce $\mathcal{M}_{e,k}|_J$ as defined in the previous part of this appendix. The first task is to define the coordinate charts in $\mathcal{M}_{e,k}|_J$ that give the latter its manifold structure when $J \in \mathcal{J}$ is from a suitable residual set. As explained below, this set can be taken to be the set of regular values for the map $\pi_\mathcal{J}$. In any event, fix a subvariety $C \in \mathcal{M}_{e,k}|_J$ and a point $(\varphi, J, j) \in \mathcal{P}^m$ such that $\varphi(\Sigma) = C$.

**Lemma A.13**: *Suppose that $k > 1$. Then there exists a neighborhood in $\pi_\mathcal{J}^{-1}(J)$ of $(\varphi, J, j)$ such that the map from the latter to $\mathcal{M}_{e,k}|_J$ defined by the rule $(\varphi', J, j') \to \varphi'(\Sigma)$ defines a 1-1 homeomorphism onto an open neighborhood of C. In the case when $k = 0$ or $k = 1$, this map is 1-1 up to the action of the group of $j$-holomorphic diffeomorphisms of $\Sigma$.*

By way of explanation with regards to the $k = 0$ and $k = 1$ cases, note that in the case $k = 0$, the vector space V is zero dimensional, and there is just a single complex structure to consider, $j$. Let $(\varphi_1, J, j)$ denote any given point in $\mathcal{P}^m$. Then composition of $\varphi_1$ with a $j$-holomorphic diffeomorphism of $S^2$ produces a new element, $(\varphi_2, J, j) \in \mathcal{P}^m$ such that $\varphi_1(\Sigma) = \varphi_2(\Sigma)$. To say something about the case $k = 1$, fix an element $(\varphi_1, J, j_1) \in \mathcal{P}^m$. There exists now a 2-dimensional torus of $j_1$-holomorphic diffeomorphisms of $\Sigma$; and any such diffeomorphism can be composed with $\varphi_1$ to produce a new element $(\varphi_2, J, j_1) \in \mathcal{P}^m$ with $\varphi_1(\Sigma) = \varphi_2(\Sigma)$. If $j_1$ is sufficiently close to $j$ in $V_1$, then the respective groups of $j_1$ and $j$-holomorphic diffeomorphisms are canonically isomorphic. Such an isomorphism is used to define the action of the space of $j$-holomorphic diffeomorphism on a neighborhood of $(\varphi, J, j)$ in $\pi_\mathcal{J}^{-1}(J)$.

*Proof of Lemma A.13*: The existence of a neighborhood of $(\varphi, J, j)$ that maps onto a neighborhood of C follows from Proposition A.6. Proposition A.6 asserts that the map in question is an open map. It remains only to establish that the neighborhood can be chosen so as to guarantee that the map is 1-1. To see why this is, suppose that $(\varphi_1, J, j_1)$ and $(\varphi_2, J, j_2)$ are both close to $(\varphi, J, j)$ and are such that $\varphi_1(\Sigma) = \varphi_2(\Sigma)$. In particular, suppose that $\varepsilon > 0$ has been specified and that $D_*(\varphi, \varphi_1) + \|j - j_1\|_\infty < \varepsilon$, and $(\varphi_2, j_2)$ obey this same inequality. This being the case, the map $\varphi_1$ can be written as $\exp_\varphi(\varphi_1')$ where $\exp: TX \to X$ is the metric's exponential map and $v_1$ is a section of $\varphi^*TX$ with norm $\|\varphi_1'\|_{1,2;v} \le c_0\varepsilon$. Likewise, the map $\varphi_2$ is defined by a corresponding section, $\varphi_2'$, of $\varphi^*TX$. Let $\mathfrak{D}_*$ denote the $(\varphi, J, j)$ version of the operator in (A.6). It follows from what is said surrounding (A.7) and (A.8) that the pair $(\varphi_2' - \varphi_1', j_2 - j_1)$ can be written as $\mathfrak{v} + \mathfrak{w}$ where $\mathfrak{v} = (\varphi', j') \in \ker(\mathfrak{D}_*)$ and $\mathfrak{w} = (\varphi'', j'')$ obeys

$$\|\varphi''\|_{1,2;\upsilon} + \|j''\|_\infty \leq c_0 \varepsilon \|\varphi'\|_{1,2;\upsilon}.$$
(A.24)

To continue, let $\eta$ denote the $(1, 0)$ part of $\varphi'$. To say that $\mathfrak{v} \in \ker(\mathfrak{D}_*)$ asserts neither more nor less than the fact that $(\eta, j')$ obeys (A.6). This understood, let $\Delta \subset \Sigma$ denote the set of critical points of $\varphi$. Write $\varphi^*TX$ on $C-\Delta$ as the direct sum $N \oplus N^\perp$ as in Part 2 of Section 2c and write $\mathfrak{D}$ in block diagonal form with respect to this splitting as done in (2.11). Write $\eta$ as $\eta = (\eta_0, \eta_1)$ with $\eta_0$ a section of $N$ and with $\eta_1$ a section of $N^\perp$. Given that $\varphi_1(\Sigma) = \varphi_2(\Sigma)$, the following is a consequence of (A.24): Fix $\delta > 0$ and let $D(\delta) \subset C$ denote the union of the disks of radius $\delta$ with centers at the points in $\Delta$. With $\delta$ fixed, then the norm of $\eta_0$ on $C-D(\delta)$ must be everywhere much less than that of $\eta_1$ if $\varepsilon$ is small. In particular, given $\delta > 0$ and $\tau > 0$, there exists $\varepsilon_{\delta,\tau} > 0$ such that if $\varepsilon < \varepsilon_{\delta,\tau}$, then it must be the case that $|\eta_0| \leq \tau|\eta_1|$ at all points in $C-D(\delta)$.

To see what to make of this, note first that (2.11) and (A.6) imply that

- $\|j'\|_\infty \leq c_\delta \tau \|\varphi'\|_{1,2;\upsilon}$.
- $|\bar\partial \eta_1| \leq c_\delta \tau |\eta_1|$ *at all points in* $C-D(\delta)$.

(A.25)

Given that the operator $\mathfrak{D}$ in (2.5) is elliptic, and given (A.6), the first bullet of (A.25) prevents $\varphi'$ from having the preponderance of its $L^2$ norm concentrated in $D(\delta)$. Given the second bullet in (A.25), this is possible when $\delta < c_0^{-1}$ and $\tau < c_0^{-1}$ only if $\eta_1$ differs by less than $c_0^{-1} \|\eta_1\|_{1,2;\upsilon}$ from its $L^2$ orthogonal projection to the kernel of $\bar\partial$. This means that $\eta_1 = 0$ when the genus, k, of $\Sigma$ is greater than 1; for $\ker(\bar\partial) = 0$ in this case. Granted this and the remarks in the preceding paragraph, this implies $(\varphi', j') = 0$ when $\varepsilon < c_0^{-1}$. This is to say that $(\varphi_1, j_1) = (\varphi_2, j_2)$.

Consider next the case when $k = 0$. In this case, $\Sigma = S^2$ and there is just a single complex structure up to diffeomorphisms. This the case, $V = \{0\}$ and $j_1 = j_2 = j$. Composing $\varphi_2$ with a j-holomorphic diffeomorphism results in a J-holomorphic map with the same image as the original. This group of holomorphic diffeomorphisms is generated by the elements in the kernel of $\bar\partial$. This understood, composition of $\varphi_2$ with a suitably small diffeomorphism produces a version of $\eta_1$ that is $L^2$ orthogonal to the kernel of $\bar\partial$ and thus is equal to zero. This version of $\varphi_2$ is therefore equal to $\varphi_1$.

Consider next the case of $k = 1$. In this case, the map $\varphi_2$ can be composed with a suitablly small normed, $j_2$-holomorphic diffeomorphism of $\Sigma$ so that the resulting version of $\eta_1$ is again $L^2$ orthogonal to the kernel of $\bar\partial$. The result of this composition makes $\varphi_2$ equal to $\varphi_1$ and so $j_2 = j_1$ as well.

This last lemma plays the central role in the ***Proof of Proposition 3.2***: What follows proves the proposition also for the cases when J $\in \mathcal{P}^m$ after modifying the assertion about the manifold structure as follows: The manifold structure is of differentiability class m-1.

Consider first the case when k > 1. Fix $(\varphi, J, j) \in \mathcal{P}^m$ so that $\varphi(\Sigma) = C$. Then $\Sigma$ with the complex structure j is the model curve. Lemma A.13 finds a neighborhood of $(\varphi, J, j) \in \pi_\mathcal{J}^{-1}(J)$ on which the map that sends an element $(\varphi', J, j') \in \pi_\mathcal{J}^{-1}$ to $\varphi'(\Sigma)$ defines a homeomorphism to a neighborhood of C in $\mathcal{M}_{e,k}|_J$. Use $\Phi_J$ to denote the latter map.

Meanwhile, it follows from what is said in Part 3 of Section a) of this appendix that there exists such a neighborhood that is homeomorphic to the zero locus of an m-1 times differentiable (or infinitely differentiable if $J \in \mathcal{P}$) map from a ball about the origin in the kernel of $\mathfrak{D}_C$ to the cokernel of $\mathfrak{D}_C$. Moreover, the latter map sends 0 to 0, and the homeomorphism in question sends 0 to $(\varphi, J, j)$. The map from the ball in kernel($\mathfrak{D}_C$) is the restriction of $f$ in (A.7) to elements of the form $(\mathfrak{v}, 0)$, and the relevant homeomorphism is the restriction to such elements of the map $\Psi$ in (A.8). This restriction is denoted in what follows by $\Psi_J$.

If $\mathfrak{D}_C$ has trivial cokernel, the domain of the homeomorphism $\Phi_J \circ \Psi_J$ is a ball about the origin in a vector space of dimension $\iota_e = e \cdot e - c \cdot e$. As a consequence, this homeomorphism defines a coordinate chart for a neighborhood of C. To see that the coordinate transition function between two such charts is suitably differentiable, remark first that the inverse function theorem gives $\pi_\mathcal{J}^{-1}(J)$ the structure of an m-1 times differentiable manifold (or infinitely differentiable if $J \in \mathcal{P}$) at the points $(\varphi, J, j)$ where the cokernel of $\mathfrak{D}_{(\cdot)}$ is trivial. The manifold coordinate charts are given by the homeomorphism $\Psi_J$. As a consequence, the coordinate transition functions for two such overlapping charts in $\pi_\mathcal{J}^{-1}(J)$ are m-1 times differentiable or infinitely differentiable as the case may be. This implies directly that the corresponding coordinate transition functions between the intersecting coordinate charts in $\mathcal{M}_{e,k}|_J$ are also m-1 times differentiable or infinitely differentiable.

The arguments in the cases when k = 0 or k = 1 are essentially identical to those above with the proviso that one must take into account that the holomorphic vector fields on $\Sigma$ define elements in the kernel of the operator $\mathfrak{D}_C$. This is done by viewing $\Phi_J \circ \Psi_J$ as a map from a ball about the origin in kernel($\mathfrak{D}_C$)/kernel($\bar{\partial}$) to a neighborhood of C in $\mathcal{M}_{e,k}|_J$. This change requires mostly cosmetic modifications to the arguments given above and so the details are omitted.

***Proof of Proposition 3.3***: This proposition is also proved for the case when $\mathcal{J}$ is replaced by $\mathcal{J}^m$ but for the modification to read that $\mathcal{M}_{e,k}|_J$ is a manifold of class m-1. The proposition follows from Proposition 3.2 if all $C \in \mathcal{M}_{e,k}|_J$ versions of $\mathfrak{D}_C$ have trivial cokernel. This is guaranteed if J is a regular value of the map $\pi_\mathcal{J}$. The Smale-Sard theorem [S] guarantees that this is so for almost complex structures from a residual set in $\mathcal{J}^m$. Proposition A.5 guarantees that this is so for almost complex structures from a residual set in $\mathcal{J}$.

### d) Proof of Proposition 3.5

The proof that follow for Proposition 3.5 also prove the assertion for the $J \in \mathcal{J}^m$ versions of $\mathcal{M}_{e,k}|_J$ with it understood that the adjective 'smooth' should be replaced by 'm-2 times differentiable'.

To start, fix $C \in \mathcal{M}_{e,k}|_J$ and an element $(\varphi, J, j) \in \mathcal{P}^m$ such that $\varphi(\Sigma) = C$. Consider first the case when $k > 1$. The arguments in the previous part of the appendix establish the following: There is a neighborhood $\mathcal{N} \subset \mathcal{P}^m$ with an embedding, $\Psi_{J,C}$, onto a neighborhood of C that maps the element $(\varphi, J, j)$ to C. The embedding sends any given $(\varphi', J, j)$ to $\varphi'(\Sigma)$. Use $\Psi_{J,C}$ to identify $\mathcal{N}$ with $\Psi_{J,C}(\mathcal{N}) \subset \mathcal{M}_{e,k}|_J$. Define next a map, $\vartheta: \mathcal{N} \times (\times_d \Sigma) \to \times_d X$ by the sending a given element $((\varphi', J, j); (p_1, \ldots, p_d))$ to $(\varphi'(p_1), \ldots, \varphi'(p_d))$. Introduce $\mathcal{G}_\mathcal{N} \subset (\mathcal{N} \times (\times_d \Sigma)) \times_d X$ to denote the graph of this map; it is a codimension 4d submanifold, thus a manifold of dimension $\iota_e + 2d$. The map from $\mathcal{G}_\mathcal{N}$ to $\mathcal{N} \times (\times_d X)$ is a suitably differentiable (m-2 times or infinitely as the case may be) map onto a neighborhood in $\mathcal{M}_{e,k,d}|_J$ of its intersection with the submanifold $\{C\} \times (\times_d X)$ in $\mathcal{M}_{e,k} \times (\times_d X)$. Note in this regard that the map from $\mathcal{G}_\mathcal{N}$ restricts as an embedding onto its image from a suitable neighborhood of any point $(C, (p_1, \ldots, p_d)), (\varphi(p_1), \ldots, \varphi(p_d))$ if no entry of $(p_1, \ldots, p_d)$ maps to a singular point of C.

The collection of the maps just defined give $\mathcal{M}_{e,k,d}|_J$ the structure of codimension 2d image variety. The fact that this image variety is a submanifold when $k = g$ follows from the final remark in the preceding paragraph because all subvarieties in $\mathcal{M}_{e,g}$ are non-singular.

The story in the case $k = 0, 1$ is analogous but for the necessity of replacing the neighborhood $\mathcal{N}$ with the a ball in the quotient of the kernel of $\mathfrak{D}_C$ by the subvector space consisting of the kernel of $\bar{\partial}$.

### e) Proof of Proposition 3.6

The $\pi_{d-2}$-inverse image in $\mathcal{M}_{e,k,d}$ of a point $\mathfrak{w} \in \times_{d-2} X$ is as asserted if $\mathfrak{w}$ is a regular value the map $\vartheta: \mathcal{P} \times (\times_{d-2} \Sigma) \to \times_{d-2} X$. This understood, then the first claim by the proposition follows because the set of regular values of a smooth map between finite dimensional manifolds is a residual set. It proves convenient to prove the remaining assertions, those of the four bullets, in reverse order.

To prove the fourth bullet, use Proposition 3.5 to conclude that $\mathcal{M}_{e,g,d-2}$ is a smooth manifold. The bullet now follows because the set of regular values of a smooth map between smooth, finite dimensional manifolds is a residual set.

To start the proof of the third bullet, reintroduce the set $\mathcal{X}$ from the proof of Proposition 1.1 in Section 4a. Given $\Xi \subset \mathcal{X}$, reintroduce also $\mathcal{M}_\Xi = \times_{(e',k') \in \Xi} \mathcal{M}_{e',k'}$ and $\mathcal{M}_{\Xi,d} \subset \mathcal{M}_\Xi \times (\times_d X)$. Recall that the latter is the subspace whose elements are of the form $((C_1, \ldots, C_n), (x_1, \ldots, x_d))$ with $\{x_m\}_{1 \leq m \leq d} \in \cup_{1 \leq i \leq n} C_i$. Here, n denotes the number of elements in $\Xi$. As noted in the proof of Proposition 1.1, the set $\mathcal{M}_{\Xi,d} \subset \mathcal{M}_\Xi \times (\times_d X)$ is a codimension 2d image variety that is given as the image of a map from a manifold whose dimension is 2d more than that of $\mathcal{M}_\Xi$. Use $\mathcal{P}_{\Xi,d}$ to denote this manifold and use $\mathfrak{p}_{\Xi,d}$ to denote the corresponding map. The set of points in $\times_{d-2} X$ that are regular values for every $\Xi \in \mathcal{X}$ versions of $\mathfrak{p}_{\Xi,d} \circ \pi_{d-2}$ is a residual set. Take $\mathfrak{w}$ is from this set. The various dimensions of the spaces $\{\mathcal{P}_{\Xi,d}\}_{\Xi \in \mathcal{X}}$, with (3.7) imply that $\pi_{d-2}^{-1}(\mathfrak{w})$ is empty unless $\Xi$ has but one entry, and this is either (e, g), (e, g-1) or (e, g-2). Given this last fact, the assertion of the third bullet of Proposition 3.6 follows directly from Proposition 3.1.

Proposition 3.1 also implies the first part of the assertion in the second bullet of Proposition 3.6. To prove the second part of the second bullet, introduce first $\underline{\mathcal{M}}_{e,g-2}$ to denote the closure in $\mathcal{M}_e$ of $\mathcal{M}_{e,g-2}$. Given an integer $N \geq 1$, use $\mathcal{O}^N \subset \underline{\mathcal{M}}_{e,g-2}$ to denote the subset of elements with $\eth$-distance greater than $N^{-1}$ from $\underline{\mathcal{M}}_{e,g-2} - \mathcal{M}_{e,g-2}$. Now suppose that $\mathfrak{w} \in \times_{d-2} X$ is such that $\pi_{d-2}^{-1}(\mathfrak{w}) \subset \mathcal{M}_{e,g-2,d}$ is finite. As noted, this is a residual set. A given point $\mathfrak{w}$ in this set has an open neighborhood, $O^N_\mathfrak{w} \subset \times_{d-2} X$, that is characterized as follows: A point $\mathfrak{w}'$ is in $O^N_\mathfrak{w}$ if $\pi_{d-2}^{-1}(\mathfrak{w}')$ has finite intersection with $\pi_\mathcal{M}^{-1}(\mathcal{O}^N)$. Denote by $O^N$ the union of all of the various versions of $O^N_\mathfrak{w}$. This is an open and dense set.

Now suppose that $\mathfrak{w} \in O^N$ is a point with the following property: If $C \in \pi_{d-2}^{-1}(\mathfrak{w})$ and $\pi_\mathcal{M}(C) \in \mathcal{O}^N$, then no singular point of C is an entry of $\mathfrak{w}$. This set is open, and it is also dense in $O^N$. Here is why it is dense: Suppose a singular point of C coincides with an entry of $\mathfrak{w}$. Moving this entry to any nearby point in C will not change C, and so not move C's singular points. Let $O^{N'} \subset O^N$ denote this open, dense set. This set is open and dense in $\times_{d-2} X$. Therefore, $O' = \cap_{N=1,2,\ldots} O^{N'} \subset \times_{d-2} X$ is a residual set.

The proof of the assertion made by the first bullet of Proposition 3.6 requires the following additional observation: There is an open and dense subset in $\times_{d-2} X$ that obeys the conclusions of the second bullet and whose $\pi_{d-2}$-inverse image in $\mathcal{M}_{e,k,d}$ for all k < g-2 is empty. To see why, let O denote this set. Meanwhile, let $\hat{O} \subset \times_{d-2} X$ denote the subset whose $\pi_{d-2}$-inverse image is respectively finite in $\mathcal{M}_{e,g-2,d}$ and empty in all k < g-2 versions of $\mathcal{M}_{e,k,d}$. As argued above, $\hat{O}$ is dense in $\times_{d-2} X$, and it is a consequence of Proposition 3.1 that this set $\hat{O}$ also open. Given that $O \subset \hat{O}$ is dense, and it is enough to prove that O is also an open subset of $\hat{O}$. To this end, remark that each $N \geq 1$ version of $O^{N'}$ as defined in the preceding paragraph is open, and so $O^{N'} \cap \hat{O}$ is open. Meanwhile, each $\mathfrak{w}$ in O has a neighborhood in $\hat{O}$ that lies in some $N \geq 1$ version of $O^{N'}$.

Consider now the assertion made by the first bullet of Proposition 3.6. To start, introduce the version of the space $\mathcal{P}|_J$ that is defined in Part a) in this appendix using the given almost complex structure J for the class e and genus g-1. Let $\Sigma$ denote the genus g-1 surface that is used to define $\mathcal{P}|_J$. If $(\varphi, j, J) \in \mathcal{P}|_J$, then there are precisely 2 points in $\Sigma$ that have the same $\varphi$-image. The map $\varphi$ embeds the complement of these two points, but embeds only some neighborhood of each. Let $\mathcal{P}_\diamond \subset \mathcal{P}|_J \times \Sigma$ denote the subset of the form $((\varphi, j, J), z)$ where z is one of these two special points. This $\mathcal{P}_\diamond$ is a smooth, 2d - 2 dimensional manifold, and the restriction to $\mathcal{P}_\diamond$ of the projection from $\mathcal{P}|_J \times \Sigma$ to its first factor is a 2-1 covering map.

For each $i \in \{1, \ldots, d-2\}$, use $\mathcal{X}_i \subset \mathcal{P}_\diamond \times (\times_{d-2}\Sigma)$ to denote the codimension 2 submanifold of points $(((\varphi, j, J), p), (p_1, \ldots, p_d))$ where $p = p_i$. Since $\mathcal{P}_\diamond \times (\times_{d-2}\Sigma)$ has dimension $4d - 6$, so each $i \in \{1, \ldots, d-2\}$ version of $\mathcal{X}_i$ has dimension $4d - 8$. Meanwhile, define $\vartheta \colon \mathcal{P}_\diamond \times (\times_{d-2}\Sigma) \to \times_{d-2} X$ by the rule

$$(((\varphi, j, J), p), p_1, \ldots, p_{d-2}) \to (\varphi(p_1), \ldots, \varphi(p_{d-2})) \,.$$
(A.26)

This is a smooth map. For each $i \in \{1, \ldots, d-2\}$, let $\vartheta_i$ denote the restriction of $\vartheta$ to $\mathcal{X}_i$. Let $O \subset \times_{d-2} X$ denote the open, dense set that is defined two paragraphs back, and let $O^\# \subset O$ denote the residual set of points that are simultaneous regular values for $\vartheta$ and for each $\vartheta_i$. In particular, if $\mathfrak{w} \in O^\#$, then there are at most a countable set of points $\vartheta^{-1}(\mathfrak{w})$ that lie in $\mathcal{X}_i$, and these points have no accumulations.

The last remark has the following consequence: Let $\mathfrak{w} \in O^\#$. Then there is at most a countable subset of subvarieties from $\mathcal{M}_{e,g-1,d-2}$ in $\pi_{d-2}^{-1}(\mathfrak{w})$ that have their singular point at an entry of $\mathfrak{w}$. Moreover, this subset must be finite. Indeed, were this set infinite, then any limit point must converge in $\mathcal{M}_e \times (\times_{d-2} X)$ to a point in $\mathcal{M}_{e,g-2,d-2}$ that lies in the $\pi_{d-2}$-inverse image of $\mathfrak{w}$ and has a singular point at an entry of $\mathfrak{w}$. But the condition that $\mathfrak{w} \in O$ forbids this event.

**f) Proof of Proposition 3.7**

The proof has thirteen parts. Part 1 proves the assertion of the first bullet. Part 2 proves Item a) of the second bullet, and Parts 2-7 prove Item b) of the second bullet. With regards to the other items in the second bullet, Part 8 proves Item c), Part 9 proves Item d), Part 10 proves Item e) and Parts 11-13 prove Item f).

*Part 1*: This part of the proof considers the first bullet in the proposition. To this end, fix $n \in \{0, 1, 2\}$ and let $\mathcal{M}_{e,g-n,2d-4} \subset \mathcal{M}_{e,g-n} \times (\times_{d-2} X) \times (\times_{d-2} X)$. In the case $n = 0$, this is a smooth manifold of dimension $6d - 8$; and for $n > 0$, it is an image variety of dimension $6d-8-2n$. Write the map $\pi_{2d-4}$ as $\pi_{d-2} \times \pi_{d-2}$. It maps a $6d-8-2n$ dimensional manifold to one of dimension $8d - 16$. Thus, if $d > 4$, then there is an open and dense subset of pairs $(\mathfrak{w}, \mathfrak{w}') \in (\times_{d-2} X) \times (\times_{d-2} X)$ such that both are regular points, they have distinct entries, and such that there is no curve in $\{\mathcal{M}_{e,g-n}\}_{n=0,1,2}$ that contains all entries of both $\mathfrak{w}$ and $\mathfrak{w}'$. This the case, it follows from Proposition 3.1 that there exists $\kappa > 1$ such that the $\partial$-distance between any curve in $\mathcal{M}^{\mathfrak{w}}$ and any curve in $\mathcal{M}^{\mathfrak{w}'}$ is at least $\kappa^{-1}$.

*Part 2*: This part begins the arguments for the assertions that are made by the second bullet. To this end, remark that Item a) is straightforward, so the first concern is that of Item b).

Introduce the version of the space $\mathcal{P}$ that is defined in Part a) in this appendix using the cohomology class e and genus g - 2. Define $\vartheta_1 \colon \mathcal{P} \times \Sigma \to X$ so as to send a given point $((\varphi, j, J), p) \to \varphi(p)$. Introduce next the map $\vartheta_{d-2} \colon (\mathcal{P} \times \Sigma) \times (\times_{d-2} \Sigma) \to \times_{d-2} X$ that sends any given $((\varphi, j, J), p), (p_1, \ldots, p_{d-2})$ to $(\varphi(p_1), \ldots, \varphi(p_{d-2}))$.

Fix $J \in \mathcal{J}$, a regular value for the map from $\mathcal{P}$ to $\mathcal{J}$. If x is a regular value of $\vartheta_1$ on $\mathcal{P}|_J \times \Sigma$, then $\vartheta_1^{-1}(x)$ is a submanifold of dimension (2d - 6) in $\mathcal{P}|_J \times \Sigma$ and so the image via $\vartheta_{d-2}$ of $\vartheta_1^{-1}(x) \times (\times_{d-2} \Sigma)$ is the image of a smooth map from a manifold of dimension 4d - 10 into one of dimension 4d - 8. As a consequence, there is a residual set in $\times_{d-2} X$ with the property that there is no subvariety in J's version of $\mathcal{M}_{e,g-2}$ that contains x and all entries of $\mathfrak{w}$. Moreover, if $\mathfrak{w}$ is a regular, value, there are but a finite number of curves in $\mathcal{M}_{e,g-2}$ that contain all entries of $\mathfrak{w}$ and so there is a positive lower bound to the distance between x and the union of all such curves. These observations establish Item b) of Proposition 3.7 when x is not a critical value of $\vartheta_1$ on $\mathcal{P}|_J \times \Sigma$.

What follows explains how to deal with critical values of $\vartheta_1$. To this end, suppose that the set of critical points of $\vartheta_1$ on $\mathcal{P}|_J \times \Sigma$ is an image variety with at least codimension 3. This means that any given critical point has a neighborhood in this critical locus that is the image of a manifold with components of dimension 2d - 5 or less. Thus, there is a residual set $\mathfrak{C} \subset \times_{d-2} X$ with the following property: Let $\mathfrak{w}$ denote a point from $\mathfrak{C}$. Then $\vartheta_{d-2}^{-1}(\mathfrak{w})$ has no points of the form $((\varphi, j, J), p, (p_1, \ldots, p_{d-2}))$ where the point $((\varphi, j, J), p)$ is a critical point of $\vartheta_1$.

With this last point in mind, suppose that a given point $((\varphi, j, J), z) \in \mathcal{P}|_J \times \Sigma$ is not in the critical point locus of $\vartheta_1$ on $\mathcal{P}|_J \times \Sigma$, but even so maps via $\vartheta_1$ to x. Such a point has a neighborhood in $\mathcal{P}|_J \times \Sigma$ that intersects $\vartheta_1^{-1}(x)$ as a codimension 2d - 6 submanifold. Granted the preceding, it follows that there exists a residual set in $\times_{d-2} X$ with the property that if $\mathfrak{w}$ is in this set, then there is point $((\varphi', j', J), z')$ from this neighborhood and no point $(p_1, \ldots, p_{d-2}) \in \times_{d-2} \Sigma$ which maps via $\vartheta_{d-2}$ to $\mathfrak{w}$. It then follows from what was said in the preceding paragraph that if $\mathfrak{w}$ is also from the residual set $\mathfrak{C}$ and is a regular value, then there is a positive lower bound to the distance from x to any curve in $\mathcal{M}_{e,g-2}$ that goes through all entries of $\mathfrak{w}$.

Given the conclusions from these last three paragraphs, Item b of Proposition 3.7 follows with a proof of the following claim:

*There exists a residual set in $\mathcal{J}$ such that if J comes from this set, then the critical point locus of the corresponding version of $\vartheta_1$ is a dimension 3 image variety.*
(A.27)

By the way, a 3-dimension critical locus is expected given the following observation: If d is large, then the subset of non-surjective linear maps from $\mathbb{R}^{2d-6}$ to $\mathbb{R}^4$ is a dimension 3 image variety in the Euclidean space of (2d-6) by 4 matrices.

*Part 3*: This part of the proof, and Parts 4-7 set up the machinery to analyze the critical point structure of the map $\vartheta_1$ for generic J. The constructions made here and in Part 4-7 to analyze the set of critical points of $\vartheta_1$ are also used subsequently with minor modifications to prove Items c)-f) of Proposition 3.7 and to prove Proposition 3.8.

The analysis that follows of $\vartheta_1$'s critical locus uses the notation from Part a) of this appendix, and in particular the spaces $\mathcal{E}_*$, $\mathcal{J}^m$, $\mathcal{P}^m$ and the vector bundle $\mathcal{V}$. Fix $n \in \{0, 1, 2, 3\}$ and introduce $\mathcal{F}_n \to \mathcal{P}^m$ to denote the restriction from $(\mathcal{E}_* \times \mathcal{J}^m \times V_1) \times \Sigma$ of the space of n dimensional orthonormal frames in the bundle $(T\mathcal{E}_* \oplus V) \oplus T\Sigma$. Note that this fiber bundle has a free action of the orthogonal group SO(2d-4-n). View $\mathcal{F}_n$ and the bundle $\mathcal{V}$ as sitting over $\mathcal{P}^m \times \Sigma$.

Define a fiber preserving map over $\mathcal{P}^m \times \Sigma$ from the fiber bundle $\mathcal{F}_n$ to the vector bundle $\oplus_{2d-4-n}(\mathcal{V} \oplus \vartheta_1^*TX)$ by the following rule: At a given $((\varphi, j, J), p)$, this map sends an n tuple $(((\varphi_1, j_1), v_1) \ldots, ((\varphi_{2d-4-n}, j_{2d-4-n}), v_{2d-4-n}))$ to the element whose component in the $\vartheta_1^*TX$ part of the i'th summand is $\varphi_i|_p + \varphi_* v_i$ and whose component in the $\mathcal{V}$ part of the i'th summand is the $(\varphi_i, j_i)$ version of the section of $\varphi^*T_{J(1,0)}X \otimes T_j^{0,1}\Sigma$ that is defined by the left hand side of (A.6). Use $\mathfrak{X}$ to denote this map.

By way of motivation for introducing $\mathcal{F}_n$ and the map $\mathfrak{X}$, note that a given element $((\varphi, J, j), p) \in \mathcal{P}|_J \times \Sigma$ is a critical point of the map $\vartheta_1$ if and only if there is a dimension 2d-7 linear subspace in $T(\mathcal{P}|_J \times \Sigma)|_{((\varphi,j,J),p)}$ with the following property: If $((\varphi', j'), v)$ is in this space, then the vector $\varphi'|_p + \varphi_* v = 0$ in $\varphi^*TX|_p$. Here, $\varphi_*$ denotes the differential of $\varphi$, viewed as a homomorphism from $T\Sigma$ to $\varphi^*TX$. Keep in mind for what follows that a given $(\varphi', j') \in C^\infty(\Sigma, \varphi^*TX) \oplus V$ lies in the tangent space at $(\varphi, j, J)$ to $\mathcal{P}|_J$ if and only if (A.6) holds.

As is explained momentarily, the locus where $\mathfrak{X}$ hits the zero section is an m-3 times differentiable manifold with a free action of SO(2d-4-n). Given that this is true, then the projection induced map from $\mathfrak{X}^{-1}(0)$ to $\mathcal{J}^m$ is a Fredholm map. Furthermore, if m is large, the Sard-Smale theorem can be invoked to conclude that there is a residual set of regular values in $\mathcal{J}^m$. An argument like that given to prove Proposition A.5 proves that there is a residual set of regular values in $\mathcal{J}$ for the map from $\mathfrak{X}^{-1}(0)|_\mathcal{J}$ to $\mathcal{J}$. Let $J_*$ denote such an almost complex structure. The inverse image over $J_*$ in $\mathcal{F}_n$ is empty when $n < 3$, and it is a manifold with a free action of SO(3) whose quotient has dimension 3 if $n = 3$. By construction, the image in $\mathcal{P}|_{J_*} \times \Sigma$ of this manifold via the projection map is the critical locus of the restriction to $\mathcal{P}|_{J_*} \times \Sigma$ of $\vartheta_1$.

*Part 4*: The locus $\mathfrak{X}^{-1}(0)$ is a submanifold near a zero of $\mathfrak{X}$ if the differential of $\mathfrak{X}$ is surjective at the point. This understood, what follows analyzes this differential. To start, let $((\varphi, j, J), p), \mathfrak{p} = (((\varphi_1, j_1), v_1) \ldots, ((\varphi_{2d-4-n}, j_{2d-4-n}), v_{2d-4-n}))$ denote a given zero of $\mathfrak{X}$. A tangent vector to the SO(2d-4-n) quotient at the orbit of this zero lifts to give the following three part data set: The first part of the data is an element $\mathfrak{w} = (\varphi', j', J') \in$ $T\mathcal{P}^m$, this a solution to (A.5). The second part is a tangent vector v to $\Sigma$ at p; and the third part is a 2d-4-n tuple $q = (((\lambda_1, \mathfrak{k}_1), x_1), \ldots, ((\lambda_{2d-4-n}, \mathfrak{k}_{2d-4-n}), x_{2d-4-n}))$ in $(T\mathcal{E}_*|_\varphi \times V) \times T\Sigma|_p$ that is orthogonal to the span of $p$.

What follows says what it means for the differential of $\mathfrak{X}$ at this given zero in $\mathcal{F}_n$ to be surjective: Fix data $\{(g_i, w_i)\}_{1 \le i \le 2d-4-n}$, with $g_i$ a section of $\varphi^* T_{J(1,0)} X \otimes T_j^{0,1}\Sigma$ and $w_i \in (\varphi^* TX)|_p$. First,

$$\lambda_i|_p + \varphi_* x_i + \nabla_v \varphi_i|_p + (\nabla_v \varphi_*) v_i = w_i \quad \text{for each } i \in \{1, \ldots, 2d-4-n\}. \tag{A.28}$$

Second, the data $\{g_i\}_{1 \le i \le 2d-4-n}$ gives the inhomogeneous term in a system of differential equations that are described next. To this end, fix a point $p' \in \Sigma$, a holomorphic coordinate u for a disk centered at p and coordinates for a ball in $\mathbb{R}^4$ centered on $\varphi(p')$. This done, then the image of $\varphi$ appears as a map $u \to \varphi(u)$ from a disk centered at the origin in $\mathbb{C}$ to a ball centered at the origin in $\mathbb{R}^4$. Use $(\varphi_i, j_i)$ and $(\lambda_i, \mathfrak{k}_i)$ to denote the $C^\infty(\Sigma; \varphi^* TX) \oplus V$ components from the respective i'th entries of $p$ and $q$. Write $\varphi_i$ as a map $u \to \varphi_i(u) \in \mathbb{R}^4$ and write $\lambda_i$ as a map $u \to \lambda_i(u) \in \mathbb{R}^4$. Meanwhile write $\varphi'$ as a map $u \to \varphi'(u)$. Then equation (A.6) for the pair $(\varphi_i, j_i)$ and for (A.5) for $(\varphi', j', J')$ and the constraint on $(\lambda_i, \mathfrak{k}_i)$ can be written on these coordinates as

- $(\bar{\partial}\varphi_i)_{J(1,0)} + (\nabla_{\varphi_i} J)\cdot(\bar{\partial}\varphi)_{J(0,1)} + j_i(\partial\varphi)_{J(1,0)} = 0.$
- $(\bar{\partial}\varphi')_{J(1,0)} + (\nabla_{\varphi'} J)\cdot(\bar{\partial}\varphi)_{J(0,1)} + j'(\partial\varphi)_{J(1,0)} + J'(\bar{\partial}\varphi)_{J(0,1)} = 0.$
- $(\bar{\partial}\lambda_i)_{J(1,0)} + (\nabla_{\lambda_i} J)\cdot(\bar{\partial}\varphi)_{J(0,1)} + \mathfrak{k}_i(\partial\varphi)_{J(1,0)} + j_i(\partial\varphi')_{J(1,0)} + j'(\partial\varphi_i)_{J(1,0)}$
  $\quad + (\nabla_{\varphi_i} J)\cdot(\bar{\partial}\varphi')_{J(0,1)} + (\nabla_{\varphi'} J)\cdot(\bar{\partial}\varphi_i)_{J(0,1)} + (\nabla_{\varphi_i} J)\cdot((\nabla_{\varphi'} J)\cdot(\bar{\partial}\varphi)_{J(1,0)})_{J(0,1)}$
  $\quad\quad + (\nabla_{\varphi'} J)\cdot((\nabla_{\varphi_i} J)\cdot(\bar{\partial}\varphi)_{J(1,0)})_{J(0,1)} + (\nabla^{\otimes 2}_{(\varphi_i,\varphi')} J)\cdot(\bar{\partial}\varphi)_{J(0,1)}$
  $\quad + J'\cdot(\bar{\partial}\varphi_i)_{J(0,1)} + (\nabla_{\varphi_i} J')\cdot(\bar{\partial}\varphi)_{J(0,1)} + (\nabla_{\varphi_i} J)\cdot((J'(\bar{\partial}\varphi')_{J(1,0)})_{J(0,1)} + j_i J'(\partial\varphi)_{J(0,1)} = g_i.$

(A.29)

Here, $(\cdot)_{J(1,0)}$ and $(\cdot)_{J(0,1)}$ denote the respective $(1, 0)$ and $(0, 1)$ parts of the indicated vector, this as defined by J.

*Part 5*: To see about solving (A.28) and (A.29) for a given set $\{(g_i, w_i)\}_{1 \le i \le 2d-4-n}$. Note first that the operator that sends $(\lambda_i, \mathfrak{k}_i)$ to the section of $\varphi^* T_{J(1,0)} X \otimes T_j^{0,1}\Sigma$ given near $p'$ as

$$(\bar{\partial}\lambda_i)_{J(1,0)} + (\nabla_{\lambda_i} J)\cdot(\bar{\partial}\varphi)_{J(0,1)} + \mathfrak{k}_i(\partial\varphi)_{J(1,0)} \tag{A.30}$$

has finite dimensional cokernel. Consider first the case where this cokernel is trivial. Granted that such is the case, then each $i \in \{1, \ldots, 2d-4-n\}$ version of the equation in the third bullet of (A.29) can be solved given any choice of tangent vector (φ´, J´, j´) to $\mathcal{P}$ at (φ, J, j). As a consequence, it is enough to consider the case of (A.28) and (A.29) when all i ∈ {1, …, 2d-4-n} versions of $g_i$ are zero. Moreover, the equation in the second bullet of (A.29) can be solved given any choice of J´. This understood, the question is whether a given value for the 2d-4-n tuple $(\lambda_1|_p, \ldots, \lambda_{2d-4-n}|_p)$ can be obtained from a solution of solution $\{g_i = 0\}_{1 \le i \le 2d-4-n}$ versions of the equation in the third bullet of (A.29) as defined using a suitable choice of J´ and then a suitable solution (φ´, J´, j´) to the second bullet in (A.29).

To see that the answer to this question is affirmative, take p´ = p so that the third bullet in (A.29) describes the constraint on the data set $(\lambda_i, \ell_i)$ near p. Fix a small disk D ⊂ ℂ in the domain of the coordinate u very near the origin on which the map φ restricts as an embedding. Take D to contain the origin in ℂ if φ's differential there is injective. Note in this regard that φ can have at most one critical point. In any event, choose J´ so that its pull-back as an endomorphism of φ*TX has support on D and such that it acts as zero on the tangent space to φ(D). This implies in particular that J´($\bar\partial$ φ) = 0 and J´(∂φ) = 0. This the case, then φ´ and j´ are taken to be zero also, and so the second bullet in (A.29) is satisfied. Meanwhile, the first and third read

- $(\bar\partial \varphi_i)_{J(1,0)} + (\nabla_{\varphi_i} J) \cdot (\bar\partial \varphi)_{J(0,1)} + j_i (\partial\varphi)_{J(1,0)} = 0$.
- $(\bar\partial \lambda_i)_{J(1,0)} + (\nabla_{\lambda_i} J) \cdot (\bar\partial \varphi)_{J(0,1)} + \ell_i (\partial\varphi)_{J(1,0)} + J´ \cdot (\bar\partial \varphi_i)_{J(0,1)} - J´ \cdot \nabla_{\varphi_i} (\bar\partial \varphi)_{J(0,1)} = 0$.

(A.31)

*Part 6*: To exploit this last equation, let o ⊂ Σ denote the critical point of φ if such a point exists; and let N ⊂ φ*(TX)|_{Σ-o} denote the local normal bundle, thus the orthogonal complement to the image of TΣ|_{Σ-o} via φ's differential. In the case o exists the subbundle N none-the-less extends over o as a subbundle in φ*TX. This extension is defined by using the n = 1 version of (2.2) near o to extend N over o as the span of $\frac{\partial}{\partial w}$. In any case, let Π: φ*(TX)|_D → N|_D denote the orthogonal projection and let $\eta_i$ denote Πφ_i. What follows explains why no constant c > 0 and no unit vector α = $(\alpha_i)_{1 \le i \le 2d-4-n}$ ∈ $\mathbb{R}^{2d-4-n}$ exist such that the following is true: Set α·η = $\sum_{1 \le i \le 2d-4-n} \alpha_i \eta_i$. Then |α·η| bounds c times the norm of Π·(($\bar\partial$ (α·η))_{J(0,1)}) on a neighborhood of the origin in ℂ if α·η|_0 = 0. Indeed, suppose that this assertion were false. Then it follows from the top bullet in (A.31) that |Π·d(α·η)| ≤ c|α·η| and so α·η vanishes on the whole of D. As a consequence, the projection of α·φ = $\sum_{1 \le i \le 2d-4-n} \alpha_i \varphi_i$ to N must vanish on the whole of Σ–o, and so α·φ defines a tangent vector on Σ–o to the image of φ. Given the top line of (A.31), this implies that α·φ defines a holomorphic vector field on Σ. It must therefore vanish identically.

What follows is a direct consequence of the preceding conclusion. There exists in any given open set in any given neighborhood of the origin a set of distinct points, none the origin, and with the following property: Let W denote this set. Let φ denote the map from $\oplus_{2d-4-n} \mathbb{R}$ to $\oplus_{u \in W} N|_u$ whose entry in the factor labeled by any given u ∈ W sends $(\alpha_1, \ldots, \alpha_{2d-4-n})$ to the point $(\bar{\partial}(\alpha \cdot \eta))|_u$. Then this map φ is an injection. Fix such a set W, but one with no element at the origin

Let (u, u´) → G(u, u´) denote the restriction to D × D of the Green's function for the operator on Σ that appears on D as in (A.29). Since (A.30) depicts a first order, elliptic operator with symbol that of $\bar{\partial}$, so G(·, ·) can be described as follows. Fix adapted coordinates (z, w) centered at φ(p). Then

$$G(u, u') = \frac{1}{2\pi} \frac{1}{u - u'} + \mathcal{O}(1)$$

(A.32)

To continue, suppose that J´ has support in a union of radius Δ << 1 disks that are centered on the points in the set W. The radius, d, should be much less than the distance between the points in W and much less than the distance to the origin. Require in addition that if u ∈ W, then J´ on the radius Δ disk centered at u is as follows: Use parallel transport from the center of the disk along the radial geodesics to identify the fibers of $\varphi^*\text{Hom}(T_{J(0,1)}X \otimes T_{J(1,0)}X)$ over this radius Δ disk with its fiber over u. This done, then

$$J' = \chi_\Delta J'|_u$$

(A.33)

where $\chi_\Delta$ is a positive function with value 1 at the center of the disk, with integral equal to 1. It follows from (A.32) that the solution to any given i ∈ {1, …, 2d-4-n} version of the equation in the second bullet of (A.31) appears in D as

$$\Pi(\lambda_i)_{J(1,0)}|_0 = \Delta^2 \left( \sum_{u \in W} \frac{1}{u} (J' \cdot (\bar{\partial}\eta_i)_{J(0,1)})|_u + r_i \right)$$

(A.34)

where

$$|r_i| \le c \Delta \sum_{u \in W} |J'|_u| \sup_{|(\cdot)-u|<\Delta} |\nabla(\bar{\partial}\eta_i)_{J(1,0)}|.$$

(A.35)

Here, c ≥ 1 depends on W, but not on Δ if the latter is small.

Given what was said about the map φ, it follows from (A.34) and (A.35) that J´ can be chosen when Δ is small so that the vector $(\lambda_1|_0, \ldots, \lambda_{4d-2-n}|_0) \in \oplus_{2d-4-n}(\varphi^*TX)|_p$ has any desired value.

*Part 7*: What follows explains how to modify the preceding argument for the case when (A.30) has a non-trivial cokernel. Again take J´ so that its pull-back as an endomorphism of φ*TX has support on D and such that it acts as zero on the tangent space to φ(D). This done, then the second equation in (A.29) is satisfied with φ´ = 0 and j´ = 0. The third equation in (A.29) has left side given by the left side of (A.31) but with the inhomogeneous term $g_i$ on the right hand side. The operator that is depicted in (A.30) is Fredholm, and so its cokernel is the kernel of its $L^2$ adjoint. Let $x$ denote an element in the kernel of this adjoint with $L^2$ norm equal to 1. Such an element constrains J´ in the sense that

$$\langle x, J´\cdot(\bar{\partial}\varphi_i)_{J(0,1)} \rangle - J´\cdot\nabla_{\varphi_i}(\bar{\partial}\varphi)_{J(0,1)}\rangle_2 = \langle x, g_i\rangle_2 .$$
(A.36)

Here, the notation has $\langle\,,\,\rangle_2$ denoting the $L^2$ inner product on $C^\infty(\Sigma;\varphi^*T_{J(1,0)}X \otimes T_j^{0,1}\Sigma)$.

To see how to procede, digress momentarily and let u´ ∈ D denote a given point. Define a linear functional on $C^\infty(\Sigma;\varphi^*(\text{Hom}(T_{J(0,1)};T_{J(1,0)})))$ by the rule

$$J´´ \to \langle x, J´´\cdot(\bar{\partial}\eta_i)_{J(0,1)}\rangle|_{u´}$$
(A.37)

Let q denote the dimension of the cokernel of the operator that is depicted in (A.30), and let $\{x_\alpha\}_{1\le\alpha\le q}$ denote an $L^2$ orthonormal basis for this cokernel. Given u´ ∈ D, denote the $x_\alpha$ and $\eta_i$ version of (A.37) by $L_{(\alpha,i);u´}$. Fix r > 0 but much smaller than the diameter of the disk D. Let $A_r \subset D$ denote the annulus with inner radius r and outer radius 2r. Fix a set W ⊂ $A_r$ as done previously. Denote this set by $W_r$ to distinguish it from a subsequent choice of another version of W. Fix again Δ = $\Delta_r$ > 0 but much smaller than the distance between the points in $W_r$. For each u ∈ $W_r$, let $D_r(u)$ denote the disk of radius $\Delta_r$ centered at u. There exists a set $W_r(u) \subset D_r(u)$ of less than $c_0 q$ points with the following property: The map from the span of $\{x_\alpha\}_{1\le\alpha\le q}$ to $\oplus_{u´\in W_r(u)}(\varphi^*T_{J(1,0)}X \otimes T_j^{0,1}\Sigma)|_{u´}$ that sends $x \to \oplus_{u´\in W_r(u)} x|_{u´}$ is injective. This understood, it follows that the set of linear functionals

$$\{L_{(\alpha,i);u´}: 1 \le \alpha \le q, 1 \le i \le 2d-4-n, \text{ and } u´ \in \cup_{u\in W_r} W_r(u)\}$$
(A.38)

are linearly independent if $\Delta_r$ is sufficiently small. Fix next a very small, but positive number D << $\Delta_r$ such that the disks of radius Δ´ centered on the points in each u ∈ $W_r$ version of $W_r(u)$ are disjoint.

With the digression now over, what follows directly explains why it is sufficient to consider the case of (A.29) where all $\{g_i\}_{1\le i\le 2d-4-n}$ are zero. The key observation in this regard is that the set in (A.38) are linearly independent. It follows as a consequence that (A.36) can be satisfied for each index i and each cokernel element $x$ with J´ chosen to have support in the various radius D disks centered at the points in $\cup_{u\in W_r} W_r(u)$. In particular, in the radius D disk centered on a point u´ from this set, J´ can be written as in (A.33) with u´ replacing u and D replacing Δ.  This is because such a choice for J´ allows the left hand side of (A.36) to be written as

$$D^2 \left( \sum_{u \in W_r} \sum_{u' \in W_r(u)} \langle x, J'\cdot(\overline{\partial}\eta_i)_{J(0,1)}\rangle|_{u'} + \mathfrak{e}_i \right)$$

(A.39)

where $\langle \, , \, \rangle$ denotes the inner product on $\varphi^*T_{J(1,0)}X \otimes T_j^{0,1}\Sigma$ and where $\mathfrak{e}_i$ obeys

$$|\mathfrak{e}_i| \leq c_0 \, D \sum_{u \in W_r} \sum_{u' \in W_r(u)} ||J'|_{u'}| \, |(\nabla\overline{\partial}\eta_i)_{J(0,1)}| \, .$$

(A.40)

Granted such a choice for J´, the obstruction to simultaneously satisfying the inhomogeneous counterpart to the various i ∈ {2d-4-n} versions of the lower equation in (A.31) vanishes.

Before leaving this topic, note that the choice for J´ as just described can be made so that

$$|J'| \leq c_0 \, D^{-2} \, \|g_i\|_1 \, c(r)$$

(A.41)

where $c(r) \geq 1$ depends on r.  Here, $\|\cdot\|_1$ denotes the $L^1$ norm on $C^\infty(\Sigma; \varphi^*T_{J(1,0)}X \otimes T_j^{0,1}\Sigma)$. It is a consequence of (A.34) that for very small D, such a choice leads to a collection $\{\lambda_i\}_{1 \leq i \leq 2d-4-n}$ that obeys

$$|\lambda_i|(0) \leq c_0 \, \|g_i\|_1 \, r^{-1} \, c(r).$$

(A.42)

Now consider solving the lower equation in (A.31) so as to obtain any given value of $\oplus_{1 \leq i \leq 2d-4-n} \varphi^*T_{J(1,0)}X|_0$.  To do this, fix s << r and a set W as before, but now chosen to lie in the annulus $A_s$ where $s \leq |u| \leq 2s$.  Fix Δ > 0 as before, but much less than s.  Write J´ as a sum of two contributions, $J' = J'_s + J'_r$ where $J'_s$ has support in the disks of radius Δ centered at the points in W, and where $J'_r$ has support in the disks of radius D centered on the points in $\cup_{u \in W_r} W_r(u)$.  This $J'_r$ term is used to cancel the cokernel obstructions given by (A.36) in the case where

$$g_i = -(J'_s\cdot(\overline{\partial}\varphi_i)_{J(0,1)} - J'_s\cdot\nabla_{\varphi_i}(\overline{\partial}\varphi)_{J(0,1)}) \, .$$

(A.43)

Meanwhile $J'_s$ is given by (A.33) as in the case of vanishing cokernel.  This done, it follows from (A.34) and (A.42) that

$$\lambda_i(0) = \Delta^2 \sum_{u \in W} \tfrac{1}{u} ((J'_s \cdot (\bar{\partial}\eta_i)_{J(0,1)})|_u + r_i) + e_i$$

(A.44)

where $r_i$ is given in (A.34) and (A.35); and where $e_i$ obeys

$$|e_i| \leq c_0\, r^{-1} c(r) \Delta^2 \sum_{u \in W} (|(J'_s \cdot (\bar{\partial}\eta_i)_{J(0,1)})|_u| + \Delta\, |J'|_u|\sup_{A_s} |\nabla \bar{\partial}\eta_i|)$$

(A.45)

What follows is a direct consequence of (A.35), (A.44) and (A.45). With r fixed, and s chosen so that $s \ll r^{-1} c(r)$, and with $\Delta$ very small, any given element in the vector space $\oplus_{1 \leq i \leq 2d-4-n} \varphi^* T_{J(1,0)} X|_0$ can be obtained be realized as $(\lambda_1|_0, \ldots, \lambda_{2d-4-n}|_0)$ using a suitable choice for the collection $\{J'_s|_u\}_{u \in W}$.

*Part 8*: This part argues for Item c) of the second bullet of Proposition 3.7. To this end, introduce the version of the space $\mathcal{P}$ that is defined in Part a) in this appendix using the class e and genus g - 1. Define $\vartheta_1: \mathcal{P} \times \Sigma \to X$ by the rule$((\varphi, j, J), p) \to \varphi(p)$. Meanwhile, use $\vartheta_{d-2}: (\mathcal{P} \times \Sigma) \times (\times_{d-2}\Sigma) \to \times_{d-2} X$ to denote the map that sends any given $((\varphi, j, J), p), (p_1, \ldots, p_{d-2})$ to $(\varphi(p_1), \ldots, \varphi(p_{d-2}))$.

Fix $J \in \mathcal{J}$, a regular value for the map from $\mathcal{P}$ to $\mathcal{J}$. If x is a regular value of $\vartheta_1$ on $\mathcal{P}|_J \times \Sigma$, then $\vartheta_1^{-1}(x)$ is a submanifold of dimension (2d - 4) in $\mathcal{P}|_J \times \Sigma$ and so the image via $\vartheta_{d-2}$ of $\vartheta_1^{-1}(x) \times (\times_{d-2}\Sigma)$ is the image of a smooth map from a manifold of dimension 4d - 8 into one of dimension 4d - 8. Let $\mathfrak{w} \in \times_{d-2} X$ denote a regular value. Then $\vartheta_{d-2}^{-1}(\mathfrak{w})$ is a discrete set with no accumulation points such that a subvariety, C, that is parametrized by this set obeys the condition set forth in Item c) on kernel($D_C$). It follows from Proposition 3.1 that $\vartheta_{d-2}^{-1}(\mathfrak{w})$ is finite if $\mathfrak{w}$ is a regular point and if there are no subvarieties in J's version of $\mathcal{M}_{e,g-2}$ that contain x and all entries of $\mathfrak{w}$. Given a regular value x, all of these conditions are satisfied by the points in residual subset of $\times_{d-2} X$. Thus, the requirements of Item c) from the second bullet of Proposition 3.7 are met if x is not a critical value for $\vartheta_1$ on $\mathcal{P}|_J \times \Sigma$.

The case when x is a critical value is handled using the same strategy as that used in Part 2 to handle the critical values for the genus g - 2 version of $\vartheta_1$. To elaborate, note that the genus g-1 version of (A.27) also holds in this case. But for minor notational changes, the proof uses the same arguments that were used to prove the analogous claim in the preceding Parts 3-7. Granted that J comes from this set, and so the critical locus of $\vartheta_1$ on $\mathcal{P}|_J \times \Sigma$ is a 3-dimensional image variety, there exists a residual subset $\mathfrak{C} \subset \times_{d-2} X$ with the following property: Take the point $\mathfrak{w}$ from $\mathfrak{C}$, and there are no points of the form $((\varphi, j, J), p, (p_1, \ldots, p_{d-2}))$ in $\vartheta_{d-2}^{-1}(\mathfrak{w})$ where $((\varphi, j, J), p)$ is a critical point of $\vartheta_1$. This understood, the analysis from the preceding paragraph can also be applied to points in the $\vartheta_1$-inverse image of a given critical value x. The conclusions in this case are the same as those in the preceding paragraph but for the following requirement: The residual set in $\times_{d-2} X$ for Item c of Proposition 3.7's second bullet must be subset of the residual set $\mathfrak{C}$.

*Part 9*: What follows here explains why, given $x \in X$, there exists a residual subset in $\times_{d-2} X$ whose points satisfy the requirements set forth by Item d) of Proposition 3.7's second bullet. To this end, first choose $J \in \mathcal{J}$, a regular value for the map from $\mathcal{P}$ to $\mathcal{J}$ such that Items b) and c) of the second bullet hold; and in particular take J from the intersection of the genus g-1 and genus g-2 versions of the residual set that is described in (A.27).

Reintroduce the space $\mathcal{P}_\diamond \subset \mathcal{P}|_J \times \Sigma$ from the proof of Proposition 3.6. Recall that the projection map to from $\mathcal{P}|_J \times \Sigma$ to $\mathcal{P}|_J$ restricts to $\mathcal{P}_\diamond$ as a 2-1 covering map. Define $\vartheta_{\diamond 1} : \mathcal{P}_\diamond \to X$ so as to send any given $((\varphi, J, j), p_\diamond) \in \mathcal{P}_\diamond$ to $\varphi(p_\diamond)$. If $x \in X$ and if a given $((\varphi, J, j), p_\diamond) \in \vartheta_{\diamond 1}^{-1}(x)$ is not a critical point of $\vartheta_{\diamond 1}$, then $\vartheta_{\diamond 1}^{-1}(x)$ is a smooth, 2d-6 dimensional manifold in a neighborhood of $((\varphi, J, j), p_\diamond)$. Let $\mathcal{U}$ denote such a neighborhood. The map $\vartheta$ in (A.26) restricts to $\mathcal{U} \times (\times_{d-2} \Sigma)$ as a map from a space of dimension 4d-10 to one of dimension 4d-8. This being the case, there is a dense, open set of points in $\times_{d-2} X$ that are not in the image via $\vartheta_{\diamond 1}$ of $\mathcal{U} \times (\times_{d-2} \Sigma)$. This observation establishes Item d) when x is not a critical value of the map $\vartheta_{\diamond 1}$.

Meanwhile, the purely cosmetic modifications to the arguments in Parts 3-7 prove the following: If $J \in \mathcal{J}$ is from a suitable residual set in $\mathcal{J}$, then the set of critical points of $\vartheta_{\diamond 1}$ is 3-dimensional image variety in $\mathcal{P}_\diamond$. Let $\mathcal{C}_\diamond$ denote the manifold that maps to this image variety and let $f_\diamond$ denote the corresponding map. The composition of the map $\vartheta$ in (A.26) with the product of $f_\diamond$ and the identity map on $(\times_{d-2} \Sigma)$ gives a map from a manifold of dimension at most 2d-1 into a manifold of dimension 4d-8. As a consequence, the complement of the image is a residual set. This understood, if J is chosen from a suitable residual set in $\mathcal{J}$, then the observations in the preceding paragraph can be applied to the case when x is a critical value of $\vartheta_{\diamond 1}$. For such J, there is, for any given $x \in X$, a residual subset of $\times_{d-2} X$ whose elements obey the requirements that are set forth in Items b), c) and d) from the second bullet of Proposition 3.7.

*Part 10*: This part of the proof considers Item e) from the second bullet of Proposition 3.7. To this end, introduce the versions of the spaces $\mathcal{P}$ and $\mathcal{P}^m$ that are defined in Part a) in this appendix using the class e and genus g. Define $\vartheta_1$ and $\vartheta_{d-2}$ to denote the respective maps from $\mathcal{P}^m \times \Sigma \times (\times_{d-2} \Sigma)$ to X and $\times_{d-2} X$ that send a given element $((\varphi, j, J), \mathfrak{p} = (p, p_1, \ldots, p_{d-2}))$ to $\varphi(p)$ and to $(\varphi(p_1), \ldots, \varphi(p_{d-2}))$.

To motivate what is to come, remark that if J is a regular value for the projection from $\mathcal{P}$ to $\mathcal{J}$, and if $\mathfrak{w} \in \times_{d-2} X$ is a regular value for the map $\vartheta_{d-2}$ on $\mathcal{P}|_J \times \Sigma \times (\times_{d-2} \Sigma)$, then $\vartheta_{d-2}^{-1}(\mathfrak{w})$ maps onto the space $\mathcal{M}^m_X$ via the rule that sends $((\varphi, j, J), (p, p_1, \ldots, p_{d-2}))$ to the pair $(\varphi(\Sigma), \varphi(p))$. Such a point maps to $\mathcal{Z}^m$ if it is a critical point of $\vartheta_1$. This is to say that there is a subspace of the tangent space to $\vartheta_{d-2}^{-1}(\mathfrak{w})$ of dimension at least 3 that is annihilated by the differential of $\vartheta_1$. Note in this regard that a tangent vector to $\vartheta_{d-2}^{-1}(\mathfrak{w})$ at $((\varphi, j, J), \mathfrak{p})$ consists of a data set of the form $((\varphi', j'), \mathfrak{p}')$ where $\varphi'$ is a section of $\varphi^*TX, j' \in V$ and $\mathfrak{p}' \in (\times_{d-1} T\Sigma)_{\mathfrak{p}}$. The pair $(\varphi', j')$ must obey (A.6). In addition, $(\varphi', \mathfrak{p}')$ are constrained at entry of $\mathfrak{p}$ as follows: Write $\mathfrak{p}' = (u, u_1, \ldots, u_{d-2})$. For $k \in \{1, \ldots, d-2\}$,

$$\varphi'|_{p_k} + \varphi_* u_k = 0 \,. \tag{A.47}$$

The differential of $\varphi_1$ annihilates such a tangent vector if

$$\varphi'(p') + \varphi_* v' = 0 \,. \tag{A.48}$$

To start the verification of Item 3), for $n = 3$ or 4, let $\mathcal{G}_n \to \mathcal{P}^m$ denote restriction from $(\mathcal{E}_* \times \mathcal{J}^m \times V_1) \times (\times_{d-1} \Sigma)$ of the fiber bundle of n dimensional, orthonormal frames in $(T\mathcal{E}_* \oplus V) \oplus T(\times_{d-1} \Sigma)$. This bundle has an action of SO(n). View this bundle as sitting over $\mathcal{P}^m \times (\times_{d-1} \Sigma)$. Define a fiber preserving map from this bundle $\mathcal{G}_n$ to the vector bundle $\oplus_n (\mathcal{V} \oplus (\vartheta_1^*TX \oplus \vartheta_{d-2}^*T(\times_{d-2}X)))$ as follows: At $((\varphi, J, j), \mathfrak{p}) \in \mathcal{P}^m \times (\times_{d-1} \Sigma)$, the map sends $p = (((\varphi_1, j_1), v_1), \ldots, ((\varphi_n, j_n), v_n))$ to the element whose component in the i'th summand has $\vartheta_1^*TX$ part equal to $\varphi_i|_p + \varphi_* v_i$ and has $\vartheta_{d-2}^*T(\times_{d-2} X)$ part equal to the d-2 tuple $(\varphi_i|_{p_1} + \varphi_* v_{i,1}, \ldots, \varphi_i|_{p_{d-2}} + \varphi_* v_{i,d-2})$ where the notation is such that $v_i \in T(\times_{d-1} \Sigma)$ is written as $(v_i, v_{i,1}, \ldots, v_{i,d-2})$. Meanwhile, the component in the $\mathcal{V}$ summand is the $(\varphi_i, j_i)$ version of the section of $\varphi^*T_{J(1,0)}X \otimes T_j^{0,1}\Sigma$ that is defined by the left hand side of (A.6). Use $\mathfrak{X}$ to denote this map.

Introduce as notation $\mathcal{G}_{n,0}$ to denote $\mathcal{G}_n$'s restriction as a fiber bundle over the subspace in $\mathcal{P}^m \times (\times_{d-1} \Sigma)$ where the entries in the factor $(\times_{d-1} \Sigma)$ are distinct. With very minor and mostly notational modifications, the arguments in Parts 3-7 from this part of the appendix prove that the differential of $\mathfrak{X}$ along $\mathfrak{X}^{-1}(0)$ in $\mathcal{G}_{n,0}$ is surjective. This being the case, $\mathfrak{X}^{-1}(0) \cap \mathcal{G}_{n,0}$ is a manifold as is its SO(n) quotient. Note in this regard that the constraint that $\mathfrak{p}$ have distinct entries is needed so as to employ the Green's function construction in Parts 6 and 7 near each of the d - 1 entries of $\mathfrak{p}$.

The map from $\mathfrak{X}^{-1}(0) \cap \mathcal{G}_{n,0}$ to $\mathcal{P}^m$ that is induced by the various projections is a Fredholm map when m is large, and so there is a residual set of regular values when m is large. The arguments used in the proof of Proposition A.5 can be applied to prove that there is a residual set of regular values in $\mathcal{P}$ as well. Let $J \in \mathcal{J}$ denote a regular value which is also a regular value for the map from $\mathcal{P}$ to $\mathcal{J}$. Then the space $(\mathfrak{X}^{-1}(0) \cap \mathcal{G}_{n,0})|_J$ is a smooth manifold with a free SO(n) action whose quotient is a smooth manifold of dimension is 4d-2 + n(2-n). Use $\mathfrak{Z}_J$ to denote this quotient.

Fix such a regular value $J \in \mathcal{J}$ that is also a regular value for the the map from $\mathcal{P}$. Let $\mathfrak{Z}_J$ denote the quotient of $(\mathfrak{X}^{-1}(0) \cap \mathcal{G}_{n,0})|_J$ by the afore-mentioned SO(n) action. Use $\vartheta_{3,d-2}: \mathfrak{Z}_J \to \times_{d-2} X$ to denote the map that is induced from the map $\vartheta_{d-2}$. Let $\mathfrak{w}$ denote a regular value for this map. The inverse image via $\vartheta_{3,d-2}$ of $\mathfrak{w}$ is empty if n = 4 and it is a smooth, 3-dimensional manifold if n = 3. Denote this manifold by $\mathfrak{Z}_J^{\mathfrak{w}}$. This manifold $\mathfrak{Z}_J^{\mathfrak{w}}$ maps to J's version of $\mathcal{M}^{\mathfrak{w}}{}_X$; this is because the space $\mathcal{M}^{\mathfrak{w}}{}_X$ is the image of $\vartheta_{d-1}{}^{-1}(\mathfrak{w}) \subset \mathcal{P}|_J \times (\times_{d-1}\Sigma)$ via the map that sends $((\varphi, J, j), (p, p_1, \ldots, p_{d-2}))$ to the pair $(\varphi(\Sigma), \varphi(p)) \subset \mathcal{M}^{\mathfrak{w}} \times X$. The latter map is denoted in what follows by $\mathfrak{f}^{\mathfrak{w}}: \mathfrak{Z}_J^{\mathfrak{w}} \to \mathcal{M}^{\mathfrak{w}}{}_X$. The image $\mathfrak{f}^{\mathfrak{w}}(\mathfrak{Z}_J^{\mathfrak{w}}) \subset \mathcal{M}^{\mathfrak{w}}{}_X$ is, by design, the $\mathcal{Z}^{\mathfrak{w}}{}_X$ portion of the critical point locus of the map $\pi^{\mathfrak{w}}{}_X: \mathcal{M}^{\mathfrak{w}}{}_X \to X$. Note in particular that the composition of this map $\mathfrak{f}^{\mathfrak{w}}$ with the map $\pi^{\mathfrak{w}}{}_{\mathcal{M}}$ from $\mathcal{M}^{\mathfrak{w}}{}_X$ to $\mathcal{M}^{\mathfrak{w}}$ maps $\mathfrak{Z}_J^{\mathfrak{w}}$ onto $\mathcal{Z}^{\mathfrak{w}}$ and thus gives $\mathcal{Z}^{\mathfrak{w}}$ the structure of a 3-dimensional image variety.

Given that $\mathfrak{w}$ can be chosen from a residual set in $\times_{d-2} X$, it follows that this last space has a residual set whose points satisfy the conditions that are set forth in Items b)-e) from the second bullet of Proposition 3.7.

*Part 11*: What follows here and in Parts 12 and 13 proves Item f) from the second bullet of Proposition 3.7. To this end, suppose that $J \in \mathcal{J}$ is chosen so as to be a regular value of the maps from $\mathfrak{X}^{-1}(0) \cap \mathcal{G}_{3,0}$, $\mathfrak{X}^{-1}(0) \cap \mathcal{G}_{4,0}$ and $\mathcal{P}$. The space $(\mathfrak{X}^{-1}(0) \cap \mathcal{G}_{3,0})|_J \times \Sigma$ maps to X by the rule that sends a point $(((\varphi, J, j), (p, p_1, \ldots, p_{d-2})), p')$ to $\varphi(p')$. As this map is SO(3) invariant, so it descends as a map, $\vartheta_{3,1}: \mathfrak{Z}_J \times \Sigma \to X$. Let $\mathfrak{Z}_{J,x} \subset \mathfrak{Z}_J \times \Sigma$ denote the inverse image via $\vartheta_{3,1}$ of a given point $x \in X$. Use $\vartheta_{3,d-2}$ to also denote the map from $\mathfrak{Z}_{J,x}$ to $\times_{d-2} X$ that sends $((\varphi, J, j), (p, p_1, \ldots, p_{d-2}), p')$ to $(\varphi(p_1), \ldots, \varphi(p_{d-2}))$. If x is a regular value for $\vartheta_{3,1}$, then $\mathfrak{Z}_{J,x}$ is a smooth manifold of dimension 4d-7. It follows as a consequence that there is a residual set in $\times_{d-2} X$ of the following sort: If $\mathfrak{w}$ is from this set, then $\vartheta_{3,d-2}{}^{-1}(\mathfrak{w})$ is a smooth, 1-dimensional manifold. The image of this manifold via the composition $\pi^{\mathfrak{w}}{}_X \circ \mathfrak{f}^{\mathfrak{w}}$ is the space $\mathcal{Z}^{(x,\mathfrak{w})}$.

The case when x is a critical value of $\vartheta_{3,1}$ is dealt with in a manner that is explained next. The upcoming Part 12 establishes

**Lemma A.14**: *There is a residual set of regular values in $\mathcal{J}$ which is characterized as follows: Take J from this set. Fix $x \in X$. There exists a residual set of regular points in $\times_{d-2} X$ whose points have the properties listed next. Let $\mathfrak{w}$ denote such a point. Then the locus $\vartheta_{3,d-2}{}^{-1}(\mathfrak{w}) \subset \mathfrak{Z}_J \times \Sigma$ contains no point of the form $(((\varphi, J, j), (p, p_1, \ldots, p_{d-2})), p')$ with $\varphi(p') = x$ and $p' = p$.*

With Lemma A.14 in hand, introduce $(\mathfrak{Z}_J \times \Sigma)' \subset \mathfrak{Z}_J \times \Sigma$ to denote the set of points of the form $(((\varphi, J, j), (p, p_1, \ldots, p_{d-2})), p')$ with $p'$ neither p nor any $k \in \{1, \ldots, d-2\}$ version of $p_k$. It is proved in the upcoming Part 13 that there is a residual set in $\mathcal{J}$ such that if J comes from the latter, then the critical locus of $\vartheta_{3,1}$ in $(\mathfrak{Z}_J \times \Sigma)'$ is an image variety of dimension 3. This being the case, and given Lemma A.14, there exists a residual set in $\times_{d-2} X$ with the following properties: If $\mathfrak{w}$ comes from this set, then there are no critical points of $\vartheta_{3,1}$ in $\vartheta_{3,d-2}{}^{-1}(\mathfrak{w})$ (which is a subset of $\mathfrak{Z}_J \times \Sigma$) other than those that have the form $(((\varphi, J, j), (p, p_1, \ldots, p_{d-2})), p')$ with $p' = p_k$ for any $k \in \{1, \ldots, d-2\}$. Given this, then the analysis employed above when x is a regular value of $\vartheta_{3,1}$ can be repeated to see that there is a residual subset in $\times_{d-2} X$ with the following property: If $\mathfrak{w}$ is from this set, then $\vartheta_{3,d-2}^{-1}(\mathfrak{w})$ is a 1-dimensional image variety, and thus so is $\mathcal{Z}^{(x,\mathfrak{w})}$.

The conclusions of this paragraph and the first paragraph of Part 11 imply that there is a residual set in $\mathcal{J}$ with the following property: Given any $x \in X$, there is a residual set $\times_{d-2} X$ (to be denoted by $\mathcal{X}_x$) whose points satisfy the conditions that are set forth in Items b)-f) of Proposition 3.7.

As a parenthetical remark, note that Lemma A.14 has the following immediate corollary:

**Lemma A.15**: *There is a residual set of regular values in $\mathcal{J}$ which is characterized as follows: Take J from this set. Fix $x \in X$. There exists a residual set of regular points in $\times_{d-2} X$ such that x is not a critical point of $\pi^{\mathfrak{w}}_X \colon \mathcal{M}^{\mathfrak{w}}_X \to X$ if $\mathfrak{w}$ comes from this set. In particular, for such $\mathfrak{w}$, the space $\mathcal{M}^{(x,\mathfrak{w})}$ is a 2-dimensional smooth manifold.*

*Part 12*: This part proves Lemma A.14. To motivate what is to follow, digress for a moment for an observation concerning transversality of maps. To set the stage, suppose that Y is a smooth 4-dimensional manifold and that $\{(Z_i, \phi_i)\}_{1 \le i \le 5}$ is a set whose elements have the form $(Z, \phi)$ with Z a smooth 3-manifold and $\phi \colon Z \to Y$ a smooth map. Dimension counting leads to the following observation: If the maps $\{\phi_i\}_{1 \le i \le 5}$ are suitably generic, then $\cap_{1 \le i \le 5} \phi_i(Z_i) = \emptyset$.

To apply this last observation, introduce now $\mathfrak{Z}$ to denote $(\mathfrak{X}^{-1}(0) \cap \mathcal{G}_{3,0})/SO(3)$. Then use $\mathfrak{Z}_{5,X} \subset \times_5 \mathfrak{Z}$ to denote the subspace of 5-tuples with the properties listed below. To set the notation for the upcoming list, write a given element in the product $\times_5 \mathfrak{Z}$ as $\{((\varphi_\alpha, J_\alpha, j_\alpha), \mathfrak{p}_\alpha, [p_\alpha])\}_{1 \le \alpha \le 5}$. Write each $\alpha \in \{1, \ldots, 5\}$ version of $\mathfrak{p}_\alpha$ as $(p_\alpha, q_\alpha)$ with $p_\alpha \in \Sigma$ and $q_\alpha \in \times_{d-2} \Sigma$. With the notation now set, here are the properties:

- $J_\alpha = J_\beta$ *for each pair* $\alpha, \beta \in \{1, \ldots, 5\}$.
- $\varphi_\alpha(p_\alpha) = \varphi_\beta(p_\beta)$ *for each pair* $\alpha, \beta \in \{1, \ldots, 5\}$.
- *If* $\alpha, \beta \in \{1, \ldots, 5\}$ *are distinct, then* $q_\alpha$ *and* $q_\beta$ *have no components in common.*
- *If* $\alpha, \beta \in \{1, \ldots, 5\}$ *are distinct, then* $\varphi_\alpha(\Sigma) \ne \varphi_\beta(\Sigma)$.

(A.49)

The space $\mathfrak{Z}_{5,X}$ maps to $\mathcal{J}^{\mathfrak{m}}$ with the image of a given point $\{(\varphi_\alpha, J, j_\alpha), \mathfrak{p}_\alpha, [p_\alpha])\}$ being J. It also maps to X with the image of this point given by the common value of $\{\varphi_\alpha(p_\alpha)\}_{1 \le \alpha \le 5}$. The next lemma says what is needed about $\mathfrak{Z}_{5,X}$.

**Lemma A.16**: *The space $\mathfrak{Z}_{5,X}$ sits in $\times_5 \mathfrak{Z}$ as a smooth manifold such that the corresponding map to $\mathcal{J}^{\mathfrak{m}}$ is a Fredholm map with index $5(4d - 8) - 1$. Moreover, this map has a residual set of regular values in $\mathcal{P}$.*

This lemma is proved momentarily.

***Proof of Lemma A.14***: Choose $J \in \mathcal{J}$ to be a regular value for Lemma A.16's map to $\mathcal{J}$, and a regular value for the maps from $\mathfrak{X}^{-1}(0) \cap \mathcal{G}_{3,0}$, $\mathfrak{X}^{-1}(0) \cap \mathcal{G}_{4,0}$, and $\mathcal{P}$ to $\mathcal{J}$. In addition, $J$ must be such that the assertion of the first bullet from Proposition 3.7 holds. The manifold $\mathfrak{Z}_{5,X}|_J$ maps to $\times_5 (\times_{d-2} X)$ via the restriction of $\times_5(\vartheta_{3,d-2})$. As the range space has dimension 1 greater than the domain, the inverse image of a regular value is empty. This has the following implication: Fix $x \in X$ and there is a regular point $\mathfrak{w} \in \times_{d-2} X$ such that $\vartheta_{3,d-2}^{-1}(\mathfrak{w}) \subset \mathfrak{Z}_J$ has no point of the form $((\varphi, J, j), (p, p_1, \ldots, p_{d-2}))$ with $\varphi(p) = x$. Indeed, there must be a residual set of such points $\mathfrak{w}$. To see why this is, let $U \subset \times_{d-2} X$ denote the subset of regular points $\mathfrak{w}$ such that $\vartheta_{3,d-2}^{-1}(\mathfrak{w})$ contains a point with $\varphi(p) = x$. Then the set $\times_5 U$ consists of points $(\mathfrak{w}_1, \ldots, \mathfrak{w}_5)$ that are not regular values for the restriction to $\mathfrak{Z}_{5,X}|_J$ of the map $\times_5(\vartheta_{3,d-2})$. The complement of this set in $\times_5(\times_{d-2} X)$ consists of points of the form $(\mathfrak{w}_1, \ldots, \mathfrak{w}_5)$ such that at least one of them is not in U. Thus, it is the union of 5 sets, the first of the form $((\times_{d-2} X) - U) \times (\times_4(\times_{d-2} X))$ and the others obtained from the latter via a cyclic permutation of its entries. According to Lemma A.16, their union must contain a residual set. This is possible only if $(\times_{d-2} X) - U$ contains such a set. This last conclusion implies what is asserted by Lemma A.14.

***Proof of Lemma A.16***: The space $\mathfrak{Z}_{5,X}$ is a smooth submanifold of $\times_5 \mathfrak{Z}$ as described by the lemma if the map from $\times_5 \mathfrak{Z}$ to $(\times_5 (\mathcal{J}^m \times X))$ is transversal to the full diagonal at all points in $\mathfrak{Z}_{5,X}$. It is argued momentarily that this is the case, and the argument implies that the map to $\mathcal{J}^m$ is Fredholm with the asserted index. The arguments used for Proposition A.5 can be employed yet again to see that this map has a residual set of regular values in the space $\mathcal{J}$.

The transversality proof has three steps. The first two supply some background.

<u>Step 1</u>: A tangent vector to a given point $((\varphi, J, j), \mathfrak{p}, [p]) \in \mathfrak{Z}$ lifts to $\mathfrak{X}^{-1}(0) \cap \mathcal{G}_{n,0}$ as a data set $((\varphi', J', j'), \mathfrak{p}', q)$ of the following sort: First, $(\varphi', J', j')$ obeys the equation in the second bullet in (A.29). In addition $\varphi'$ and the tangent vector $\mathfrak{p}' = (u, u_1, \ldots, u_{d-2}) \in T(\times_{d-2}\Sigma)|_{\mathfrak{p}}$ are constrained at each entry of $\mathfrak{p}$ via (A.47) and (A.48). Second, what is written as $q$ has the form $q = (((\lambda_1, \mathfrak{k}_1), x_1), ((\lambda_2, \mathfrak{k}_2), x_2), ((\lambda_3, \mathfrak{k}_3), x_3))$ where each $(\lambda_i, \mathfrak{k}_i)$ obeys the $g_i = 0$ version of the third bullet of (A.29) with it understood that the corresponding $(\varphi_i, j_i)$ comes from the i'th entry of $p$. Here, this i'th entry of $p$ should be written as $((\varphi_i, j_i), v_i)$ with $v_i = (v_i, v_{i1}, \ldots v_{id-2})$ denoting a d-1 tuple such that any given entry is a tangent vector to $\Sigma$ at the corresponding entry of $\mathfrak{p} = (p, p_1, \ldots, p_{d-2})$. The pair $(\varphi_i, j_i)$ obeys the equation in the top bullet of (A.29), and $\varphi_i, v_i$ are constrained at each entry of $\mathfrak{p}$ via the version of (A.47) and (A.48) that has $\varphi_i$ in lieu of $\varphi'$, and the entries of $v_i$ in lieu of $(u, u_1, \ldots, u_{d-2})$. Each $x_i = (x_i, x_{i1}, \ldots, x_{id-2})$ from $\mathfrak{q}$ is a d-1 tuple where any given component is tangent vector to $\Sigma$ at the corresponding component of $\mathfrak{p}$. The pair $(\lambda_i, x_i)$ with $v_i$ and $\mathfrak{p}'$ are constrained at the entries of $\mathfrak{p}$ to obey

- $\lambda_i|_p + \varphi_* x_i + \nabla_u \varphi_i|_p + (\nabla_u \varphi_*) v_i = 0$
- $\lambda_i|_{p_k} + \varphi_* x_{ik} + \nabla_{u_k} \varphi_i|_{p_k} + (\nabla_{u_k} \varphi_*) v_{ik} = 0$ *for each* $k \in \{1, \ldots, d-2\}$.

(A.50)

The vectors that comprise $\mathfrak{q}$ must also be orthogonal to the span of $p$.

Step 2: Fix a point $\{((\varphi_\alpha, J_\alpha, j_\alpha), \mathfrak{p}_\alpha = (p_\alpha, \mathfrak{p}'_\alpha), [p_\alpha])\}_{1 \leq \alpha \leq 5} \in \times_5 \mathfrak{Z}$ denote a point in $\mathfrak{Z}_{5,X}$. Let J denote their common $\mathcal{J}^m$ component and let $x \in X$ denote the one point comprising the set $\{\varphi_\alpha(p_\alpha)\}_{1 \leq \alpha \leq 5}$. The map to the full diagonal in $\times_5 (\mathcal{J}^m \times X)$ at this given point is transversal provided that the following is true: Let $(J_1, \ldots, J_4)$ denote any given 4-tuple in $T\mathcal{J}^m|_J$ and let $(y_1, \ldots, y_4)$ denote any given 5-tuple in $TX|_x$. There must exist a tangent vector to $\times_5 \mathfrak{Z}$ such that any given $\alpha \in \{1, \ldots, 4\}$ entry has J' component equal to $J_\alpha + J$ where J is the $\alpha = 5$ entry for the J' component. Meanwhile the $(\varphi', u)$ entry for any $\alpha \in \{1, \ldots, 4\}$ must be such that

$$\varphi'|_{p_\alpha} + \varphi_{\alpha *} u = y_\alpha + y$$

(A.51)

where y is given by the $\alpha = 5$ version of the left hand side of this last expression.

Step 3: To see that such a vector can be found, note first that each $\varphi_\alpha$ is an embedding, so the component of $y_\alpha + y$ along $\varphi_\alpha(\Sigma)$ can be realized as $\varphi_{\alpha *} u$ for some $u \in T\Sigma|_{p_\alpha}$. This understood, at issue is the normal projection of the various versions of (A.51) and the various constraints on the J' entries. These can all be met for the following reason: Going back to what is said in Parts 3-7, the constraints on the J' component and the normal projections of (A.47) can all be met by a suitable choice of J because the constraints involve only the restriction of the J' component to a neighborhood of each $\varphi_\alpha(\Sigma)$ that is disjoint from the point x. As these curves are all distinct, there is no obstruction to setting J as needed. The Green's function construction from Parts 6 and 7 can be used near each $p_\alpha$ with the corresponding version of the middle equation in (A.29) so as to satisfy each of the $\alpha \in \{1, \ldots, 5\}$ versions of (A.51). Here again, J is constrained only near each of the 5 images of $\Sigma$, but not at precisely the point, x, where they intersect.

*Part 13*: To motivate the up coming constructions, suppose that $J \in \mathcal{J}$ is a regular value for the map from $\mathfrak{X}^{-1}(0) \times \mathcal{G}_{3,0}$, $\mathfrak{X}^{-1}(0) \times \mathcal{G}_{3,0}$ and from $\mathcal{P}$. A tangent vector to a point $((\varphi, J, j), \mathfrak{p}, [p], p') \in \mathfrak{Z}_J \times \Sigma$ lifts to $(\mathfrak{X}^{-1}(0) \cap \mathcal{G}_{3,0})|_J \times \Sigma$ as a data set $((\varphi', j'), \mathfrak{p}', q, v')$ such that $(((\varphi', j'), \mathfrak{p}', q)$ are as describe in Step 1 of the proof of Lemma A.16 with the proviso that $J' = 0$. Meanwhile, $v'$ is a tangent vector to $p'$. The differential of the map $\vartheta_{3,1}$ sends such a tangent vector to the vector $\varphi'(p') + \varphi_* v' \in (\varphi^* TX)|_p$. Since $\mathfrak{Z}_J \times \Sigma$ has dimension $2d-3$, a point in $\mathfrak{Z}_J \times \Sigma$ is a critical point of $\vartheta_{3,1}$ if and only if there is a $2d-6$ dimensional subspace in its tangent space whose vectors obey (A.48).

To see about the existence of such a space, return to a very large m version of the space $\mathfrak{Z} = (\mathfrak{X}^{-1}(0) \cap \mathcal{G}_{3,0})/SO(3)$. Introduce $(\mathfrak{Z} \times \Sigma)_0 \subset \mathfrak{Z} \times \Sigma$ to denote the subspace of points $((\varphi, J, j), \mathfrak{p}, [p], p')$ where $p'$ is distinct from all entries of $\mathfrak{p}$, and where $J$ is a regular value for both the projection from $\mathfrak{X}^{-1}(0) \cap \mathcal{G}_{3,0}$, $\mathfrak{X}^{-1}(0) \cap \mathcal{G}_{4,0}$ and from $\mathcal{P}^m$ to $\mathcal{J}^m$. Note that this is an open set in $\mathfrak{Z} \times \Sigma$.

Let $\mathcal{K} \to (\mathfrak{Z} \times \Sigma)_0$ denote the bundle whose fiber at a point $((\varphi, J, j), \mathfrak{p}, [p], p')$ is of the space of $2d-6$ dimensional, orthonormal frames in $T(\mathfrak{Z}|_J \times \Sigma)$. The latter is a smooth manifold with a smooth action of $SO(2d-6)$. Introduce the vector bundle, $\mathfrak{E}$, over $(\mathfrak{Z} \times \Sigma)_0$ whose fiber at the given point is $\oplus_{1 \leq \alpha \leq 2d-6} \varphi^* TX$.

Use $\wp$ to denote the fiber preserving map from $\mathcal{K}$ to $\mathfrak{E}$ that is defined over the point $((\varphi, J, j), \mathfrak{p}, [p], p')$ by the rule that follow. Let $\alpha \in \{1, \ldots, 2d-6\}$ and introduce $\mathit{k}_\alpha$ to denote the index $\alpha$ component of the given frame. This component defines the corresponding index $\alpha$ summand $\varphi^* TX$ in $\mathfrak{E}$. To say more, write $\mathit{k}_\alpha$ as $((\varphi', j'), \mathfrak{p}', q, v')$. The index $\alpha$ component of $\wp$ is $\varphi'(p') + \varphi_* v'$; this the left hand side of (A.48).

As is argued momentarily, the locus $\wp^{-1}(0)$ is a manifold such that the map to $\mathcal{J}^m$ is Fredholm with index 3. If this is the case, then the arguments for Proposition A.5 can be employed to prove that there is a residual set of regular values in $\mathcal{J}$ for this map. The inverse image of such a regular value is a 3-dimensional manifold whose image in the subspace $(\mathfrak{Z}|_J \times \Sigma)_0$ of $\mathfrak{Z}_J \times \Sigma$ is the intersection of the critical locus of $\vartheta_{3,1}$ with the subset of elements of the form $(((\varphi, J, j), \mathfrak{p}, [p], p')$ such that $p'$ is not an entry of $\mathfrak{p}$.

To prove that $\wp^{-1}(0)$ is a manifold, it is sufficient to prove that the differential of $\wp$ along $\wp^{-1}(0)$ is surjective. The arguments to prove this imply directly that the corresponding map to $\mathcal{J}^m$ is Fredholm with index 3. To see about this surjectivity, note that $\wp$'s differential maps a tangent vector to $\mathcal{K}$ at a given point to the fiber of $\mathfrak{E}$ at the image of this point via the projection $\mathcal{K} \to (\mathfrak{Z} \times \Sigma)_0$. To say more, write the $\varphi'$ and $v'$ components of the index $\alpha$ component of the frame over $(((\varphi, J, j), \mathfrak{p}, [p], p')$ as $\varphi'_\alpha$ and $v'_\alpha$. Note that $\varphi'_\alpha$ has a corresponding $j'_\alpha$ from $V$. The differential of $\wp'$ involves only the component in $T\mathcal{K}$ that corresponds to a first order change in $\varphi'_\alpha$ and $v'_\alpha$. It also involves a tangent vector, $u'$, at $p'$ that gives the first order change of the latter point. To say more, denote this first order change in the data $(\varphi'_\alpha, j'_\alpha)$ and $v'_\alpha$ as $(\lambda'_\alpha, \mathit{k}'_\alpha)$ and $x'_\alpha$. Then the index $\alpha$ component of the differential of $\wp'$ gives the element

$$\lambda'_\alpha|_{p'} + \varphi_* x'_\alpha + \nabla_u \varphi'_\alpha + (\nabla_u \varphi_*) v'_\alpha \in (\varphi^* TX)|_{p'} .$$
(A.52)

Meanwhile, that data $(\lambda'_\alpha, \mathfrak{k}'_\alpha)$ is constrained by the requirement that the index $\alpha$ component of the frame over $(((\varphi, J, \mathfrak{j}), \mathfrak{p}, [p], p')$ define a tangent vector to $\mathfrak{Z}_J \times \Sigma$.

With the preceding understood, here is the key observation: Let J´ denote the component of the tangent vector to $\mathcal{K}$ that comes from a first order change in J. The requirement that J be a regular value for both the maps from $\mathcal{P}^m$ and $\mathfrak{Z}$ to $\mathcal{J}^m$ implies that J´ is *not yet constrained* by the requirement that the domain for the linear map in (A.50) consist of tangent vectors to $\mathcal{K}$. As a consequence, J´ can be chosen so the various 2d-6 versions of (A.50) can have any given value in $\oplus_{1 \le \alpha \le 2d-6} (\varphi^* TX)|_{p'}$. The argument for this employs once again the Green's function construction from Part 6 to obtain any desired value for the part of $\lambda'_\alpha|_{p'}$ that is normal to $\varphi(\Sigma)$ at p´. The details are messy, but in the end straight forward and so omitted.

### g) Proof of Proposition 3.8

Let $\mathcal{J}_{e2}$ denote the residual set from Proposition 3.6. The proof of Proposition 3.8 starts with the following lemma.

**Lemma A.17**: *There exists a residual set in $\mathcal{J}_{e2}$ which is characterized by what follows. Fix an almost complex structure J from this set to define $\mathcal{M}_{e,g}$. There is a residual subset of regular points $\mathfrak{w} \in \times_{d-2} X$ with the property that if $\mathfrak{w}$ comes from this set, and if w is an entry of $\mathfrak{w}$, then $\mathcal{Y}^w$ is a 1-dimensional image variety in $\mathcal{M}^{\mathfrak{w}}$.*

Grant this lemma for the moment. One might expect that a pair of distinct and suitably generic 1-dimensional image varieties in a 4-dimensional manifold will have no points in common. The proof of Proposition 3.8 argues that this expectation is met for distinct pairs from any given version of $\{\mathcal{Y}^w : w$ *is an entry of* $\mathfrak{w}\}$ if $\mathcal{M}_{e,g}$ is first defined using an almost complex structure from a certain residual set in $\mathcal{J}_{e2}$, and then $\mathfrak{w}$ is chosen from a certain residual set in $\times_{d-2} X$.

The proof of Proposition 3.8 including Lemma A.17 has four parts.

*Part 1*: To motivate what is to come, fix for the moment, an almost complex structure $J \in \mathcal{J}$ which is a regular value for the genus g and class e version of $\mathcal{P}$. Use J to define the space $\mathcal{M}_{e,g}$. Fix a regular point, $\mathfrak{w} \in \times_{d-2} X$. This guarantees that $\mathcal{M}^{\mathfrak{w}}$ is a smooth, 4-dimensional manifold. Let w denote an entry of $\mathfrak{w}$. To consider the structure of $\mathcal{Y}^w$, let $C \in \mathcal{M}^{\mathfrak{w}}$ denote a given curve and let $N \to C$ denote its normal bundle in X. Introduce the operator $D_C$ as defined by (2.12). The tangent space at C consists of the subspace $\ker_{C,\mathfrak{w}} = \{\eta \in \ker(D_C): \eta = 0$ *at each entry of* $\mathfrak{w}\}$. Fix an adapted coordinate chart, $(z, w)$, centered at $w$ so that the $w = 0$ locus is tangent to C at the origin. This identifies the fiber of N at $w$ with the span of $\frac{\partial}{\partial w}$ and it identifies this same span with the tangent space to $\mathbb{CP}^1$ at $\varphi_w(C)$. Fix a holomorphic coordinate, u, for C centered at $w$ with $du = dz$ at $w$. The differential of $\varphi_w$ at C sends a given $\eta \in \ker_{C,\mathfrak{w}}$ to the vector $\partial \eta|_w$, here viewed using these identifications as an element in $T\mathbb{CP}^1|_{\varphi_w(C)}$. Use $\varphi_{w*}$ to denote this differential. A curve $C \in \mathcal{M}^{\mathfrak{w}}$ is a critical point of $\varphi_w$ if and only if the kernel of $\varphi_{w*}$ contains a 3-dimensional, orthonormal frame.

Granted the preceding, Lemma A.17 and Proposition 3.8 are proved via an analysis of the space whose elements consist of data sets $(J, (C, \mathfrak{w}), w, \mathfrak{p})$ of the following sort: First, J is an almost complex structure. Second $(C, \mathfrak{w})$ are in J's version of $\mathcal{M}_{e,g,d-2}$. Third, $w$ is an entry of $\mathfrak{w}$. Finally, $\mathfrak{p}$ consists of an orthonormal frame $(\eta_1, \eta_2, \eta_3)$ in the space $\ker_{C,\mathfrak{w}}$ such that each frame vector is in the kernel of $\varphi_{w*}$.

*Part 2*: The space just described is studied by reintroducing the genus g and bundle e versions of the manifold $\mathcal{P}^{\mathfrak{m}}$ and $\mathcal{P}$. The use of the latter is not strictly speaking necessary, but it does allow for direct appeal to various constructions from the preceding parts of this appendix. In any event, to make contact with the problem at hand, consider now the open submanifold of points $((\varphi, J, \mathfrak{j}), \mathfrak{p}) \in \mathcal{P}^{\mathfrak{m}} \times (\times_{d-2} \Sigma)$ where J is a regular value for the map from $\mathcal{P}^{\mathfrak{m}}$ to $\mathcal{J}^{\mathfrak{m}}$ and where the point $\mathfrak{p}$ has distinct entries. For n = 3 and n = 4 define a fiber bundle, $\mathcal{Q}_n$, over this subset of $\mathcal{P}^{\mathfrak{m}} \times (\times_{d-2}\Sigma)$ as follows: The fiber over a given point $((\varphi, J, \mathfrak{j}), \mathfrak{p})$ consists of the space of orthonormal n-frames where each frame vector has the form $((\varphi´, \mathfrak{j}´), \mathfrak{p}´)$ with $(\varphi´, \mathfrak{j}´)$ and $\mathfrak{p}´ = (u_1, \ldots, u_{d-2}) \in T(\times_{d-2}\Sigma)|_{\mathfrak{p}}$ such that (A.6) and each $k \in \{1, \ldots, d-2\}$ version of (A.47) holds. Note that this fiber bundle has a free, fiber preserving action of SO(n).

To connect this with what was said in Part 1, let $((\varphi´, \mathfrak{j}´), \mathfrak{p}´)$ denote a vector as just described. Then the projection of $\varphi´$ to the pull-back via $\varphi$ of the normal bundle to the curve $C = \varphi(\Sigma)$ is a vector in $\ker_{C,\mathfrak{w}}$ where $\mathfrak{w} = (\varphi(p_1), .., \varphi(p_{d-2}))$.

Fix $\alpha \in \{1, \ldots, d-2\}$. Define a vector bundle over this same part of $\mathcal{P}^{\mathfrak{m}} \times (\times_{d-2}\Sigma)$ as follows: The fiber over the point $((\varphi, J, \mathfrak{j}), \mathfrak{p})$ is that of $(\varphi^*T_{J(1,0)}X/\varphi_*T_{\mathfrak{j}(1,0)}\Sigma) \otimes T_\mathfrak{j}^{(1,0)}\Sigma$ at the component $p_\alpha$ of P. Use $\mathfrak{N}_\alpha$ to denote the latter. Define a fiber preserving map

$$\mathfrak{f}_\alpha: \mathcal{Q}_n \to \oplus_n \mathfrak{N}_\alpha$$

(A.53)

by the following rule: For $i \in \{1, \ldots, n\}$, let $((\varphi´, \mathfrak{j}´), \mathfrak{p}´)$ denote the i'th frame vector of a frame in the fiber over a given point $((\varphi, J, \mathfrak{j}), \mathfrak{p})$. The latter defines the entry of $\mathfrak{f}_\alpha$ in the i'th summand of the fiber of $\oplus_n \mathfrak{N}_\alpha$ over $((\varphi, J, \mathfrak{j}), \mathfrak{p})$, this being the projection at $p_\alpha$ to the $(\varphi^*T_{J(1,0)}X/\varphi_*T_{\mathfrak{j}(1,0)}\Sigma) \otimes T_\mathfrak{j}^{(1,0)}\Sigma$ of the covariant derivative of $\varphi´$.

***Proof of Lemma A.17***: The locus $\mathfrak{f}_\alpha^{-1}(0)$ is a manifold if the differential of $\mathfrak{f}_\alpha$ is transversal along this locus. It is argued in Part 3 that this is indeed the case. These upcoming arguments imply directly that the map from $\mathfrak{f}_\alpha^{-1}(0)$ to $\mathcal{J}^m$ is a Fredholm map whose index is $4d - 4$ when $n = 3$ and $4d - 6$ when $n = 4$. Granted that such is the case, the arguments used to prove Proposition A.5 can be used again to prove that there is a residual set of regular values for this map in $\mathcal{J}$ if $m$ is large.

Let $J \in \mathcal{J}$ denote such a regular value. Then $\mathfrak{f}_\alpha^{-1}(0)|_J$ is a smooth manifold with a smooth action of $SO(n)$ that preserves the induced map to $\mathcal{P}_J$. Let $\mathfrak{Y}_J$ denote this quotient. The latter is a smooth manifold of dimension $4d - 7$ if $n = 3$ and of dimension $4d - 12$ if $n = 4$. Use $\vartheta_{\mathfrak{N},d-2}$ to denote the map from $\mathfrak{N}_J$ to $\times_{d-2} X$ that is obtained by composing first the tautological map to $\mathcal{P}_J \times \Sigma_{d-2}$ with the map $\vartheta_{d-2}$ from the latter space to $\times_{d-2} X$. Let $\mathfrak{w}$ denote a regular value. In the case $n = 4$, the inverse image via $\vartheta_{\mathfrak{N},d-2}$ of $\mathfrak{w}$ is a empty. In the case $n = 3$, it is a 1-dimensional submanifold of $\mathfrak{N}_J$. Let $w$ denote the index $\alpha$ entry of $\mathfrak{w}$. By design, the manifold $\vartheta_{\mathfrak{N},d-2}^{-1}(\mathfrak{w}) \subset \mathfrak{N}_J$ maps onto $\mathcal{Y}^w$.

*Part 3*: This part explains why $\mathfrak{f}_\alpha$ has transversal differential along its zero locus. To set the stage, let $((\varphi, J, \mathfrak{j}), \mathfrak{p})$ denote a point in $\mathcal{P}^m \times (\times_{d-2}\Sigma)$ where $\mathcal{Q}_n$ is defined. Let $p = \{((\varphi_i, \mathfrak{j}_i), v_i = (v_{i1}, \ldots, v_{id-2}))\}_{1 \le i \le n}$ denote a point in the fiber of $\mathcal{Q}_n$ over this given point. A tangent vector to the point in $\mathcal{Q}_n/SO(n)$ defined by $p$ lifts to the frame $p$ as a data set $((\varphi', J', \mathfrak{j}'), \mathfrak{p}', q)$ where $(\varphi', J', \mathfrak{j}')$ is tangent to $\mathcal{P}^m$ at $(\varphi, J, \mathfrak{j})$, where $\mathfrak{p}' = (u_1, \ldots, u_{d-2}) \in T(\times_{d-2}\Sigma)|_\mathfrak{p}$, and where $q$ is as follows. It has $n$ components, with the $i$'th having the form $((\lambda_i, \mathfrak{k}_i), x_i)$. The pair $(\lambda_i, \mathfrak{k}_i)$ obeys the $g_i = 0$ version of the third bullet of (A.29) with it understood that the corresponding $(\varphi_i, \mathfrak{j}_i)$ comes from the $i$'th entry of $p$. Note that the latter must obey the top equation in (A.29), and $\varphi_i$ with $v_i$, it must obey the version of (A.47) that has $\varphi_i$ in lieu of $\varphi'$ and the entries of $v_i$ in lieu of $(u_1, \ldots, u_{d-2})$. Meanwhile $x_i = (x_{i1}, \ldots, x_{id-2}) \in T(\times_{d-2}\Sigma)|_\mathfrak{p}$. The pair $(\lambda_i, x_i)$ with $v_i$ and $\mathfrak{p}'$ are constrained at the entries of $\mathfrak{p}$ via the equation that appears in the second bullet of (A.50). Finally, each of the $n$ components of $q$ is orthogonal to the span of $p$.

Suppose now that $((\varphi, J, \mathfrak{j}), \mathfrak{p}, p)$ is a zero of $\mathfrak{f}_\alpha$. The differential of $\mathfrak{f}_\alpha$ is surjective at this zero if the following is true: Let $N$ denote the normal bundle to $\varphi(\Sigma)$. Fix a vector in $\oplus_n (\varphi^*(N) \otimes T_\mathfrak{j}^{1,0}\Sigma)|_{u_\alpha}$. Then, there exists a tangent vector of the sort just described which is such that the projection to $\varphi^*(N) \otimes T_\mathfrak{j}^{1,0}\Sigma$ at $u_\alpha$ of each $i \in \{1, \ldots, n\}$ version of $\partial \lambda_i$ is equal to $i$'th component of the given vector in $\oplus_n (\varphi^*(N) \otimes T_\mathfrak{j}^{1,0}\Sigma)|_{u_\alpha}$.

The proof that such is the case uses a version of the Green's function construction from Part 6 of the proof of Proposition 3.7. What follows gives the details. The first observation is that the assumption that $J$ is a regular value for the map from $\mathcal{P}^m$ to $\mathcal{P}$ implies that there no constraints on J´ are required to solve the relevant versions of (A.29), (A.50) and (A.47); this because the operator depicted in (A.30) is invertible. With the preceding understood, the plan is to chose J´ so as to achieve any given desired value for the normal bundle projection of $(\partial \lambda_i)_{1 \leq i \leq n}$ at $u_\alpha$. To this end, the first step is to restrict J´ so that the equations in the first and third bullets of (A.29) are given by what is written in (A.31) and so that the equation in the second bullet of (A.29) is solved by taking $(\varphi´, j´) = 0$. Take $\mathfrak{p}´ = 0$ to solve (A.47). A choice now for J´ must be made so that the corresponding solution to the equation in the second bullet of (A.31) obeys the equation in the second bullet of (A.50) and is such that the normal bundle projection of $\partial \lambda_i$ has the desired value at $u_\alpha$.

To achieve the goal just described, fix a j-holomorphic coordinate for $\Sigma$ centered at $u_\alpha$. This coordinate, u, takes values in a small radius disk centered about the origin. Introduce the notation from Part 6 of the proof of Proposition 3.7. Note in what follows that the normal bundle projection in this Part 6 is denoted by $\Pi$, and $\eta_i$ is used to denote each $i \in \{1, \ldots, n\}$ version of the normal bundle projection of $\varphi_i$. Part 6 of the proof of Proposition 3.7 defines a finite set, W, in the u-coordinate disk. Copy what is done there to define here such a set of points, all very near the origin but none at $u = 0$. Use $W_0$ to denote this set. Fix $d > 0$ but much less then the distance from any point in $W_0$ to the origin, and much less then the distance between any two distinct elements in $W_0$. Let $W_d$ denote the set obtained from $W_0$ by translating each element a distance d along the direction of the real axis in $\mathbb{C}$. Set $W = W_0 \cup W_d$.

Fix $\Delta > 0$ but much smaller than the distance between any point in W and 0, and much smaller than the distance between distinct points in W. Take J as in (A.33). The analysis in Part 6 of the proof of Proposition 3.7 can be repeated in this context to prove that

- $\Pi(\lambda_i)_{J(1,0)}|_0 = \Delta^2 \left( \sum_{u \in W_0} \frac{1}{u} (J´\cdot \bar{\partial} \eta_i)|_u + \frac{1}{u+d} (J´\cdot \bar{\partial} \eta_i)|_{u+d} + \mathfrak{r}_i \right)$ ,
- $\Pi(\partial \lambda_i)_{J(1,0)}|_0 = - \Delta^2 \left( \sum_{u \in W_0} \frac{1}{u^2} (J´\cdot \bar{\partial} \eta_i)|_u + \frac{1}{(u+d)^2} (J´\cdot \bar{\partial} \eta_i)|_{u+d} + \mathfrak{r}_i´ \right)$ ,

(A.54)

where $|\mathfrak{r}_i|$ and $|\mathfrak{r}_i´|$ are bounded by what is written on the right hand side of (A.35).

What is said in Part 6 of the proof of Proposition 3.7 implies the following: Given $\{J´|_u\}_{u \in W_0}$, then the various $u \in W_0$ versions of $J´|_{u+d}$ can be chosen so that the following conditions hold when $\Delta$ is small. First, each $i \in \{1, \ldots, n\}$ version of what is written on the right hand side of the top line of (A.54) is zero. Moreover, this can be achieved by so that

$$\frac{1}{u+d}(J´\cdot \bar{\partial} \eta_i)|_{u+d} = -\frac{1}{u}(J´\cdot \bar{\partial} \eta_i)|_u + \mathfrak{e}_{u,i} \text{ for each } u \in W_0 ,$$

(A.55)

where $|e_{u,i}| \leq c\Delta \sum_{u \in W_0} |J'|_u|$ with c a constant that depends on $W_0$ but not on d nor on $\Delta$ if $\Delta$ is much less than d.

Granted this choice for $\{J'|_{u+d}\}_{u \in W_0}$, then the right hand side of the bottom bullet in (A.54) has the form

$$\Pi(\partial \lambda_i)_{J(1,0)}|_0 = \Delta^2 (\sum_{u \in W_0} \frac{d}{u^3} (J' \cdot \bar{\partial} \eta_i)|_u + o_i)$$
(A.56)

where $|o_{u,i}| \leq c(\Delta + d^2) \sum_{u \in W_0} |J'|_u|$ with c again depending only on $W_0$, but neither on d nor on $\Delta$ if $\Delta \ll d$. With (A.56) in hand, the argument from the last paragraph in Part 6 of the proof of Proposition 3.7 can be applied here to see that $\{J'|_u\}|_{u \in W_0}$ can be chosen when d is very small and $\Delta$ is much less than d so as to make the n-tuple $(\Pi(\partial \lambda_i)_{J(1,0)}|_0)_{i=1,\ldots,n}$ equal to any given n-tuple. This last fact implies the asserted surjectivity at $((\varphi, J, \mathfrak{j}), \mathfrak{p}, p)$ of the differential of $\mathfrak{f}_\alpha$.

*Part 4*: To complete the proof of Proposition 3.8, consider now the open set of points $((\varphi, J, \mathfrak{j}), \mathfrak{p}) \in \mathcal{P}^m \times (\times_{d-2} \Sigma)$ where J is a regular value for the projections from $\mathcal{P}^m$ and from each $n \in \{3, 4\}$ and each $\alpha \in \{1, \ldots, d-2\}$ version of the submanifold $\mathfrak{f}_\alpha^{-1}(0) \subset \mathcal{Q}_n$. Require also that $\mathfrak{p}$ have distinct entries. Use $\mathcal{Q}^2$ to denote the fiber bundle over this open set in $\mathcal{P}^m \times (\times_{d-2} \Sigma)$ whose fiber over a given $((\varphi, J, \mathfrak{j}), \mathfrak{p})$ consists of the product of two copies of the fiber of $\mathcal{Q}_3$. For indices $\alpha \neq \beta$ from $\{1, \ldots, d-2\}$, set $\mathfrak{f}_{\alpha\beta}$ to be the fiber preserving map from $\mathcal{Q}^2$ to the vector bundle $(\oplus_3 \mathfrak{N}_\alpha) \oplus (\oplus_3 \mathfrak{N}_\beta)$ that sends the fiber point $(p, p')$ over the given point $((\varphi, J, \mathfrak{j}), \mathfrak{p})$ to $(\mathfrak{f}_\alpha(p), \mathfrak{f}_\beta(p'))$. The bundle $\mathcal{Q}^2$ has a fiber preserving, free action of $SO(3) \times SO(3)$, and this action restricts to an action on $\mathfrak{f}_{\alpha\beta}^{-1}(0)$.

Since the entries of $\mathfrak{p}$ are distinct, the Green's function constructions given in Part 3 above can be applied near the index $\alpha$ entry and the index $\beta$ entry of $\mathfrak{p}$ simultaneously to prove that the differential of $\mathfrak{f}_{\alpha\beta}$ is transverse along its zero locus. This being the case, it follows that $\mathfrak{f}_{\alpha\beta}^{-1}(0)$ is a manifold with a free action of $SO(3) \times SO(3)$. Let $\mathfrak{Q}_{\alpha\beta}$ denote its quotient by this action. The latter is a manifold also. Moreover, the associated map to $\mathcal{J}^m$ has everywhere Fredholm differential, with index equal to 4d - 10. The arguments used to prove Proposition A.5 can be used here to prove that there is a residual set of regular values for this map in $\mathcal{J}$. Let J denote such a regular value. Then $\mathfrak{Q}_{\alpha\beta}|_J$ is a smooth manifold of dimension 4d-10. The manifold $\mathfrak{Q}_{\alpha\beta}|_J$ maps to $\times_{d-2} X$ via the composition of the map to $\mathcal{P}|_J \times (\times_{d-2} \Sigma)$ with the map $\vartheta_{d-2}$. The inverse image of a regular value of the latter map is empty since the domain space has smaller dimension than the range.

Suppose that $J \in \mathcal{J}$ is a regular value for all $\alpha \neq \beta$ versions of $\mathfrak{f}_{\alpha\beta}$. Suppose in addition that $\mathfrak{w} \in \times_{d-2} X$ is a regular point and also a regular value for all of the associated maps from $\{\mathfrak{Q}_{\alpha\beta}|_J\}_{1\le\alpha<\beta\le d-2}$. Then the conclusions of the preceding paragraph imply the following: If $w \ne w'$ are entries of $\mathfrak{w}$, then $\mathcal{Y}^w$ is disjoint from $\mathcal{Y}^{w'}$.